\documentclass[10pt]{article}

\usepackage{latexsym, pictexwd}
\usepackage[a4paper]{geometry}
\usepackage[english]{babel} 
\usepackage{color}
\usepackage{amsmath, amsfonts,amssymb,mathtools,amsthm,mathrsfs}
\usepackage{bbm}
\usepackage{enumerate}
\usepackage{hyperref}

\newtheorem{proposition}{Proposition}[section]
\newtheorem{theorem}[proposition]{Theorem}
\newtheorem{lemma}[proposition]{Lemma}
\newtheorem{corollary}[proposition]{Corollary}
\newtheorem{definition}[proposition]{Definition}
\newtheorem{assumption}[proposition]{Assumption}

\newtheorem{setting}[proposition]{Setting}
\theoremstyle{definition}
\newtheorem{remark}[proposition]{Remark}
\newtheorem{example}[proposition]{Example}
\newtheoremstyle{step}{3pt}{0pt}{}{}{\bf}{}{.5em}{}
\theoremstyle{step} 

\DeclareMathAlphabet{\mathpzc}{OT1}{pzc}{m}{it}

\numberwithin{equation}{section}

\providecommand{\N}{{\ensuremath{\mathbbm{N}}}}
\providecommand{\Q}{{\ensuremath{\mathbbm{Q}}}}
\providecommand{\R}{{\ensuremath{\mathbbm{R}}}}

\newcommand{\E}{{\ensuremath{\mathbb{E}}}}
\renewcommand{\P}{{\ensuremath{\mathbb{P}}}}
\providecommand{\1}{{\ensuremath{\mathbbm{1}}}}
\providecommand{\al}      {{\ensuremath{\alpha}}}

\providecommand{\Var}{\operatorname{Var}}

\providecommand{\del}{{\ensuremath{\partial}}}

\providecommand{\ra}{\rightarrow}
\providecommand{\C}{{\ensuremath{\mathcal{C}}}}
\providecommand{\D}{{\ensuremath{\mathcal{D}}}}

\newcommand{\cadlag}{c\`adl\`ag}

\providecommand{\eps}     {{\ensuremath{\varepsilon}}}
\newcommand{\Dom}{{\ensuremath{\mathrm{Dom}}}}
\providecommand{\Law}[2][]{{\ensuremath{\mathcal L^{#1}\left(#2\right)}}}
\providecommand{\ceil}[1]{{\ensuremath{\lceil {#1} \rceil}}}
\providecommand{\floor}[1]{{\ensuremath{\left\lfloor {#1} \right\rfloor}}}
\providecommand{\sn}[1][0]    {{\ensuremath{\sigma^n_{#1}}}}
\providecommand{\tn}[1][0]    {{\ensuremath{\tau^n_{#1}}}}

\providecommand{\snh}[1][0]    {{\ensuremath{\hat \sigma^n_{#1}}}}
\providecommand{\tnh}[1][0]    {{\ensuremath{\hat \tau^n_{#1}}}}
\providecommand{\snb}[1][0]    {{\ensuremath{\bar \sigma^n_{#1}}}}
\providecommand{\tnb}[1][0]    {{\ensuremath{\bar \tau^n_{#1}}}}

\providecommand{\Xnb}[1][0]    {{\ensuremath{\bar X^n_{#1}}}}

\providecommand{\T}{{\ensuremath{\mathcal{T}}}}
\providecommand{\F}{{\ensuremath{\mathcal{F}}}}

\providecommand{\limn}   {{\ensuremath{{\displaystyle \lim_{n \ra \infty}}}}}
\providecommand{\limsupn}   {{\ensuremath{{\displaystyle \limsup_{n \ra \infty}}}}}

\begin{document}
\thispagestyle{empty}

\title{\Large Stochastic averaging for multiscale Markov processes
	with an\\
  application to a Wright-Fisher model with fluctuating selection}
\author{\sc by Martin Hutzenthaler, Peter Pfaffelhuber and Clemens Printz}

\date{\today}

\maketitle
\begin{abstract}  
	Let $Z = (Z_t)_{t\in[0,\infty)}$ be an ergodic Markov process and,
	for every $n\in\N$, let
	$Z^n = (Z_{n t})_{t\in[0,\infty)}$ drive a process~$X^n$.
	Classical results show
	under suitable conditions that the sequence
	of non-Markovian processes $(X^n)_{n\in\N}$
	converges to a Markov process and give its infinitesimal characteristics.
	Here, we consider a general sequence
	$(Z^n)_{n\in\N}$.
	Using a general
	result on stochastic averaging from [Kur92], we derive
	conditions 
	which ensure that the sequence $(X^n)_{n\in\N}$ converges as in the classical
	case. As an application, we consider the diffusion limit of
	a Wright-Fisher model with fluctuating selection.
\end{abstract}
\makeatletter
\let\@makefnmark\relax
\let\@thefnmark\relax

\@footnotetext{\emph{AMS 2010 subject classification:} 60F05 (Primary)
	60K37, 60J80 (Secondary)}

\@footnotetext{\emph{Key words and phrases:} Stochastic averaging,
	random walk in random environment, martingale problem} \makeatother


\section{Introduction}

Stochastic averaging is a well-known concept and
has been
introduced a while ago (see e.g.\ \cite{Kash1966}).
Consider a sequence of bivariate Markov processes $(X^n,Z^n)_{n\in\N}$.
The general idea is that
the processes $(Z^n)_{n\in\N}$ (subsequently denoted as fast variables)
converge quickly to an equilibrium and
that
the non-Markovian processes $(X^n)_{n\in\N}$ (subsequently denoted as slow variables)
evolve on a slower timescale
and
only sense this equilibrium
in the limit as $n\to\infty$
and, thus, converge to a Markov process.
Stochastic averaging results in the literature include, e.g., Thm.~1.7.6 in~\cite{EthierKurtz1986}
and the references
\cite{EthierNagylaki1980,
	EthierNagylaki1988,AV2010, PardouxVere1, PardouxVere2,VereKulik2012}.
Theorem 2.13
in the recent paper~\cite{KangKurtzPopovic2014} also treats processes
with three different timescales under different assumptions.
All of these references assume that the fast variables converge in a suitable
sense to an equilibrium process or to an equilibrium distribution 
(depending on the current state of the slow variables).
~

This paper is motivated by the observation that in many applications
the fast variables do not converge to an equilibrium process or an equilibrium distribution.
Still, the slow variables can be approximated by a Markov process.
Our intuition is that the slow variables only depend on the fast variables
through certain functions and for the processes to converge it suffices
that these functions of the fast variables converge suitably.
We consider three timescales since there are often three types of dynamics involved, namely
dynamics depending on and affecting the slow variables only,
dynamics by which the fast variables affect the slow variables and 
dynamics depending on and affecting the fast variables only.
More precisely, we assume for every $n\in\N$ 
that the pre-generator $L_n$ of the Markov process $(X^n,Z^n)$
satisfies
for all $f\in\Dom(L_n)$ that
\begin{align}\label{eq:genform}
L_n f & = L_{0,n} f + n\cdot L_{1,n} f +
n^2\cdot L_{2,n} f,
\end{align}
where $\Dom(L_n)$ is the domain of the pre-generator $L_n$.
For every $n\in\N$, we think of $n^2L_{2,n}$ as the pre-generator of the fast variable $Z^n$
evolving on timescale $O(n^2)$
and we think of $L_{0,n}+nL_{1,n}$ as the pre-generator of the slow variable $X^n$
given the fast variable $Z^n$.
We will show in our main result,
Theorem~\ref{T:stochastic.averaging}  below,
under suitable assumptions that the non-Markov processes $(X^n)_{n\in\N}$ converge to
a Markov process. 
Theorem~\ref{T:stochastic.averaging}  below
is an application of Theorem 2.1 in~\cite{Kurtz1992} which
is a general result on  stochastic averaging.
The main contribution of our paper is to demonstrate
how to apply the abstract result of~\cite{Kurtz1992}
to settings
where the occupation measures of the driving processes $(Z^n)_{n\in\N}$
might not converge.
In particular Theorem~\ref{T:stochastic.averaging} enables us to derive
the diffusion approximation of Wright-Fisher processes with fluctuating selection 
which is an important model in population genetics.

We explain our approach with a simple example. Let a random walker on the real line
move at constant speed ($\in\R$ indicating positive or negative direction)
for an exponentially distributed  time period, choose then a new speed according to a given distribution
and continue so forth. If the exponential waiting times become shorter and shorter and the
distributions of the random speeds are suitable then these processes
converge to a Brownian motion.
More formally,
let $N$ be a Poisson process with rate~1 and, for every $n\in\N$,
let $\bar{Z}_1^n, \bar{Z}_2^n,...$ be independent and identically distributed real-valued
random
variables with distribution $\pi_n$ having mean $\mu_n\in\R$ and
variance $\sigma_n^2/2\in[0,\infty)$.
We assume that $\lim_{n\to\infty}n\mu_n=a\in\R$ and that
$\lim_{n\to\infty}\sigma_n^2=\sigma^2\in(0,\infty)$.
For every $n\in\N$ define $Z^n = (Z_t^n)_{t\in[0,\infty)}$ by
$Z_t^n = \bar{Z}^n_{N_{n^2t}}$
and define $X^n =
(X_t^n)_{t\in[0,\infty)}$ for every $t\in[0,\infty)$ by
\begin{align}
\label{eq:exXn}
X_t^n := n \int_0^{t} Z_s^n \,ds.
\end{align}
For each $n\in\N$ the pre-generator of the bivariate Markov process
$(X^n,Z^n)$ satisfies
for all $f\in \C_c^\infty(\R^2,\R)$ that
$L_n f = n \cdot L_{1,n}f + n^2 \cdot L_{2,n} f$
where
\begin{align}\label{eq:321}
(L_{1,n}f)(x,z) = z \frac{\partial f}{\partial x}(x,z), \qquad
(L_{2,n}f)(x,z) = \int_{\R} f(x,y) \pi_n(dy) - f(x,z)
\end{align}
for all $(x,z)\in\R^2$.
Of course 
a corollary of the celebrated Lindeberg-Feller theorem
shows that the finite-dimensional distributions of
$(X^n)_{n\in\N}$ converge
to a Brownian motion
if and only if
Lindeberg's condition is satisfied
or, equivalently,
if for all $\eps\in(0,1)$ it holds that
$\lim_{n\to\infty}\E\left[(Z_0^n)^2\1_{\{|Z_0^n|>\eps n\}}\right]=0$.
Using our stochastic averaging result,
Theorem~\ref{T:stochastic.averaging} below, we will obtain convergence in distribution on the space of cadlag 
functions and 
we will assume that there exists $\delta\in(0,1]$ such that
$\sup_{n\in\N}\E\left[|Z_0^n|^{2+\delta}\right]<\infty$.
The following heuristic then explains with a pre-generator
calculation
why the only possible
limit process of the sequence $(X^n)_{n\in\N}$
is $(at+\sigma W_t)_{t\in[0,\infty)}$ where $W$ is a real-valued standard Brownian motion.
Since $\sup_{s\in[0,\infty)}\sup_{n\in\N}\E[(Z_s^n)^2]<\infty$,
the occupation measures (see section~\ref{sect:notations} below)
of the processes $(Z^n)_{n\in\N}$
are relatively compact.
Moreover, the law of large numbers together with $\sup_{n\in\N}\E\big[|Z_0^n|^{2+\delta}\big]<\infty$
implies for alle $t\in(0,\infty)$ that $\lim_{n\to\infty}\E\big[|\int_0^t(Z_s^n)^2-\E[(Z_s^n)^2]\,ds|\big]=0$,
and then Lemma~\ref{l:occupation.measure} below implies almost surely that
$\int_0^\infty z^2\Gamma(ds,dz)=ds\,\lim_{n\to\infty}\E[(Z_0^n)^2]=ds\,\tfrac{\sigma^2}{2}$.
Consequently,
using for every $f\in \C_c^2(\R,\R)$ and $n\in\N$ that $L_{2,n}f\equiv 0$,
we get for all $f\in \C_c^2(\R,\R)$ approximately in the limit $n\to\infty$
that
\begin{align}
\begin{split}
&[0,\infty)\ni t\mapsto
(f+\tfrac{1}{n}L_{1,n}f)(X_t^n,Z_t^n)-\int_0^t\left(L_n(f+\tfrac{1}{n}L_{1,n}f)\right)(X_s^n,Z_s^n)\,ds
\\&
=f(X_t^n)+\tfrac{1}{n}Z_t^nf'(X_t^n)-\int_0^t \left(nL_{1,n}f+L_{1,n}L_{1,n}f+nL_{2,n}L_{1,n}f\right)(X_s^n,Z_s^n)\,ds
\\&
=f(X_t^n)+\tfrac{1}{n}Z_t^nf'(X_t^n)-\int_0^t \left(L_{1,n}L_{1,n}f\right)(X_s^n,Z_s^n)+n\int_{\R}(L_{1,n}f)(X_s^n,y)\,\pi_n(dy)\,ds
\\&
=f(X_t^n)+\tfrac{1}{n}Z_t^nf'(X_t^n)-\int_0^t f''(X_s^n)(Z_s^n)^2\,ds-\int_0^tf'(X_s^n)n\mu_n\,ds
\\&
\approx f(X_t^n)-\int_0^t f''(X_s^n)\int_0^\infty z^2 \Gamma(ds,dz)-\int_0^tf'(X_s^n)n\mu_n\,ds
\\&
\approx f(X_t^n)-\int_0^t f''(X_s^n)\tfrac{\sigma^2}{2}\,ds-\int_0^t f'(X_s^n)a\,ds
\end{split}
\end{align}
is a local martingale.
So we recognize the pre-generator of the Brownian motion
$(at+\sigma W_t)_{t\in[0,\infty)}$.
Note that the second derivative appears
as $\tfrac{\sigma^2}{2}f''(x)=\lim_{n\to\infty}\E\left[\left(L_{1,n}L_{1,n}f\right)(x,Z_0^n)\right]$ where $x\in\R$
and $f\in \C_c^2(\R,\R)$.
An analogous iterated operator appears also in the Wright-Fisher model with fluctuating selection;
see Remark~\ref{r:informal.deduction.SDE.KLM} for more details.
Moreover, we emphasize that the strength of Theorem~\ref{T:stochastic.averaging} is that
 the sequence $(\pi_n)_{n\in\N}$ does not need to have any convergence properties
except for suitable convergence of the first and second moments.
Stochastic averaging results in the literature typically assume that $(\pi_n)_{n\in\N}$ converges suitably to a measure; e.g. Theorem~1.7.6 in~\cite{EthierKurtz1986} assumes that the limit of $(\pi_n)_{n\in\N}$ exists and that $\lim_{n\to\infty}L_{2,n}$
generates a strongly continuous contraction semigroup $(S_t)_{t\in[0,\infty)}$ and that the limit $\lim_{\lambda\to 0+}y \int_{0}^{\infty}e^{-\lambda t}S_t\,dt$ exists in a weak sense.

Next we explain our approach in the abstract setting of the second paragraph of this introduction.
For simplicity,
we assume for every $n\in\N$ that $L_{0,n}=0$.
For this, fix
a
function $f$ in a dense subset of the continuous and bounded functions on the state
space of the limiting Markov process.
We assume 
for every $n\in\N$
-- identifying $f$ with a function in the domain of $L_n$ which is constant in the second argument --
that 
there exists a function $h_n\in\Dom(L_n)$
and
a measure $\pi_n$ on the state space 
of $Z^n$ (typically the ergodic equilibrium of $Z^n$)
such that for all $(x,z)$ in the state space of $(X^n,Z^n)$ and all $n\in\N$ it holds that
\begin{equation}  \begin{split}
(L_{2,n}f)(x,z)&=0\\
L_{2,n}h_n(x,z) & = \int L_{1,n} f(x, .) d\pi_n - L_{1,n} f(x, z).
\end{split}     \end{equation}
Then for all $n\in\N$ it follows from $L_n$ being the pre-generator of $(X^n,Z^n)$ that
\begin{equation} \begin{split}
&[0,\infty)\ni t
\mapsto
(f+\tfrac{1}{n}h_n)(X_t^n,Z_t^n)-\int_0^t\left(L_n(f+\tfrac{1}{n}h_n)\right)(X_s^n,Z_s^n)\,ds
\\& = (f+\tfrac{1}{n}h_n)(X^n_t, Z^n_t)
- \int_0^t 
n(L_{1,n}f) (X^n_s, Z^n_s)
+
(L_{1,n}h_n) (X^n_s, Z^n_s)
+
n(L_{2,n}h_n) (X^n_s, Z^n_s)
\,ds
\\& = (f+\tfrac{1}{n}h_n)(X^n_t, Z^n_t) - \int_0^t
n \int (L_{1,n} f)(X^n_s, \cdot) d\pi_n
+
(L_{1,n}h_n) (X^n_s, Z^n_s) 
\,ds
\label{eq:intro.local.martingale}
\end{split}
\end{equation}
is a local martingale. Moreover we assume that the sequence $\tfrac{1}n h_n$ converges suitably to $0$.
Now, as in our application in Section~\ref{S:KLM},
the sequence of functions $(n\cdot L_{1,n}f)_{n\in\N}$ might not converge
but the sequence of averaged functions does. So we additionally
assume for every $x$ in the state space of the limiting Markov process that
the limit
\begin{align}
\label{eq:prov2} 
\lim_{n\to\infty} n \int (L_{1,n} f)(x, \cdot) d\pi_n &
=: (A_1 f)(x)
\end{align}
exists.
Moreover we assume 
for every $n\in\N$
that there exists
a function $g_n$ on the state space of $Z^n$
and a suitable function $A_2f$
such that
for all $t\in[0,\infty)$ in the limit $n\to\infty$ it holds that
\begin{equation}  \begin{split}
\label{eq:prov3} \int_0^t L_{1,n}h_n(X^n_s, Z^n_s) ds & \approx
\int_0^t A_2 f(X^n_s,g_n(Z_s^n)) ds 
= \int_0^t \int A_2 f(X^n_s,z))
\Gamma_{g_n(Z^n)}(ds,dz),
\end{split}     \end{equation}
where we used for each $n\in\N$
the occupation measure $\Gamma_{g_n(Z^n)}$ of $g_n(Z^n)$.
The reason for introducing the functions $(g_n)_{n\in\N}$ is that
the occupation measures of the processes $(Z^n)_{n\in\N}$ might not converge
but the occupation measures of $(g_n(Z^n))_{n\in\N}$ (which possibly have a much
smaller state space) could converge.
Finally we assume 
that the sequence $(X^n)_{n\in\N}$ satisfies the compact containment condition
and 
for every $t\in[0,\infty)$
that the family $\{g_n(Z_s^n)\colon n\in\N,s\in[0,t]\}$ is tight.
Then the sequence $(X^n,\Gamma_{g_n(Z^n)})_{n\in\N}$ is tight
and~\eqref{eq:intro.local.martingale}, \eqref{eq:prov2} and~\eqref{eq:prov3}
suggest that every limit point $(X,\Gamma)$
satisfies that
\begin{align}
[0,\infty)\ni t
\mapsto
f(X_t) - \int_0^t A_1f(X_s)ds - \int_0^t \int A_2 f(X_s,z))
\Gamma(ds,dz)
\end{align}
is a local martingale, suggesting the form of the pre-generator for
$X$. Before we carry out the technicalities of this reasoning, we will fix some notations including the occupation measure of a stochastic process.

\subsection{Notation} \label{sect:notations}
	Throughout this section let $(E, d_E)$ be a metric space.
	\begin{enumerate}
	\item 	We write $\N=\{1,2,3,\dots\}$ and $\N_0=\N\cup\{0\}$.
	\item For all $x,y\in\R$ we write $x\vee y:=\operatorname{max}\{x,y\}$, $x\wedge y:=\operatorname{min}\{x,y\}$, $x^+:=x\vee0$  and $x^-:=-(x\wedge0)$.
	\item	We use the convention that $0^0:=1$, that $\inf \emptyset:=\infty$ and that $\sup \emptyset:=-\infty$.	
	\item For every countable set $\mathcal D$, $z\in\R^{\mathcal D}$ and $\beta\in\N_0^{\mathcal D}$ with $\#\{i\in\mathcal D\colon\beta_i\neq0  \}<\infty$ let $z^\beta:=\prod_{i\in\mathcal D}z_i^{\beta_i}$.
	\item For every $t\in[0,\infty)$ and every set $\T\subseteq [0,\infty)$ we write
	$ \floor{t}_{\T} :=\sup\{s\in\T\cup\{0\}\colon s \leq t \}$,  $\lceil t \rceil_{\T} :=\inf \{s\in\T\colon s\geq t\}$ and write the floor function as $ \floor{t} := \floor{t}_{\N_0}$.
	\item For a random variable X and a probability measure $\mu$ we use the notation $X\sim\mu$ to denote that $X$ is distributed according to $\mu$.
	For $x\in E$ we denote by $\delta_x$ the Dirac-probability-measure.
	For $n\in\N$, $p\in[0,1]$ we denote by $\operatorname{Bin}(n,p)$ the binomial distribution, in particular $\operatorname{Bin}(n,0)=\delta_0$ and $\operatorname{Bin}(n,1)=\delta_n$.
	For $\lambda\in(0,\infty)$ we denote by $\operatorname{Poi}(\lambda)$ the Poisson distribution with parameter $\lambda$
.
	
	\item We denote by
			$\mathcal B(E,\R)$ (resp. $\mathcal C(E,\R)$/$\mathcal C_b(E,\R)$/$\mathcal{C}_c(E,\R)$) 
			the set of
			Borel-measurable 
			(resp.	continuous/bounded and continuous/compactly supported and continuous)
			functions $f\colon E\to \mathbb R$ and we denote by
			$\D([0,\infty),E)$ the set of \cadlag-functions $f\colon
			[0,\infty) \to E$.
	\item For a function $A\colon\Dom(A)\subseteq \mathcal C(E,\R)
			\to \mathcal B(E,\R)$, we say that an $E$-valued stochastic process $X =
			(X_t)_{t\in[0,\infty)}$ solves the (local)
			$\mathcal{D}([0,\infty),E)$-martingale problem for
			$A$ with respect to a set $\T\subseteq[0,\infty)$ if $X$ has \cadlag-paths,
      for all $u\in[0,\infty)$ it holds that $\int_{0}^{\floor{u}_{\T}}|(Af)(X_s)|ds<\infty$ almost surely and
			\begin{equation} \begin{split}
			\Big( f(X_t) - f(X_0) - \int_0^{\floor{t}_{\T}} (Af)(X_s)ds\Big)_{t\in[0,\infty)}
			\end{split}     \end{equation}
			is a (local) martingale
			for all $f\in\Dom(A)$. In this case, we say that $A$ is a
			pre-generator for the process $X$.
      Finally, we say that $(X_t)_{t\in[0,\infty)}$ solves the
			(local) $\mathcal{D}([0,\infty),E)$-martingale problem if
       $(X_t)_{t\in[0,\infty)}$ solves the
			(local) $\mathcal{D}([0,\infty),E)$-martingale problem with respect to the set $[0,\infty)$.
	\item A sequence $(X_t^n)_{t\in[0,\infty)}$, $n\in\N$,
			of $E$-valued stochastic processes
			is said to satisfy the \emph{compact containment condition}, if for every
			$\varepsilon, t>0$, there exists a compact set $K\subseteq E$ with
			$$ \inf_{n\in\N} \P(X_s^n \in K \text{ for all }s\in[0,t]) > 1-\varepsilon.$$
	\item
			Let $\mathcal B(E)$ be
			the Borel $\sigma$-algebra and $\mathcal{M}(E)$ be the set of Borel-measures on $(E,\mathcal{B}(E))$. Let $\mathcal{M}_1(E)\subset \mathcal{M}_f(E) \subset\mathcal{M}(E)$ be the subsets of probability measures and of finite Borel-measures, both endowed with the weak topology (denote weak convergence by $\Rightarrow$).
	\item Denote
			the set of occupation measures by
			\begin{equation}
			\mathcal{L}_{m}(E)
			:=
			\left\{\Gamma\in\mathcal{M}([0,\infty)\times E)\colon\Gamma([0,t]\times E)=t
			\text{ for all }t\in[0,\infty)
			\right\}.
			\end{equation}
		  For every $t\in[0,\infty)$ let $\rho_t$ be the Prokhorov metric on $[0,t]\times E$. We endow $ \mathcal{L}_{m}(E)$ with the metric $\rho$ that satisfies for all $\mu,\nu \in \mathcal{L}_{m}(E)$ that $\rho(\mu,\nu)=\int_0^\infty e^{-t}(1\wedge \rho_t(\mu|_{\mathcal{B}([0,t]\times E)}, \nu|_{\mathcal{B}([0,t]\times E)}))dt$.
	\item For an $E$-valued stochastic process
			$X=(X_t)_{t\in[0,\infty)}$ with \cadlag-paths, its
			\emph{occupation measure} is the unique $\mathcal L_m(E)$-valued
			random variable $\Gamma_{X}$ such that for all $t\in[0,\infty)$
			and all $B\in \mathcal{B}(E)$ it holds that
			$$ \Gamma_{X}([0,t]\times B) = \int_0^t \1_B(X_s)ds\,.$$
  \item
      We adopt the convention that
      zero times an undefined quantity is zero.
      Thereby e.g. the expression $0\cdot f(-1)$ is defined for every
      function $f\colon[0,\infty)\to\R$.
\end{enumerate}

\section{Main result}
We now describe the setting we are working in as well as some basic
assumptions for our main result.

\begin{assumption}\label{ass:main}
  \begin{enumerate}
  \item Let $(\theta_n)_{n\in\N}\subset(0,\infty)$ be a sequence of real numbers with
    $\theta_n\xrightarrow{n\to\infty} \infty$ and let $(\T_n)_{n\in\N} \subseteq [0,\infty)^\N$ 
    be a sequence of subsets of $[0,\infty)$ such that
    $\sup_{s\in[0,\infty)}  \left(s-\floor{s}_{\T_n}\right)
    \xrightarrow{n\to\infty} 0$.
  \item Let $(S, d_S)$ and $(E, d_E)$ be complete and separable metric
    spaces, and
    for all $n\in\N$ let $S_n$ and $E_n$ be Borel measurable sets and let  $\tilde S_n$ and $\tilde E_n$ be sets such that $S_n\subseteq \tilde S_n \subseteq S$ and
    $E_n\subseteq \tilde E_n\subseteq E$.
    \\
    For every $n\in\N$,
    let
    $L_n\colon \Dom(L_{n})\subseteq \C(S_n\times E_n,\R) \to \mathcal B(S_n\times E_n,\R)$
    be a linear function,
    let $D_n\subseteq\{f\colon \tilde S_n\times \tilde E_n \to \R\colon f|_{S_n\times E_n}\in\Dom(L_n)\}$
    and let 
    $L_{0,n},L_{1,n},L_{2,n}\colon D_n \to \mathcal B(S_n\times E_n,\R)$
    be functions
    such that for all $f\in D_n$ it holds that
    \begin{equation}  \begin{split} \label{eq:form.Ln}
    L_n(f|_{S_n\times E_n})=\left(L_{0,n}+\theta_n L_{1,n}+\theta_n^2L_{2,n}\right)(f).
    \end{split}     \end{equation}
    
  \item For every $n\in\N$,
    let $
    (X^n, Z^n)=
    (X_t^n,Z_t^n)_{t\in[0,\infty)}$ be a solution
    of the $\D([0,\infty), S_n\times E_n)$-martingale problem for $L_n$  
    with respect to the set $\T_n$.

  \item The sequence $(X^n)_{n\in\N}$ of $S$-valued stochastic processes satisfies the
    compact containment condition.
  \item Let $(H,d_H)$ be a complete and separable metric space
    and let $g_n\colon E_n\to H$, $n\in\N$, be Borel measurable functions
    such that the family
    $(\Gamma_{g_n(Z^n)})_{n\in\N}$ is tight
    in $\mathcal{L}_m(H)$.
  \end{enumerate}
\end{assumption}

\begin{remark} \label{r:sufficient.for.occupation.time}
  Lemma~1.3 in \cite{Kurtz1992} and Prokhorov's theorem imply that
    Assumption~\ref{ass:main}.5 is fulfilled if for every $t\in(0,\infty)$
    the family
    $\{g_n(Z_s^n)\colon n\in\N, s\in[0,t]\}$ is tight.
\end{remark}

The following theorem, Theorem~\ref{T:stochastic.averaging}, is our main result on stochastic averaging.
\begin{theorem}[Stochastic averaging for solutions of martingale problems]
  \label{T:stochastic.averaging} 
  Let the setting from
  Assumption~\ref{ass:main} be given,
  let $D_0 \subseteq \C_b(S,\R)$ be a
  dense set in the topology of uniform convergence on compact sets and
  let ${A}_1\colon D_0 \to \C_b(S,\R)$ and ${A}_2\colon D_0 \to
  \C(S\times H,\R)$ be functions.
  Suppose for every $f\in D_0$ that there exist $f_n,h_n \in
  D_n$, $n\in\N$,
  such that for all $n\in\N$ it holds that
  $L_{2,n}f_n=0$,
  such that for all  $t\in[0,\infty)$ it holds that
  \begin{align}
    \label{eq:fN.converges.to.f}
    \lim_{n\to\infty} \E\bigg[\sup_{s\in[0,t]}
    \left|f(X_s^n)-(f_n+\tfrac{1}{\theta_n}h_n)(X_s^n,Z_s^n)\right|\bigg]
    &=0\,,
  \end{align}
  such that for all $t\in[0,\infty)$ there exists $p\in(1,\infty)$ with
  \begin{equation} \begin{split} \label{eq:Lp.assumption}
      \sup_{n\in\N}\int_0^t\E\big[
        \left|\left({A}_2f\right)\left(X_s^n,g_n\left(Z_s^n\right)\right)\right|^p\big]\,ds<\infty\,,
    \end{split}
  \end{equation} 
  such that all integrals in \eqref{eq:generator.A1} and \eqref{eq:generator.A2} are well-defined and such that for all $t\in[0,\infty)$ it holds that
  \begin{align}
    \lim_{n\to\infty} \E\bigg[ \sup_{s\in[0,t]} \bigg| \int_0^{s}
    \left[
    \left({A}_1 f\right)(X_r^n) -
    \left( \theta_nL_{2,n}h_n+\theta_nL_{1,n}f_n+L_{0,n}f_n\right)(X_r^n,Z_r^n)
        \right]
     \,dr \bigg|\bigg]
    &=0\,,
    \label{eq:generator.A1}
    \\
    \lim_{n\to\infty} \E\bigg[ \sup_{s\in[0,t]}\bigg| \int_0^{s}
    \left[
    \left({A}_2 f\right)(X_r^n,g_n(Z_{r}^n))
    -\left(L_{1,n}h_n+\tfrac{1}{\theta_n}L_{0,n}h_n\right)(X_r^n,Z_r^n)
    \right]
    \,dr\bigg|\bigg] &=0\,.
    \label{eq:generator.A2}
  \end{align}
  Then 
  $(X^n, \Gamma_{g_n(Z^n)})_{n\in\N}$ is relatively compact in $\mathcal D ([0,\infty),S)\times \mathcal L_m(E)$ and for every limit
  point $((X_t)_{t\in[0,\infty)}, \Gamma)$ 
  and
  for every
  $f\in D_0$ it holds that all integrals in~\eqref{eq:resulting.martingale.main.result} are well-defined and
  \begin{equation} \label{eq:resulting.martingale.main.result}
    \left(f(X_t)-\int_0^t ({A}_1 f)(X_s)\,ds
      -\int_0^t \int_{H} ({A}_2 f)(X_s,y)\Gamma(ds, dy)\right)_{t\in[0,\infty)}
  \end{equation}
  is a martingale.
\end{theorem}

\begin{proof}[Proof of Theorem~\ref{T:stochastic.averaging}]
	We will apply Theorem 2.1 in~\cite{Kurtz1992} to the sequence
	$((X^n, g_n( Z^n)))_{n\in\N}$ and first check the assumptions.
	By
	Assumption~\ref{ass:main}.4, the sequence $(X^n)_{n\in\N}$ satisfies the compact
	containment condition. Note also that the proof of Theorem 2.1
	in~\cite{Kurtz1992} only requires the relative compactness of
	$\Gamma_{g_n(Z^n)}$ and that the stronger assumption of
	relative compactness of the family $\{g_n(Z_t^n)\colon n\in\N, t\in[0,\infty)\}$ is obsolete.
	
	Next, fix $f\in D_0$ for the rest of the proof.
	By assumption, there exist
	$f_n,h_n \in
	D_n$, $n\in\N$
	such that $L_{2,n}f_n=0$ for all $n\in\N$
	and such that~\eqref{eq:fN.converges.to.f} 
	holds.
	For all $n\in\N$
	$\Dom(L_n)$ is a vector space so that
	$(f_n+\tfrac{1}{\theta_n}h_n)|_{S_n\times E_n}\in \Dom(L_n)$.
	Define
	\begin{align}
	\eps_t^n & :=(f_n+\tfrac{1}{\theta_n}h_n)(X_t^n,Z_t^n)-f(X_t^n)
	+ \int_{\floor{t}_{\T_n}}^t
	\Big[({A}_1f)(X_s^n)+({A}_2f)(X_s^n,g_n(Z_s^n))\Big]ds 
	\\ &  \qquad + \int_0^{\floor{t}_{\T_n}}
	\Big[
	({A}_1f)(X_s^n)+({A}_2f)(X_s^n,g_n(Z_s^n))
	-\Big(L_n((f_n+\tfrac{1}{\theta_n}h_n)|_{S_n\times E_n})\Big)(X_s^n,Z_s^n)\Big]\,ds
	\nonumber
	\end{align}
	for all $t\in[0,\infty)$ and all $n\in\N$.
	Then Assumption~\ref{ass:main}.3 implies
	for every $n\in\N$ that the process
	\begin{equation} \begin{split} 
	&\Big(f(X_t^n)  -\int_0^t
	\Big[({A}_1f)(X_s^n)+({A}_2f)(X_s^n,g_n(Z_s^n))\Big]ds+\eps_t^n\Big)_{t\in[0,\infty)} \\ =
	&\Big(\left(f_n+\tfrac{1}{\theta_n}h_n\right)(X_t^n,Z_t^n) -\int_0^{\floor{t}_{\T_n}}
	\left(L_n\big((f_n+\tfrac{1}{\theta_n}h_n)|_{S_n\times E_n}\big)\right)(X_s^n,Z_s^n)\,ds\Big)_{t\in[0,\infty)}
	\end{split} \end{equation} 
	is a martingale. 
	By
	Assumption~\eqref{eq:Lp.assumption} and global boundedness of
	${A}_1f$, for every $t\in[0,\infty)$ there exists a real number
	$p\in(1,\infty)$ such that
	\begin{equation} \begin{split} \label{eq:p-th.momentA1+A2}
	\sup_{n\in\N}\int_0^t\E\left[\left|
	({A}_1f)(X_s^n)+({A}_2f)(X_s^n,g_n(Z_s^n))\right|^p\right]\,ds<\infty.
	\end{split} \end{equation}
	H\"older's inequality,  Jensen's inequality, Fubini,
	Assumption~\ref{ass:main}.1  and \eqref{eq:p-th.momentA1+A2} imply for every $t\in[0,\infty)$ that
	\begin{equation}\begin{split}
	\label{eq:EsupIntA1+A2}
	&\limsup_{n\to\infty}\E\bigg[\sup_{s\in[0,t]} \Big| \int_{\floor{s}_{\T_n}}^s
	({A}_1f)(X_r^n)+({A}_2f)(X_r^n,g_n(Z_r^n))dr\Big| \bigg]
	\\     
	\leq&\limsup_{n\to\infty}
	\E\bigg[\sup_{s\in[0,t]}\left|s-{\floor{s}_{\T_n}}\right|^\frac{p-1}p \Big(
	\int_0^s\left|
	({A}_1f)(X_r^n)+({A}_2f)(X_r^n,g_n(Z_r^n))\right|^pdr \Big)^\frac1p \bigg] 
	\\  	
	\leq&\limsup_{n\to\infty}
	\Big(\sup_{s\in[0,\infty)}  \left(s-\floor{s}_{\T_n}\right)\Big)^\frac{p-1}p 
	\Big( \int_0^t\E\left[\left|
	({A}_1f)(X_r^n)+({A}_2f)(X_r^n,g_n(Z_r^n))\right|^p\right] dr \Big)^\frac1p  
	= 0\,.
	\end{split}\end{equation}
	Recall for all $n\in\N$ that $L_n$ is linear, \eqref{eq:form.Ln}, and that $L_{2,n} f_n = 0$. Hence for all  $(x,z)\in S_n\times E_n$
	\begin{equation} 
	\begin{split}
	&\left(L_n\left((f_n+\tfrac{1}{\theta_n}h_n)|_{S_n\times E_n}\right)\right)(x,z)
	=\left(L_n(f_n|_{S_n\times E_n})\right)(x,z)+\tfrac{1}{\theta_n}\left(L_n(h_n|_{S_n\times E_n})\right)(x,z)
	\\
	&=\left(L_{0,n}f_n+\theta_nL_{1,n}f_n\right)(x,z)
	+\left(L_{1,n}h_n+\tfrac{1}{\theta_n}L_{0,n}h_n\right)(x,z)
	+\theta_n(L_{2,n}h_n)(x,z)
	\,.
	\end{split}     
	\end{equation}
	Therefore, we infer
	for all $t\in[0,\infty)$ and all $n\in\N$
	that
	\begin{align} \begin{split} \label{eq:inequality.eps}
	\E\bigg[\sup_{s\in[0,t]}\left|\eps_s^n\right|\bigg] 
	&
	\leq
	\E\bigg[\sup_{s\in[0,t]}\left|\left(f_n+\tfrac{1}{\theta_n}h_n\right)(X_s^n,Z_s^n)-f(X_s^n)\right|\bigg]
	\\ &+
	\E\left[\sup_{s\in[0,t]} \left| \int_{\floor{s}_{\T_n}}^s ({A}_1f)(X_r^n)
	+\left({A}_2 f\right)(X_r^n,g_n(Z_{r}^n)) \,dr\right|\right]
	\\ & +\E\left[\sup_{s\in[0,t]} \left| \int_0^s ({A}_1f)(X_r^n)
	- \left(\theta_nL_{2,n}h_n+\theta_nL_{1,n}f_n+L_{0,n}f_n\right) (X_r^n,Z_r^n)
	 \,dr\right|\right]
	\\ &  + \E\left[ \sup_{s\in[0,t]} \left|\int_0^{s} 
	\left({A}_2 f\right)(X_r^n,g_n(Z_{r}^n))
	-\left(L_{1,n}h_n+\tfrac{1}{\theta_n}L_{0,n}h_n\right)(X_r^n,Z_r^n)
	\,dr\right|\right]
\,.
	\end{split}\end{align}
	Then
	the assumptions~\eqref{eq:fN.converges.to.f}, 
	\eqref{eq:generator.A1} and
	\eqref{eq:generator.A2} 
	together with 
	the calculation in \eqref{eq:EsupIntA1+A2} 
	yield for every $t\in[0,\infty)$
	that the left-hand side
	of~\eqref{eq:inequality.eps} converges to $0$ as $n\to\infty$.
	Having checked all assumptions of Theorem 2.1 in \cite{Kurtz1992}, the
	assertion now follows from Theorem 2.1 in \cite{Kurtz1992}.  This
	finishes the proof of Theorem~\ref{T:stochastic.averaging}.
\end{proof}

In many applications
for each $n\in\N$
there exists a suitable measure 
 $\pi_n\in\mathcal{M}_1(E_n)$
 (typically the ergodic equilibrium of $Z^n$)
 such that one can chose $f_n:=(\tilde S_n \times \tilde E_n\ni (x,z)\mapsto f(x)\in\R)$ and $h_n$ as the solution of the Poisson equation
\begin{equation} \label{eq:Poisson} \begin{split}
L_{2,n}h_n=\int_{E_n} \left(L_{1,n}f_n+\tfrac{1}{\theta_n}L_{0,n}f_n\right)(\cdot,y)\pi_n(dy)
-\left(L_{1,n}f_n+\tfrac{1}{\theta_n}L_{0,n}f_n\right)\,.
\end{split}     \end{equation}
 Corollary~\ref{C:stochastic.averaging} below specializes Theorem~\ref{T:stochastic.averaging} to a situation where the Poisson equation~\eqref{eq:Poisson} has an explicit solution. 
 General Poisson equations have been frequently studied e.g.~in the context of Stein's method
 and there exist conditions implying existence of a solution;
 see, e.g., \cite{GlynnMeyn1996,PardouxVere1,PardouxVere2,VereKulik2012}.
 We also refer to the literature on Stein's method
 where Poisson equations are frequently solved.
 
 Moreover if proving tightness of $(\Gamma_{Z^n})_{n\in\N}$ in
 $\mathcal{L}_m(E)$ is feasible, then one can choose $H=E$
 and $(g_n)_{n\in\N}$ to be the identity functions
 in Assumption~\ref{ass:main}.5.
 In our application of Theorem~\ref{T:stochastic.averaging}
 in Sections~\ref{S:KLM} below,
 informally speaking,
 the processes
 $(X^n)_{n\in\N}$ sense the equilibria of the processes $(Z^n)_{n\in\N}$ only via
 certain real-valued functions
 $(g_n)_{n\in\N}$.
 Proving tightness of $(\Gamma_{g_n(Z^n)})_{n\in\N}$ in $\mathcal{L}_m(\R)$
 is in our application easier than proving tightness
 of $(\Gamma_{Z^n})_{n\in\N}$ in $\mathcal{L}_m(E)$.

In Corollary~\ref{C:stochastic.averaging} below 
we also  include assumptions such that $\Gamma$ in Theorem~\ref{T:stochastic.averaging} is a multiple of the Lebesgue measure on $[0,\infty)$. More precisely we assume that $H=[0,\infty)$, that there exists an operator $\tilde A_2$ such that  $(A_2f)=(S\times [0,\infty) \ni (x,r) \mapsto r(\tilde A_2 f)(x) \in \R)$, $f\in D_0$, and that the following lemma is applicable to $(Y_n)_{n\in\N}=(g_n(Z^n))_{n\in\N}$.

\begin{lemma}[The limiting occupation measure] \label{l:occupation.measure}
	Let $p\in(1,\infty)$, $\mu\in[0,\infty)$ and let $(Y^n)_{n\in\N}=((Y^n_t)_{t\in[0,\infty)})_{n\in\N}$ be a sequence of stationary, stochastic processes with paths in $\mathcal D([0,\infty),[0,\infty))$ such that for all $t\in(0,\infty)$ 
	\begin{equation}
	\sup_{n\in\N} \E[|Y_0^n|^p]<\infty
	\,,\quad
	\lim_{n\to\infty} \E[Y_0^n]=\mu
	\,\text{ and}\quad
	\lim_{n\to\infty} \E\left[\left|\int_{0}^{t}Y_s^n-\E[Y_s^n]ds\right|\right]=0\,.
	\end{equation} 
	Then
	the sequence of occupation measures $(\Gamma_{Y^n})_{n\in\N}$ is tight in $\mathcal L_m([0,\infty))$ and for each limit point~$\Gamma$ of~$(\Gamma_{Y^n})_{n\in\N}$ it holds almost surely  that $\int_{[0,\infty)}y\Gamma(ds,dy)=ds \mu$.
\end{lemma}

\begin{proof}[Proof of Lemma \ref{l:occupation.measure}]
	Stationarity of $(Y_t^n)_{t\in[0,\infty)}$, $n\in\N$, and the Markov inequality imply for all $k\in\N$ that
	\begin{equation}  \begin{split}
	\textstyle\sup_{s\in[0,\infty)}\textstyle\sup_{n\in\N}\P\big(Y_s^n\geq k\big) 
	=\textstyle\sup_{n\in\N}\P\big(Y_0^n\geq k\big)
	\leq\tfrac{1}{k}\textstyle\sup_{n\in\N}\E\big[Y_0^n\big]<\infty
	\,.
	\end{split}     
	\end{equation}
	Therefore $\{Y_s^n\colon n\in\N,\, s\in[0,\infty)\}$ is tight and Remark~\ref{r:sufficient.for.occupation.time} implies that the family
	$(\Gamma_{Y^n})_{n\in\N}$ is tight in $\mathcal{L}_m([0,\infty))$.
	Next observe that due to Fubini's theorem,
	stationarity of $Y^n$ for every $n\in\N$ and Markov's inequality 
	it holds for all $t\in[0,\infty)$ that
	\begin{align} \nonumber
	&\limsup_{K\to\infty}\sup_{n\in\N}
	\E\left[\bigg|\int_{[0,t]\times[0,\infty)}(y-(y\wedge K))\Gamma_{Y^n}(ds,dy)\bigg|
	\right]
	\leq\limsup_{K\to\infty}\sup_{n\in\N}
	\E\left[\int_{[0,t]\times[0,\infty)}y\1_{[K,\infty)}(y)\Gamma_{Y^n}(ds,dy)
	\right]
	\\
	&=
	\limsup_{K\to\infty}\sup_{n\in\N}
	\E\left[\int_0^t Y_s^n \1_{[K,\infty)}(Y_s^n) ds
	\right]
	=
	\limsup_{K\to\infty}\sup_{n\in\N}
	t\E\left[Y^n_0 \1_{\{Y^n_0\geq K\}} \right]
	\leq
	t
	\sup_{n\in\N}
	\E\left[|Y^n_0 |^p \right]
	\lim_{K\to\infty}	K^{1-p}
	=0.
	\label{eq:prep.calc.occup.meas}
	\end{align}
	Now let $\Gamma$ be a limit point of $(\Gamma_{Y^n})_{n\in\N}$ and let $(\Gamma_{Y^{n_k}})_{k\in\N}$ be a subsequence converging weakly to $\Gamma$. Then
	the monotone convergence theorem, Fatou's lemma,
	the fact that $\Gamma=w-\lim_{k\to\infty}\Gamma_{Y^{n_k}}$ together with the fact that for all $m,t\in[0,\infty)$ it holds that $(\mathcal{L}_m([0,\infty))\ni\tilde\Gamma\mapsto\int_{[0,t]\times [0,\infty)}y\wedge m \tilde{\Gamma}(ds,dy)\in\R)\in \mathcal C_b(\mathcal L_m([0,\infty)),\R)$, \eqref{eq:prep.calc.occup.meas} and
	the definition of occupation measures
	imply for all $t\in[0,\infty)$
	that
	\begin{align}  
	&\E\!\left[\left|\int_{[0,t]\times[0,\infty)} y\Gamma(ds,dy)-t\mu\right|\right]
	\leq \liminf_{c,m\to\infty}\E\!\left[\left|\int_{[0,t]\times[0,\infty)} y\wedge m\Gamma(ds,dy)-t\mu\right|\wedge c\right]
	\\& \nonumber
	\leq \liminf_{c,m\to\infty}\limsup_{k\to\infty}\E\!\left[\left|\int_{[0,t]\times[0,\infty)} y\wedge m
	\Gamma_{Y^{n_k}}(ds,dy)-t\mu\right|\wedge c\right]
	\leq \limsup_{k\to\infty}\E\!\left[\left|\int_{[0,t]\times[0,\infty)} y
	\Gamma_{Y^{n_k}}(ds,dy)-t\mu\right|\right]
	\\& \nonumber
	= \limsup_{k\to\infty}\E\!\left[\left|\int_{0}^t 
	Y^{n_k}_s-\mu\,ds\right|\right]
	\leq
	t\limsupn\,\big|\E\big[Y_0^n\big]-\mu\big|
	+
	\limsupn\,\E\bigg[  \bigg| \int_0^{t}
	Y_s^n-\E\big[Y_s^n\big]
	\,ds \bigg|\bigg]
	=0.
	\end{align}
	This implies that it holds a.s.\ for all $t\in[0,\infty)\cap\Q$ that
	$\int_0^t\int_{[0,\infty)}y\Gamma(ds,dy)=t\,\mu$.
	Since $\{[0,t]\colon t\in[0,\infty)\cap\Q\}\subseteq\mathcal{B}([0,\infty))$
	is measure determining, it follows that a.s.\ it holds that
	$\int_{[0,\infty)}y\Gamma(ds,dy)=ds\,\mu$.
\end{proof}
 	
\begin{corollary}\label{C:stochastic.averaging} 
	Let the setting from Assumption~\ref{ass:main} be given with 5.~replaced by 
	\begin{itemize}
		\item[$\tilde 5$.] Let $p\in(1,\infty)$, $\mu\in[0,\infty)$ and for all $n\in\N$ let $g_n\colon E_n\to [0,\infty)$ be Borel measurable functions 
		such that for all $n\in\N$ the process
		$(g_n(Z_r^n))_{r\in[0,\infty)}$ is stationary,
		such that
		$\sup_{n\in\N}\E\left[	\left|g_n(Z_0^n)\right|^p\right]<\infty$,
		and such that for all $t\in[0,\infty)$ it holds that
		$\lim_{n\to\infty}		|\E\left[	g_n(Z_0^n)\right]-\mu|
		+\E\big[\big|\int_0^t	g_n(Z_s^n)-\E[g_n(Z_s^n)]\,ds\big|\big]=0$.
	\end{itemize}
	Let $D_0\subseteq\C_b(S,\R)$ be a dense set in the topology of uniform convergence on compact sets and 
	let $A_1,A_2\colon D_0 \to \C_b(S,\R)$ be functions.
	Suppose for every $f\in D_0$ that $f_n:=\big(\tilde S_n\times \tilde E_n\ni (x,z)\mapsto f(x)	\in\R	\big)\in D_n$, $n\in\N$, that there exists $h_n \in
	D_n$, $n\in\N$ and $\pi_n\in\mathcal{M}_1(E_n)$, $n\in\N$
	such	that $h_n|_{S_n \times E_n}=L_{1,n}f_n$,
	such that for all $n\in\N$, $\phi\in D_n$ and $(x,z)\in S_n\times  E_n$ it holds that
		\begin{equation} \label{eq:def.L_2n.cor}
	(L_{2,n}\phi)(x,z)=\int_{E_n}\phi(x,y)\pi_n(dy)-\phi(x,z)
		\, ,
		\end{equation}
	such that for all  $t\in[0,\infty)$ it holds that
	\begin{align}
	\label{eq:hn.converges.to.0}
	\lim_{n\to\infty} \E\bigg[\sup_{s\in[0,t]}
	\left|\left(\tfrac{1}{\theta_n}h_n\right)(X_s^n,Z_s^n)\right|\bigg]
	&=0\,,
	\end{align}
	such that for all $t\in[0,\infty)$ the integrals in \eqref{eq:L0n-E[L0n]} are well-defined and it holds that
	\begin{equation}
	\begin{split} \label{eq:L0n-E[L0n]}
	&\lim_{n\to\infty} \E\bigg[ \sup_{s\in[0,t]} \bigg| \int_0^{s}
	\bigg[
	\left(L_{0,n}f_n\right)(X_r^n,Z_r^n)
	-\int_{E_n} \left(L_{0,n}f_n\right)(X_r^n,y)\pi_n (dy)
	\bigg]
	\,dr \bigg|\bigg]
	=0\,,
	\end{split}
	\end{equation}
	and such that the integrals in \eqref{eq:generator.A1.stronger} and \eqref{eq:generator.A2.stronger} are well-defined and for all $t\in[0,\infty)$ it holds that
	\begin{align}
	\lim_{n\to\infty}  \sup_{x\in S_n}\bigg| 
	\left({A}_1 f\right)(x) -\int_{E_n} \left(\theta_n
	L_{1,n}f_n+L_{0,n}f_n\right)(x,y)\pi_n (dy)
	\bigg|
	&=0\,,
	\label{eq:generator.A1.stronger}
	\\
	\lim_{n\to\infty} \E\left[\int_{0}^{t}\Big| 
	g_n(Z_{r}^n)\left({A}_2 f\right)(X_r^n)
	-\left(L_{1,n}h_n+\tfrac{1}{\theta_n}L_{0,n}h_n\right)(X_r^n,Z_r^n)
	\Big|\,dr
	\right]
	&=0\,.
	\label{eq:generator.A2.stronger}
	\end{align}
	Then 
	$(X^n)_{n\in\N}$ is relatively compact in $\mathcal D ([0,\infty),S)$ and every limit
	point $(X_t)_{t\in[0,\infty)}$ 
	is a solution of the $\D([0,\infty),S)$-martingale problem for the pre-generator $D_0\ni f\mapsto A_1f+\mu A_2f \in\C_b(S,\R)$.
\end{corollary}

\begin{proof}
	We will apply Theorem~\ref{T:stochastic.averaging} with $H=[0,\infty)$, with $d_H$ being the Euclidean distance, with
	\begin{equation} \label{eq:def.A2.cor}
	A_2^{\text{(Thm)}}=\Big(D_0\ni \phi\mapsto \Big(S\times [0,\infty)\ni(x,r)\mapsto r\big(A_2\phi\big)(x)	\in\R\Big)\in\C(S\times [0,\infty),\R)\Big)
	\end{equation}
	 and all other objects 
	 being defined and named in the same way as in Corollary~\ref{C:stochastic.averaging}.	
	Lemma~\ref{l:occupation.measure} with $Y^n=g_n(Z^n)$ yields that $\tilde 5$. implies Assumption~\ref{ass:main}.$5$.
	For all  $n\in\N$ it holds, as a consequence of \eqref{eq:def.L_2n.cor} and of $f_n$ being a function of the first coordinate only, that $L_{2,n}f_n=0$. Next,  \eqref{eq:hn.converges.to.0} implies \eqref{eq:fN.converges.to.f} and $A_2\colon D_0 \to \C_b(S,\R)$, \eqref{eq:def.A2.cor} and $\tilde 5.$ imply \eqref{eq:Lp.assumption}.
	Due to  $h_n|_{S_n \times E_n}=L_{1,n}f_n$, $n\in\N$, \eqref{eq:def.L_2n.cor} and \eqref{eq:L0n-E[L0n]}
	for all $t\in[0,\infty)$ it holds that
	  \begin{equation}
	    \begin{split} \label{eq:L2.inverse}
	    &\lim_{n\to\infty} \E\bigg[ \sup_{s\in[0,t]} \bigg| \int_0^{s}
	    \bigg[
	    \left( \theta_nL_{2,n}h_n+\theta_nL_{1,n}f_n+L_{0,n}f_n\right)(X_r^n,Z_r^n)
	    \\&
	    \qquad\qquad
	    \qquad\qquad
	    -\int_{E_n} \left(\theta_n
	      L_{1,n}f_n+L_{0,n}f_n\right)(X_r^n,y)\pi_n (dy)
	     \bigg]
	     \,dr \bigg|\bigg]
	    =0\,,
	    \end{split}
	  \end{equation}
	  which implies together with \eqref{eq:generator.A1.stronger} that \eqref{eq:generator.A1} holds. 
	  Finally, \eqref{eq:generator.A2.stronger} implies \eqref{eq:generator.A2} 
	   consequently we have checked all assumptions of  Theorem~\ref{T:stochastic.averaging}. Theorem~\ref{T:stochastic.averaging} now implies that  
	  $(X^n, \Gamma_{g_n(Z^n)})_{n\in\N}$ is relatively compact in $\mathcal D ([0,\infty),S)\times \mathcal L_m(E)$ and for every limit
	  point $((X_t)_{t\in[0,\infty)}, \Gamma)$ 
	  and
	  for every
	  $f\in D_0$ that \eqref{eq:resulting.martingale.main.result} is a martingale. 
	  Finally, due to $\tilde5.$ we can apply Lemma~\ref{l:occupation.measure} with $Y^n=g_n(Z^n)$ to rewrite the double integral in \eqref{eq:resulting.martingale.main.result} as $\int_0^t  ({A}_2 f)(X_s) \, ds\,\mu$, which finishes the proof of Corollary~\ref{C:stochastic.averaging}.	  
\end{proof}

\section{Karlin-Levikson model}
\label{S:KLM}
The well-known neutral Wright-Fisher model describes an
haploid population with non-overlapping generations and a constant size of $n\in\N$ individuals,
each being either of haplotype $A$ or $a$.  
The individuals of generation $k+1$ are sampled independently of each other
from an infinitely large pool in which $A$ and $a$ types have the same proportion as in generation $k\in\N_0$.
Charles Darwin proposed the idea that evolution is due to natural selection.
Thus we aim to study a Wright-Fisher model
with selection.
Gillespie writes in the preface of his book~\cite{Gillespie1993}:
``It is my conviction that the only viable model of selection is one based on temporal and spatial fluctuations
in the environment.''
So it is important to understand the Wright-Fisher model with fluctuating selection for large population sizes.
We will focus on temporally fluctuating selection modeled by a sequence $(\sigma_k,\tau_k)_{k\in\N_0}$ of
$(-1,\infty)^2$-valued random variables.
Generation $k+1$ is sampled in the same way as described above from an infinitely large pool 
to which type $A$ contributes $\frac{1+\sigma_k}{1+\tau_k}$-times as much as its proportion in generation $k\in\N_0$.
This kind of extension of the classical Wright-Fisher model has been studied before 1974, see the references in \cite{KarlinLevikson1974} for details.
Since \cite{KarlinLevikson1974} was the first publication to study the model in full generality,
we follow \cite{Durrett2008} to call it the \emph{Karlin-Levikson model} (KLM).

The purpose of this section is to 
rigorously derive the diffusion limit for the Karlin-Levikson model,
where we allow for autocorrelated selection coefficients.
More precisely, for some $p_n\in(0,1]$ we assume that in each generation with probability
$(1-p_n)$ the selection coefficients
are identical to the coefficients of the previous generation
and with probability $p_n$ the coefficients are sampled from the distribution of the previous selection coefficients.
Thereby,
we have independent selection regimes that last for a geometrically distributed number of generations. 
In the iid-case (that is, $(p_n)_{n\in\N}\equiv 1$), Karlin and Levikson~\cite{KarlinLevikson1974}
derive the coefficients of the SDE~\eqref{eq:SDE.CKL} without giving a formal proof,
cf.~also Theorem 7.12 in~\cite{Durrett2008}.
In the non-iid case (that is, if the sequence of selection coefficients is autocorrelated),
\cite{HuertaEtAl2008} conjecture a diffusion limit 
which is different from~\eqref{eq:SDE.CKL}.
More precisely, 
in order to obtain the diffusion approximation for the non-iid case from the iid-case,
\cite{HuertaEtAl2008} replace $\lim_{n\to\infty}n\E[(\bar{\sigma}_0^n-\bar{\tau}_0^n)^2]/2=\beta/2$ (in the case $p=1$) by
$\lim_{n\to\infty}n\sum_{k=0}^{\infty}\E[(\bar{\sigma}_0^n-\bar{\tau}_0^n)(\bar{\sigma}_k^n-\bar{\tau}_k^n)
=\lim_{n\to\infty}n\E[(\snb-\tnb)^2]\sum_{k=0}^{\infty}(1-p_n)^k=\lim_{n\to\infty}\frac{n}{p_n}\E[(\snb-\tnb)^2]
=\beta$
which differs from $(2-p)\beta/2$ by the factor $(2-p)/2$.
To the best of our knowledge there exists no rigorous proof of the diffusion approximation
of the LKM model in the literature.
The closest result is 
\cite{Taylor2008} where the diffusion approximation is derived
for a sequence of Moran models with fluctuating selection.
Finally, readers interested in  simulations of the KLM are referred to, e.g.,~\cite{KarlinLieberman1974}, \cite{Gillespie1993}
and \cite{GossmannEtAl2014}.
%
%
%
%

The section is structured as follows: We will start by constructing the model, characterize its diffusion approximation in Theorem \ref{T:KLM} and relate our result to some other publications in more detail. After giving an example to show that our assumptions may be fulfilled we will complete the section with proving Theorem \ref{T:KLM}.\\


\begin{definition}[Karlin-Levikson model]\label{def:KLM} \
	Let $(\Omega,\mathcal A,\P)$ be a probability space, 
	let $n\in\N$, 
	let $p\in(0,1]$, 
	let $(\bar X_k,\bar\sigma_k,\bar \tau_k)_{k\in\N_0}$ be a Markov process with state space $\{0,\frac1n,\dots,1\}\times (-1,\infty)^2$
	and let $q:=\big([0,1]\times(-1,\infty)^2\ni (x,\sigma,\tau)\mapsto \tfrac{(1+\sigma)x}{(1+\sigma)x+(1+\tau)(1-x)}\in [0,1]\big)
	$.
	We call $(\bar X_k,\bar\sigma_k,\bar \tau_k)_{k\in\N_0}$ a \emph{Karlin-Levikson model (KLM)} of size $n$ with environmental change probability $p$ if and only if the following conditions hold:
	\begin{enumerate}
		\item
		For all $k\in\N$ it holds that
		\begin{equation} \label{eq:Law.of.sigma,tau}
		\Law{(\bar\sigma_k,\bar \tau_k)\big|(\bar\sigma_{k-1},\bar \tau_{k-1})}
		= p \cdot\Law{\bar\sigma_0,\bar \tau_0} +(1-p) \cdot \delta_{\{(\bar\sigma_{k-1},\bar \tau_{k-1})\}}\,,
		\end{equation}
		i.e.~with probability $p$ the random variable $(\bar\sigma_k,\bar \tau_k)$ is sampled independently of $(\bar\sigma_{k-1},\bar \tau_{k-1})$ according to $\Law{\bar\sigma_0,\bar \tau_0}$ and with probability $(1-p)$ it holds that $(\bar\sigma_k,\bar \tau_k)=(\bar\sigma_{k-1},\bar \tau_{k-1})$.
		\item 
		For all $k\in\N$ it holds that 
		\begin{equation}
		\Law{ n\cdot \bar X_{k} \big| (\bar X_{k-1},\bar\sigma_{k-1},\bar \tau_{k-1})}
		= \operatorname{Bin}(n,q(\bar X_{k-1},\bar\sigma_{k-1},\bar \tau_{k-1}))\,.
		\end{equation}
		\end{enumerate}
\end{definition}

\noindent
	We need to make some assumptions on the law of the selection parameters to establish our convergence result.

\begin{assumption} \label{ass:KLM}
	Let  $p\in[0,1]$, $(p_n)_{n\in\N}\subset (0,1]$,  $\al$, $\gamma\in\R$, $\beta\in[0,\infty)$.
	For every $n\in\N$ let $(\Xnb[k],\snb[k],\tnb[k])_{k\in\N_0}$ be a KLM of size $n$ with environmental change probability $p_n$.
	We assume that $\lim_{n\to\infty} p_n=p$, $\lim_{n\to\infty} np_n=\infty$,
	\begin{align}
	\label{eq:sigma-tau->alpha}
	\lim_{n\to\infty}
	n\E[\snb-\tnb]=
	\al,\,\,
	\lim_{n\to\infty}
	\tfrac n {p_n}\E\big[(\snb-\tnb)^2\big] =  \beta,\,\,
	\lim_{n\to\infty}
	n\E\big[(\snb)^2-(\tnb)^2\big]
	=\gamma\,,
	\displaybreak[0]
	\displaybreak[0]
	\\
	\label{eq:moments.3.4.sigma.tau}
	\lim_{n\to\infty}
	\E\bigg[
	\tfrac{n}{p_n}\sum_{l=2}^{3}\bigg(
	\tfrac{|\snb-\tnb|(|\snb|+|\tnb|)^l}{(1+\snb\wedge\tnb\wedge0)^2}
	+\Big(\tfrac{|\snb-\tnb|}{1+\snb\wedge\tnb}\Big)^{l+1}\bigg)
	\bigg]
	= 0\,,
	\displaybreak[0]
	\\ 
	\label{eq:Lp.cond.sigma.tau}
	\text{that there exists a number } u\in(1,\infty)
	\text{ with }
	\sup_{n\in\N}
	\E\left[\big|	\tfrac{n}{p_n}(\snb-\tnb)^2 \big|^u\right]
	<\infty\,,
	\displaybreak[0]
	\\
	\label{eq:convergence speed|sigma-tau|}
	\text{and for all } \eps\in(0,\infty)
	\text{ that }
	\lim_{n\to\infty}
	\E \Big[n\tfrac{|	\snb-\tnb |}{1+\snb\wedge\tnb}
	\1_{\left\{|	\snb-\tnb |>\eps p_n(1+\snb\wedge\tnb)\right\}}\Big]
	=0\,.
	\end{align}
\end{assumption}
	Note that \eqref{eq:convergence speed|sigma-tau|} may be weakened: Whenever $p\neq0$ one can deduce it  from \eqref{eq:moments.3.4.sigma.tau} using H\"olders inequality and Markov's inequality.
	In Example~\ref{example:sigma.tau} we will construct random variables that fulfill Assumptions~\ref{ass:KLM}.
Now we are prepared to formulate the main result of this section.

\begin{theorem}[Diffusion limit of the KLM] \label{T:KLM}
	Let $((\Xnb[k],\snb[k],\tnb[k])_{k\in\N_0})_{n\in\N}$ be a sequence of Karlin-Levikson models for which Assumption \ref{ass:KLM}  holds with parameters $p\in[0,1]$, $(p_n)_{n\in\N}\subset (0,1]$,  $\al$, $\gamma\in\R$, $\beta\in[0,\infty)$ 
	and assume that $\Xnb[0]$ converges weakly to an $[0,1]$-valued random variable $\xi$ as $n\to\infty$. Then there exists a unique weak solution $(X,W)$ of
	\begin{align}  \nonumber
	dX_t&=X_t(1-X_t)\big[\al-\tfrac \gamma 2 +(2-p)\beta(\tfrac12-X_t)\big]dt
	+\sqrt{X_t(1-X_t)+(2-p)\beta (X_t)^2(1-X_t)^2}dW_t,
	\\ \label{eq:SDE.CKL}
		X_0&=\xi\,,
	\end{align}
	where  $W$ is a standard Brownian motion, and
	$X^n=(X^n_t)_{t\in[0,\infty)}
	:=(\Xnb[\floor{nt}])_{t\in[0,\infty)}
	\xRightarrow{n\to\infty} X$ as $[0,1]$-valued stochastic processes with \cadlag-paths.
\end{theorem}
\noindent
For the sake of simplicity we gather the above assumptions and introduce some notation.
\begin{setting} \label{setting.KLM}
	\begin{enumerate}
		\item 
		Let $(\Xnb[\text{~}],\snb[\text{~}],\tnb[\text{~}])_{n\in\N} = ((\Xnb[k],\snb[k],\tnb[k])_{k\in\N_0})_{n\in\N}$ be a sequence of KLMs for which Assumption \ref{ass:KLM} holds with $p\in[0,1]$, $(p_n)_{n\in\N}\subset (0,1]$,  $\al$, $\gamma\in\R$, $\beta\in[0,\infty)$.\\
		For all $n\in\N$ write $\pi_n:=\Law{\snb,\tnb}$, $S_n:=\{0,\frac1n,\dots,1\}$, $E_n:=(-1,\infty)^2$.
		\item
		For all $n\in\N$, $k\in\N_0$ let
		$
		((S_n\times E_n) \times (\mathcal B (S_n\times E_n))^{\otimes \N}) \ni ((x,\sigma,\tau),A)\mapsto \P\big((\Xnb[\text{~}],\snb[\text{~}],\tnb[\text{~}])\in A 
		\big|
		(\Xnb[k],\snb[k],\tnb[k])=(x,\sigma,\tau)
		\big)
		=:\P_{x,\sigma,\tau}^{n,k}(A)$
		be a regular conditional distribution.
		\\
		For each $(x,\sigma,\tau)\in S_n\times E_n$ we denote by $\E^{n,k}_{x,\sigma,\tau}$ the expectation corresponding to $\P_{x,\sigma,\tau}^{n,k}$ and use the abbreviation $\E^{n,0}_{x,\sigma,\tau}=:\E^{n}_{x,\sigma,\tau}$.
		\item
		Without loss of generality we assume finiteness of all terms in \eqref{eq:sigma-tau->alpha}-\eqref{eq:convergence speed|sigma-tau|} with $\lim_{n\to\infty}$ replaced by $\sup_{n\in\N}$. (Otherwise some statements would only hold for all but finitely many $n\in\N$ instead of for all $n\in\N$, which does not influence the assertion of Theorem~\ref{T:KLM}). 
	\end{enumerate}
\end{setting}

\begin{remark} \label{r:informal.deduction.SDE.KLM}
  We show that
  the coefficients of the SDE~\eqref{eq:SDE.CKL} cannot be derived through expectation and variance of
  the change in frequency in one generation times the time-speedup $n$ in the case $p<1$.
  For all $x\in[0,1]$ the expected change in frequence satisfies in the limit $n\to\infty$ that
	(see \eqref{eq:calc.mom1} and \eqref{eq.mom1.convergence} for a formal proof)
	\begin{equation}
	\begin{split}\label{eq:infinitesimal.mean}
	n\E\big[\Xnb[1]-x\big]
  &=
  n\E\big[ q(x,\snb,\tnb)-x]
  =n\E\big[\tfrac{(1+\snb) x}{(1+\snb) x+(1+\tnb) (1-x)}-x\big]
  =x(1-x)n\E\big[\tfrac{\snb -\tnb}{1+\snb x+\tnb (1-x)}\big]
  \\
	&\approx
	x(1-x)n\E\big[(\snb[0]-\tnb[0])\big(1-\snb x-\tnb (1-x)\big)\big]
  \\
  &
	=
	x(1-x)n\E\big[(\snb[0]-\tnb[0])-\tfrac12((\snb[0])^2-(\tnb[0])^2)+(\snb[0]-\tnb[0])^2(\tfrac12-x)\big]
	\\
	&
	\xrightarrow{n\to\infty}
	x(1-x)(\alpha-\tfrac{\gamma}2+p\beta(\tfrac{1}{2}-x))\,,
	\end{split}
	\end{equation} 
	whereas the variance of the change in frequency satisfies
  in the limit $n\to\infty$ that
	(see \eqref{eq:calc.nE[X-x]^2.to.p.beta}-\eqref{eq:dom.cvg.for.q(1-q)} for a formal proof)
	\begin{equation}
	\begin{split}\label{eq:infinitesimal.variance}
  n\Var(\Xnb[1]-x)
  &=\tfrac{1}{n}\E\big[\Var(n\Xnb[1]|\snb,\tnb)\big]+n\Var\big(\E\big[\Xnb[1]-x|\snb,\tnb\big]\big)
	\\
	&=
	\tfrac{1}{n}\E[nq(x,\snb[0],\tnb[0])(1-q(x,\snb[0],\tnb[0]))]
  +n\Var\big(x(1-x)\tfrac{\snb -\tnb}{1+\snb x+\tnb (1-x)}\big)
  \\
	&\approx
	x(1-x)+n\E\big[(x(1-x)(\snb[0]-\tnb[0]))^2\big]
	\xrightarrow{n\to\infty}
	x(1-x)+p\beta x^2(1-x)^2 \,.
	\end{split}
	\end{equation}
  In the case $p=1$,~\eqref{eq:infinitesimal.mean} and~\eqref{eq:infinitesimal.variance}
  give the drift and the diffusion coefficient, respectively, of the SDE~\eqref{eq:SDE.CKL}.
  In the case $p<1$ of autocorrelated selection coefficients,
  however, the drift coefficient has an additional term $(2-2p)\beta(\tfrac{1}{2}-x)x(1-x)$
  and the diffusion coefficient has an additional term $(2-2p)\beta x^2(1-x)^2$.
  
  Our proof of Theorem~\ref{T:KLM} shows that the drift coefficient of the SDE~\eqref{eq:SDE.CKL}
  can be derived as in~\eqref{eq:infinitesimal.mean}
  but with the identity function replaced by $S_n\times E_n\ni (x,\sigma,\tau)\mapsto x+
  \tfrac{1}{p_n}\E_{x,\sigma,\tau}\big[\Xnb[1]-\Xnb[0]\big]\in\R$, $n\in\N$,
  and similarly for the diffusions coefficient.
  More precisely, for all $x\in[0,1]$ it holds in the limit $n\to\infty$ that
  \begin{equation}  \begin{split}
    &n\E\big[\Xnb[1]+\tfrac{1}{p_n}\E_{\Xnb[1],\snb[1],\tnb[1]}\big[\Xnb[1]-\Xnb[0]\big]
            -x-\tfrac{1}{p_n}\E_{x,\snb[0],\tnb[0]}\big[\Xnb[1]-x\big]\big]
    \\&
    \approx n\E\big[\Xnb[1]-x\big]+\tfrac{n}{p_n}\E\Big[(\Xnb[1]-x)
    \tfrac{\E_{\Xnb[1],\snb[0],\tnb[0]}[\Xnb[1]-\Xnb[0]]-\E_{x,\snb[0],\tnb[0]}[\Xnb[1]-x]}{\Xnb[1]-x}\Big]
            \P\left((\snb[1],\tnb[1])=(\snb,\tnb)\right)
    \\&
    \approx n\E\big[\Xnb[1]-x\big]+\tfrac{n}{p_n}\E\Big[\E_{x,\snb,\tnb}\big[\Xnb[1]-x\big]
    \tfrac{\partial}{\partial x}\E_{x,\snb[0],\tnb[0]}\big[\Xnb[1]-x\big]\Big](1-p_n)
    \\&
    = n\E\big[\Xnb[1]-x\big]+\tfrac{n}{p_n}
    \E\big[\tfrac{\snb-\tnb}{1+\snb x+\tnb(1-x)}x(1-x)\cdot\tfrac{\snb-\tnb}{1+\snb x+\tnb(1-x)}(1-2x)\big]
    (1-p_n)
    \\&
	\xrightarrow{n\to\infty}
	x(1-x)(\alpha-\tfrac{\gamma}2+p\beta(\tfrac{1}{2}-x))+\beta x(1-x)(1-2x)(1-p)
  =
	x(1-x)(\alpha-\tfrac{\gamma}2+(2-p)\beta(\tfrac{1}{2}-x))
  \end{split}     \end{equation}
  and that
  \begin{equation}  \begin{split}
    &n\Var(\Xnb[1]-x)+2n\textup{Cov}\Big(\Xnb[1]-x,\tfrac{1}{p_n}\E_{\Xnb[1],\snb[1],\tnb[1]}\big[\Xnb[1]-\Xnb[0]\big]\Big)
    \\
   & =
    n\Var(\Xnb[1]-x)+\tfrac{2n}{p_n}\textup{Cov}\Big(\Xnb[1]-x,
      \tfrac{\snb[1]-\tnb[1]}{1+\snb[1]\Xnb[1]+\tnb[1](1-\Xnb[1])}\Xnb[1](1-\Xnb[1])\Big)
    \\
   & \approx
    n\Var(\Xnb[1]-x)+\tfrac{2n}{p_n}\E\Big[(\Xnb[1]-x)
      (\snb[1]-\tnb[1])x(1-x)\Big]
    \\
   & =
    n\Var(\Xnb[1]-x)+\tfrac{2n}{p_n}(1-p_n)\E\Big[(\Xnb[1]-x) (\snb[0]-\tnb[0])\Big]x(1-x)
    +\tfrac{2n}{p_n}p_n\E\big[\Xnb[1]-x\big]\E\big[ \snb[1]-\tnb[1]\big]x(1-x)
    \\
   & =
    n\Var(\Xnb[1]-x)+\tfrac{2n}{p_n}(1-p_n)\E\Big[\tfrac{(\snb-\tnb)^2}{1+\snb x+\tnb(1-x)} \Big]x^2(1-x)^2
    +\tfrac{2n}{p_n}p_n\E\big[\Xnb[1]-x\big]\E\big[ \snb[1]-\tnb[1]\big]x(1-x)
    \\
   & \approx
	x(1-x)+p\beta x^2(1-x)^2
    +2(1-p)\beta x^2(1-x)^2+0
	=x(1-x)+(2-p)\beta x^2(1-x)^2\,.
  \end{split}     \end{equation}
	This corresponds to the coefficients of the SDE~\eqref{eq:SDE.CKL}.
\end{remark}

	In the following example we construct random variables for which  Assumption~\ref{ass:KLM} holds.

\begin{example}[Construction of random variables complying with Assumption~\ref{ass:KLM}]
	\label{example:sigma.tau} 
	Let $\al,\gamma\in\R$, 
	$\beta\in[0,\infty)$,
	$a\in[0,\frac13\1_{\{\gamma\neq0\}}+\1_{\{\gamma=0\}})$.
	Let $T,U,V$ and $W$ be random variables with finite fourth moments taking values in $(-\frac12,\infty)$ ($T$,$V$ and $W$) respectively in $\R$ ($U$).
	Assume that $\E[U]=0$, $\E[V-W]=\al$, $\E[U^2]=\beta$, $\E[TU]+\1_{a=0}\E[U|U|]=\gamma$ and if $\gamma=0$ assume $T\equiv0$.
	For all $n\in\N$ let $p_n:=n^{-a}$ and
	define $(-\tfrac34,\infty)^2$-valued random variables
	\begin{equation} \label{eq:example.sigma.tau}
	(\sn,\tn) :=\big(n^{-\frac{1+a} {2}}U^++n^{-\frac{1-a} {2}} \tfrac T2+n^{-1}V,\,
	n^{-\frac{1+a} {2}}U^-+n^{-\frac{1-a} {2}} \tfrac T2 +n^{-1}W\big)
	\,.
	\end{equation}
	We will show that this construction satisfies Assumption~\ref{ass:KLM}.
	For all $n\in\N$ it holds that
	\begin{align}
	&\E[n(\sn-\tn)]
	=\E[n(n^{-\frac{1+a}2}U+n^{-1}(V-W))]
	=n^{\frac{1-a}2}\E[U]+\E[V-W]
	=\al\,.
	\end{align}
	Furthermore, it holds that
	\begin{align}
	\begin{split} \label{eq:example.CKL.beta}
	\lim_{n\to\infty}
	\E\left[\tfrac{n}{p_n}(\sn-\tn)^2\right]
	=
	\lim_{n\to\infty}
	\E\left[n^{1+a}\big(n^{-\frac{1+a}2}U+n^{-1}(V-W)\big)^2\right]
	\\
	=\,
	\E\left[U^2\right]
	+\lim_{n\to\infty}
	\Big[
	2n^{\frac{a-1}2}\E\left[U(V-W)\right]
	+n^{a-1} \E\left[(V-W)^2\right]
	\Big]
	=\beta,
	\end{split}
	\end{align}
	and that
	\begin{align}
	\begin{split}
	&\lim_{n\to\infty}
	\E\left[n\left((\sn)^2-(\tn)^2\right)\right]
	=
	\lim_{n\to\infty}
	\E\left[n\left(\sn-\tn\right)\left(\sn+\tn\right)\right]
	\\
	=&	\lim_{n\to\infty}
	\E\left[	n\left(n^{-\frac{1+a}2}U+n^{-1}(V-W) \right)
	\left( n^{-\frac{1+a}2}|U|+n^{-\frac{1-a}2}T+n^{-1}(V+W) \right)
	\right]
	\\
	=\,&\limn\Big[\E\big[TU\big]
	+ n^{-a}\E\big[U|U|\big]
	+O\Big(n^{-\frac{1-a}{2}}\Big)\Big]
	=\gamma,
	\end{split}
	\end{align}
	which proves \eqref{eq:sigma-tau->alpha}. 
	For all $n\in\N$ it holds that $1+\sn\wedge\tn\geq1+\sn\wedge\tn\wedge0\geq\frac14$ and hence, due to finiteness of fourth moments, for all $l\in\{2,3\}$ that
	\begin{align}\label{eq:example.sigma.tau.higher.moments}
	\begin{split}
	0\leq
	& \,\,
	\tfrac1{4^2+4^4}\limsup_{n\to\infty}\tfrac{n}{p_n}\E\Big[
	\tfrac{|\sn-\tn|(|\sn|+|\tn|)^l}{(1+\sn\wedge\tn\wedge0)^2} +
	\big(\tfrac{|\sn-\tn|}{1+\sn\wedge\tn}\big)^{l+1}
	\Big]
	\\
	\leq
	&\limsup_{n\to\infty}\E\left[\tfrac{n}{p_n}|\sn-\tn|(|\sn|+|\tn|)^l\right]
	\\ 
	\leq
	&\limsup_{n\to\infty}\E\left[n^{1+a}\big(n^{-\frac{1+a}{2}}(|U|+|V|+|W|)\big)
	\big(n^{-\frac{1-a}{2}}|T|+n^{-\frac{1+a}{2}}(|U|+|V|+|W|)\big)^l\right]
	\\ 
	\leq
	&\limsup_{n\to\infty} n^{\frac{1+a}{2}-(1-a)\1_{\{T\not\equiv0\}}-(1+a)\1_{\{T\equiv0\}}}\E\left[\big(|T|+|U|+|V|+|W|\big)^{l+1}\right]
	\end{split}
	\\ \nonumber
	=
	&\lim_{n\to\infty}O\big(\1_{\{T\not\equiv0\}}n^{\frac{-1+3a}{2}}+\1_{\{T\equiv0\}}n^{-\frac{1+a}{2}}\big)
	=0\,,
	\end{align}
	which proves  \eqref{eq:moments.3.4.sigma.tau}.
	Finite absolute third moments of $U$, $V$ and $W$ imply 
	\begin{equation}
	\begin{split}
	\sup_{n\in\N}\E\left[\big|\tfrac n{p_n} (\sn-\tn)^2\big|^\frac32\right]
	&=
	\sup_{n\in\N}
	\E\left[n^{\frac32(1+a)} \big|n^{-\frac{1+a}{2}}U+n^{-1}(V-W)\big|^3\right]
	\\
	&\leq
	\E\Big[(|U|+|V-W|)^3\Big]
	<\infty
	\,, 
	\end{split}
	\end{equation}
	which proves
	\eqref{eq:Lp.cond.sigma.tau} with $u=\frac32$.
	Finally, for all $\eps\in(0,\infty)$ H\"older's inequality, $1+\sn[0]\wedge\tn[0]\geq\frac14$ and Markov's inequality imply
	\begin{equation}\label{eq:convergence.speed.example}
	\begin{split}
	&\qquad\limsup_{n\to\infty}
	\E \left[n\tfrac{|	\sn[0]-\tn[0] |}{1+\sn[0]\wedge\tn[0]}
	\1_{\{|	\sn[0]-\tn[0] |>\eps p_n(1+\sn[0]\wedge\tn[0])\}}\right]	
	\\	&
	\leq 4
	\limsup_{n\to\infty} n
	\left(\E \big[\left(	\sn[0]-\tn[0] \right)^4\big]\right)^{\frac14}
	\Big(\P\big(4\left|	\sn[0]-\tn[0] \right|>\eps p_n\big) \Big)^{\frac34}
	\\
	&\leq
	4 \limsup_{n\to\infty}
	n
	\left(\E \big[\left|	\sn[0]-\tn[0] \right|^4\big]\right)^{\frac14}
	\left(\tfrac{4^4\E [\left|	\sn[0]-\tn[0] \right|^4]}{(p_n)^4\eps^4}
	\right)^{\frac34}
	= 4^4 \eps^{-3}
	\limsup_{n\to\infty}
	np_n^{-3}\E \big[\left|	\sn[0]-\tn[0] \right|^4\big]
	\\
	&
	\leq
	4^4
	\eps^{-3}
	\lim_{n\to\infty}n^{1+3a-4\frac{1+a}{2}}
	\E\big[\left(|U|+|V-W|\right)^4\big]
	=0\,.
	\end{split}
	\end{equation}
	\noindent
	Therefore \eqref{eq:convergence speed|sigma-tau|} is fulfilled as well and we have shown that
	the construction in \eqref{eq:example.sigma.tau}
	is suitable to fulfill Assumption~\ref{ass:KLM}.
\end{example}

The first lemma illustrates what the martingale problem for the KLM looks like and how to find the decomposition of the form \eqref{eq:form.Ln} of the pre-generator. It can easily be adapted to other models constructed via discrete-time Markov processes.

\begin{lemma}[Martingale problem for the KLM] \label{l:(X,s,t)solves.MP(Ln)}
	Let Setting~\ref{setting.KLM} be given,
	let $n\in\N$,
	let	$\Dom(L_n):=$\\
	$\Big\{f \in \C(S_n\times E_n,\R)\colon
	\exists K\in(0,\infty)
	\forall (x,\sigma,\tau)\in S_n\times E_n \,\,
	|f(x,\sigma,\tau)|\leq K\big(1\vee\tfrac{|\sigma-\tau|}{1+\sigma\wedge\tau}\big) \Big\}$,
	 let
	$D_n:=\Big\{f\in\C^{3,0}([0,1]\times (-1,\infty)^2,\R)\colon \exists K\in(0,\infty)
		\forall (x,\sigma,\tau)\in S_n\times E_n
	\forall m\in\{0,1,2,3\}
	\,\,
	\big|\frac{\partial^m}{\partial x^m}f(x,\sigma,\tau)\big|
	\leq K\Big(1\vee \sum_{l=1}^{m\vee1} \big(\frac{|\sigma-\tau|}{1+\sigma\wedge\tau}\big)^l \Big)
	\Big\}\,,$
	let $L_n:=\big(\Dom(L_n) \ni f \mapsto  n\E^{n}_{\cdot,\cdot,\cdot}[f(\Xnb[1],\snb[1],\tnb[1])-f(\cdot,\cdot,\cdot)]\in \C(S_n\times E_n,\R)\big)$, 
	and for $i\in\{0,1,2\}$ let $L_{i,n}:D_n \to \C(S_n\times E_n,\R))$ such that for all $f\in D_n$ and $(x,\sigma,\tau)\in S_n\times E_n$ it holds that
	\begin{align} 
	\begin{split}
	(L_{0,n}f)(x,\sigma,\tau)
	&:=\tfrac n2 \E^n_{x,\sigma,\tau}\big[(\bar X_1^n-x)^2\big] \int_{E_n}(1-p_n)\tfrac{\partial^2 f}{\partial x^2} (x,\sigma,\tau)+ p_n\tfrac{\partial^2 f}{\partial x^2} (x,\zeta,\eta) \pi_n(d(\zeta,\eta))
	\\
	&\hspace{-45pt}+
	n\E^n_{x,\sigma,\tau}\bigg[\int_x^{\bar X_1^n}
	\tfrac12 (\bar X_1^n-t)^2 \int_{E_n}
	(1-p_n)\tfrac{\partial^3 f}{\partial t^3} (t,\sigma,\tau)
	+p_n\tfrac{\partial^3 f}{\partial t^3} (t,\zeta,\eta) \,
	\pi_n(d(\zeta,\eta))\,dt\bigg]
	,
	\\
	(L_{1,n}f)(x,\sigma,\tau)
	&:=\sqrt{\tfrac n {p_n}}\E^n_{x,\sigma,\tau}[\bar X_1^n-x]
	\int_{E_n}
	(1-p_n)\tfrac{\partial f}{\partial x} (x,\sigma,\tau) 
	+p_n\tfrac{\partial f}{\partial x} (x,\zeta,\eta)
	\, \pi_n(d(\zeta,\eta))
	\,,
	\\ 
	 \label{eq:L_{i,n}.KLM}
	(L_{2,n}f)(x,\sigma,\tau)
	&:=\int_{E_n} f(x,\zeta,\eta)-f(x,\sigma,\tau)\,\pi_n(d(\zeta,\eta))\,.
	\end{split}
	\end{align}	
	Then 
	 $(X_t^n,\sigma_t^n,\tau_t^n)_{t\in[0,\infty)}
	 :=(\Xnb[\floor{tn}],\snb[\floor{tn}],\tnb[\floor{tn}])_{t\in[0,\infty)}
	 $ 
	 solves the martingale problem for the pre-generator $L_n$	
	with respect to the set $\T_n:=\frac{\N_0}{n}$.
	Furthermore,  for all 
	$f\in D_n$
	it holds that $f|_{S_n\times E_n}\in \Dom(L_n)$ and that
	$L_n(f|_{S_n\times E_n})=(L_{0,n}+\sqrt{np_n}L_{1,n}+np_nL_{2,n})(f)$.
\end{lemma}

\begin{proof}
	First we prove that $L_n$, $L_{0,n}$, $L_{1,n}$ and $L_{2,n}$ are well-defined linear operators.
	Throughout the whole proof let $F:=\Big\{	 \phi\in\C\big(([0,1]\times(-1,\infty)^2)^2,\R\big) 
	\colon \exists k\in(0,\infty) \,\forall (y,\zeta.\eta,x,\sigma,\tau)\in \big([0,1]\times(-1,\infty)^2\big)^2
	\,\,\, 
	|\phi(y,\zeta.\eta,x,\sigma,\tau)|<k \Big(\Big(\sum_{l=0}^{3} \big( \frac{|\sigma-\tau|}{1+\sigma\wedge\tau} \big)^l\Big)+\Big(\sum_{l=0}^{3} \big( \frac{|\zeta-\eta|}{1+\zeta\wedge\eta} \big)^l\Big)\Big)	\Big\}$.
	Due to \eqref{eq:moments.3.4.sigma.tau} it holds for all 
	$f\in F$ that
	\begin{equation} \label{eq:x.to.int(f(x,s,t)d(s,t).in.C}
	\Big([0,1]^2\times (-1,\infty )^2\ni (x,y,\sigma,\tau)\mapsto\int_{E_n}f(x,\zeta,\eta, y,\sigma,\tau)\pi_n(d(\zeta,\eta))\in\R\Big) \in \C([0,1]^2\times (-1,\infty )^2,\R)\,.
	\end{equation}
	Moreover, Definition~\ref{def:KLM} implies for all $g\in F$ and $(x,\sigma,\tau)\in S_n\times E_n$ that 
	\begin{equation} \label{eq:E[g(.,.,.,.,.,.)]<...}
	\begin{split}
	&\E^n_{x,\sigma,\tau}[g(\Xnb[1], \snb[1],\tnb[1],x,\sigma,\tau)]
	\\
	&=\sum_{k=0}^{n}\tbinom nk q(x,\sigma,\tau)^k(1-q(x,\sigma,\tau))^{1-k}
	\Big[p_n\int_{E_n} g\big(\tfrac kn ,\zeta,\eta, x,\sigma,\tau\big) \pi_n(d(\zeta,\eta)) 
	+(1-p_n) g\big(\tfrac kn ,\sigma,\tau, x,\sigma,\tau\big) \Big]\,.
	\end{split}
	\end{equation}
	Combining \eqref{eq:E[g(.,.,.,.,.,.)]<...}, continuity of $q$ and \eqref{eq:x.to.int(f(x,s,t)d(s,t).in.C} we infer for all $g\in F$ that
	\begin{equation} \label{eq:(x,s,t)toE[g(...)].in.C}
	\big(S_n\times E_n \ni (x,\sigma,\tau)\mapsto\E^n_{x,\sigma,\tau}[g(\Xnb[1], \snb[1],\tnb[1],x,\sigma,\tau)]\in\R\big)
	\in \C(S_n\times E_n,\R)\,.
	\end{equation}
	As a consequence of \eqref{eq:(x,s,t)toE[g(...)].in.C} and \eqref{eq:x.to.int(f(x,s,t)d(s,t).in.C} it holds that $L_{i,n}\colon D_n \to \C(S_n\times E_n,\R)$, $i\in\{0,1,2\}$, and $L_{n}\colon \Dom(L_n) \to \C(S_n\times E_n,\R)$ are well-defined linear operators.
	\\
	Note for all $f\in\Dom (L_n)$ and $k\in\N_0$ that $(L_nf)$ evaluated at $(\Xnb[k],\snb[k],\tnb[k])$ satisfies
	\begin{equation} \label{(Lnf)(Xk,sk,tk)}
	(L_nf)(\Xnb[k],\snb[k],\tnb[k])
	=n\E^{n,k}_{\Xnb[k],\snb[k],\tnb[k]}[f(\Xnb[k+1],\snb[k+1],\tnb[k+1])-f(\Xnb[k],\snb[k],\tnb[k])]\,.
	\end{equation}
	Stationarity of $(\sn[t],\tn[t])_{t\in[0,\infty)}$ and \eqref{eq:moments.3.4.sigma.tau} together with Setting~\ref{setting.KLM}.3 implies for all $f\in\Dom(L_n)$, $K_f\in(0,\infty)$ for which the inequality in the definition of $\Dom(L_n)$ is satisfied, and $t\in[0,\infty)$ that 
	\begin{align} \nonumber
	\E\big[|(L_nf)(X_t^n,\sn[t],\tn[t])|\big] 
	&\leq
	 n \E\big[|\E^{n,\floor{tn}}_{X_t^n,\sn[t],\tn[t]}[f(\Xnb[\floor{tn}+1],\snb[{\floor{tn}+1}],\tnb[{\floor{tn}+1}])]|+|f(X_t^n,\sn[t],\tn[t])|\big]
	 \\
	&\leq
	2nK_f \E\Big[1\vee \tfrac{|\sn[0]-\tn[0]|}{1+\sn[0]\wedge\tn[0]}\Big]
	< \infty\,.
	\end{align}
	Hence, due to Fubini, for all $f\in\Dom(L_n)$ it holds that $\int_0^{\floor{t}_{\T_n}}
	|\left(L_nf\right)(X_s^n,\sigma_s^n,\tau_s^n)|\,ds<\infty$ $\P$-a.s. and that
	\begin{align} 
	\begin{split}
	\Big(f(X_t^n,\sigma_t^n,\tau_t^n) &-\int_0^{\floor{t}_{\T_n}}
	\left(L_nf\right)(X_s^n,\sigma_s^n,\tau_s^n)\,ds\Big)_{t\in[0,\infty)}
	\\ = 
	\Big(f\big(\bar{X}_{\lfloor n t\rfloor}^n,\bar{\sigma}_{\lfloor n t\rfloor}^n,\bar{\tau}_{\lfloor n t\rfloor}^n\big) 
	&-\int_0^{\lfloor n t\rfloor} \tfrac1n \left( L_nf \right)
	(\bar{X}_{\lfloor s\rfloor}^n,\bar{\sigma}_{\lfloor s\rfloor}^n,\bar{\tau}_{\lfloor s\rfloor}^n)\,ds\Big)_{t\in[0,\infty)}
	\end{split} 
	\\ \nonumber
	 = 
	\Big(f\big(\bar{X}_{\lfloor n t\rfloor}^n,\bar{\sigma}_{\lfloor n t\rfloor}^n, \bar{\tau}_{\lfloor n t\rfloor}^n\big) 
	&-\sum_{k=0}^{\lfloor n t\rfloor-1} 
	\E^{n,k}_{\bar{X}_k^n,\bar{\sigma}_k^n, \bar{\tau}_k^n}[f(\bar{X}_{k+1}^n,\bar{\sigma}_{k+1}^n, \bar{\tau}_{k+1}^n) -f(\bar{X}_k^n,\bar{\sigma}_k^n, \bar{\tau}_k^n)]
	\,\Big)_{t\in[0,\infty)} 
	\\
	\nonumber
	= 
	\Big(f\big(\bar{X}_{\lfloor n t\rfloor}^n,\bar{\sigma}_{\lfloor n t\rfloor}^n, \bar{\tau}_{\lfloor n t\rfloor}^n\big) 
	&-\sum_{k=0}^{\lfloor n t\rfloor-1} 
	\E\big[f(\bar{X}_{k+1}^n,\bar{\sigma}_{k+1}^n, \bar{\tau}_{k+1}^n)
	\big|
	(\bar{X}_k^n,\bar{\sigma}_k^n, \bar{\tau}_k^n)\big]
	 -f(\bar{X}_k^n,\bar{\sigma}_k^n, \bar{\tau}_k^n)
	\,\Big)_{t\in[0,\infty)} 
	\end{align}
	is a martingale.
	Therefore $(X_t^n,\sigma_t^n,\tau_t^n)_{t\in[0,\infty)}$ solves the martingale problem for $L_n$ with respect to the set \mbox{$\T_n=\frac{\N_0}{n}$}.
	For all $f\in D_n$ we can use Definition~\ref{def:KLM}.1, Taylor's theorem and Fubini to obtain for all $(x,\sigma,\tau) \in S_n\times E_n$ that
	\begin{align} \nonumber
	(L_n (f&|_{S_n\times E_n}))(x,\sigma,\tau)
	= n \int_{E_n}\E^n_{x,\sigma,\tau}[f(\bar X_1^n,\zeta,\eta)-f(x,\sigma,\tau)]
	\left(p_n\pi_n+(1-p_n)\delta_{\{(\sigma,\tau)\}}\right) 	(d(\zeta,\eta))
	\displaybreak[0]
	\\ \nonumber
	%
	&=  np_n\int_{E_n} \E^n_{x,\sigma,\tau}\left[
	f(\bar X_1^n,\zeta,\eta)-f(x,\zeta,\eta)
	+f(x,\zeta,\eta)-f(x,\sigma,\tau)\right]
	\pi_n(d(\zeta,\eta))
	\\
	\begin{split}
	&\quad+
	n(1-p_n) \E^n_{x,\sigma,\tau}\left[
	f(\bar X_1^n,\sigma,\tau)-f(x,\sigma,\tau)\right]
	\\
	&=
	np_n\int_{E_n} \E^n_{x,\sigma,\tau}
	\Big[
	(\bar X_1^n-x)\tfrac{\partial f}{\partial x} (x,\zeta,\eta) 
	+ \tfrac12(\bar X_1^n-x)^2\tfrac{\partial^2 f}{\partial x^2} (x,\zeta,\eta)\\
	&\quad
	+\int_x^{\bar X_1^n} \tfrac12 (\bar X_1^n-t)^2	\tfrac{\partial^3 f}{\partial t^3} (t,\zeta,\eta) dt
	+ f(x,\zeta,\eta)-f(x,\sigma,\tau)\Big] \,\pi_n(d(\zeta,\eta))
	\\
	&\quad+
	n(1-p_n) \E^n_{x,\sigma,\tau}\Big[
	(\Xnb[1]-x)\tfrac{\partial f}{\partial x} (x,\sigma,\tau) 
	+ \tfrac12(\Xnb[1]-x)^2\tfrac{\partial^2 f}{\partial x^2} (x,\sigma,\tau)
	\\
	&\quad
	+\int_x^{\bar X_1^n} \tfrac12 (\bar X_1^n-t)^2	\tfrac{\partial^3 f}{\partial t^3} (t,\sigma,\tau) dt\Big]
	\\
	&=  n \E^n_{x,\sigma,\tau}[\bar X_1^n-x]
	\int_{E_n}p_n\tfrac{\partial f}{\partial x} (x,\zeta,\eta) +(1-p_n)\tfrac{\partial f}{\partial x} (x,\sigma,\tau) \,\pi_n(d(\zeta,\eta))
	\\
	&\quad+ 
	\tfrac n2\E^n_{x,\sigma,\tau}\big[(\bar X_1^n-x)^2\big]\int_{E_n} p_n\tfrac{\partial^2 f}{\partial x^2} (x,\zeta,\eta) +(1-p_n)\tfrac{\partial^2 f}{\partial x^2} (x,\sigma,\tau)\, \pi_n(d(\zeta,\eta)) 
	\end{split}
	\\ \nonumber
	&\quad
	+
	n\E^n_{x,\sigma,\tau}\bigg[\int_x^{\bar X_1^n} \tfrac12 (\bar X_1^n-t)^2	\int_{E_n} p_n\tfrac{\partial^3 f}{\partial t^3} (t,\zeta,\eta)+(1-p_n)\tfrac{\partial^3 f}{\partial t^3} (t,\sigma,\tau)\, \pi_n(d(\zeta,\eta)) 
	\, dt\bigg]
	\\\nonumber
	&\quad
	+
	np_n\int_{E_n} f(x,\zeta,\eta)-f(x,\sigma,\tau) \,\pi_n(d(\zeta,\eta))
	\displaybreak[0]
	\\ \nonumber
	&=((\sqrt{np_n}L_{1,n}+L_{0,n}+np_nL_{2,n})f)(x,\sigma,\tau)\,.
	\end{align}
	 which completes the proof of Lemma \ref{l:(X,s,t)solves.MP(Ln)}.
\end{proof}
Before proceeding with lemmas that are rather specific for the KLM we provide an auxiliary result that will be helpful for proving condition~\eqref{eq:fN.converges.to.f} or condition~\eqref{eq:hn.converges.to.0} in several examples. 

\begin{lemma}
	\label{l:E.max.iid}
	For every $n\in\N$ let $X_0^n,X_1^n,X_2^n\ldots$ be independent and identically distributed $[0,\infty)$-valued	random variables,
	let $\rho_n\in(0,\infty)$
	and let $\xi^n$ be an independent $\N_0$-valued random variable.
	If for all $\eps\in(0,\infty)$ there exist $\delta\in(0,\infty)$ and $n_0\in\N$ such that for all $n\in\N\cap[n_0,\infty)$ and for all $s\in[\P(X_0^n\leq\eps\rho_n),1]$ it holds that
	\begin{equation}\label{eq:l.max.iid.cond.i=>ii}
	\E\big[1-s^{1+\xi^n}\big]
	\leq 
	\delta^{-1}(\rho_n+1)(1-s)
	\end{equation}
	then i) implies that ii) holds, where
	\begin{enumerate}[i)]
		\item 
		for all $\varepsilon\in(0,\infty)$ it holds that \quad
		$\lim_{n\to\infty}
		\left[\tfrac1{\rho_n}\E[X_0^n]+
		\E[X_0^n\1_{\{X_0^n>\varepsilon\rho_n\}}]
		\right]
		=0\,,
		$
		
		\item 
		it holds that \quad
		$\lim_{n\to\infty}\tfrac1{\rho_n}
		\E\left[\max_{i\in\{0,1,\dots,
			\xi^n
			\}}X_i^n\right]=0\,.
		$
	\end{enumerate} 
	Moreover, if for 
	all $\eps\in(0,\infty)$ there exist $\delta\in(0,\infty)$ and $n_0\in\N$ such that for all $n\in\N\cap[n_0,\infty)$ and for all $s\in[\P(X_0^n\leq\eps\rho_n),1]$ it holds that
	\begin{equation}\label{eq:l.max.iid.cond.ii=>i}
	\E\big[1-s^{1+\xi^n}\big]\geq \delta(\rho_n+1)(1-s)
	\end{equation}
	then ii) implies that i) holds.
	
	In particular, 
	if for all $n\in\N$ it holds that $\rho_n\in(0,\infty)$ and $\xi^n\sim\operatorname{Poi}(\rho_n)$ 
	or if for all $n\in\N$ it holds that $\rho_n\in(0,n]$ and $\xi^n\sim \operatorname{Bin}(n,\frac{\rho_n}{n})$, 
	then in each case
	it holds that
	i) implies ii) and that ii) together with $\sup_{n\in\N}\E[X_0^n]<\infty$ implies i).
\end{lemma}
\begin{proof}
	Independence of all involved random variables implies
	for all $x\in(0,\infty)$, $n\in\N$
	that
	\begin{equation}  \begin{split}
	\P\left(\max_{i\in\{0,1,\ldots,\xi^n\}}X_i^n>x\right)
	=1- \P\left(\max_{i\in\{0,1,\ldots,\xi^n\}}X_i^n \leq x\right)
	=1- \E\left[\P\left(X_0^n \leq x\right)^{1+\xi^n}\right]\,.
	\end{split}
	\end{equation}
	Consequently, using Fubini, the expectation of the maximum is given for all $n\in\N$ by
	\begin{equation}\begin{split} \label{eq:Emax}
	\E\left[\max_{i\in\{0,1,\dots,\xi^n\}}X_i^n\right]
	=\int_0^\infty
	\P\left(\max_{i\in\{0,1,\ldots,\xi^n\}}X_i^n>x\right)
	dx
	=
	\int_0^\infty
	\E\big[1-\P(X_0^n\leq x)^{1+\xi^n}\big]
	dx. 	\end{split}
	\end{equation}
	First we prove that 	ii)  together with \eqref{eq:l.max.iid.cond.ii=>i} implies i).
	For all $n\in\N$, $\eps\in(0,\infty)$ and $x\in(\eps\rho_n,\infty)$ it holds that $\P(X_0^n\leq x)\geq \P(X_0^n\leq \eps\rho_n)$, hence by \eqref{eq:l.max.iid.cond.ii=>i} for all $\eps\in(0,\infty)$ there exist $\delta\in(0,\infty)$, $n_0\in\N$ such that for all $n\in\N\cap [n_0,\infty)$, $x\in(\eps\rho_n,\infty)$ it holds that
	$\E\big[1-\P(X_0^n\leq x)^{1+\xi^n}\big] \geq \delta(\rho_n+1)(1-\P(X_0^n\leq x))$. 
	Therefore
	Fubini, \eqref{eq:Emax} and ii) imply for all $\eps\in(0,\infty)$ that there exists $\delta\in(0,\infty)$ with which it holds that
	\begin{equation}
	\begin{split}
	\limsup_{n\to\infty}
	\E\big[X_0^n\1_{\{X_0^n>\eps\rho_n\}}\big]
	=
	\limsup_{n\to\infty}
	\tfrac{1}{\delta( \rho_n+1)}\int_{\eps\rho_n}^{\infty} \delta (\rho_n+1)(1-\P(X_0^n\leq x))dx
	\\
	\leq
	\limsup_{n\to\infty}
	\tfrac{1}{\delta( \rho_n+1)}\int_{\eps\rho_n}^{\infty} \E\Big[1-\P(X_0^n\leq x)^{1+\xi^n}\Big]dx
	\leq
	\limn
	\tfrac{1}{\delta \rho_n}
	\E\left[\max_{i\in\{0,1,\dots,
		\xi^n
		\}}X_i^n\right]
	=0\,.
	\end{split}
	\end{equation}
	Next we prove that i) together with \eqref{eq:l.max.iid.cond.i=>ii} implies ii). 
	For all $n\in\N$, $\eps\in(0,\infty)$ and $x\in(\eps\rho_n,\infty)$ it holds that $\P(X_0^n\leq x)\geq \P(X_0^n\leq \eps\rho_n)$, hence by \eqref{eq:l.max.iid.cond.i=>ii}  for all $\eps\in(0,\infty)$ there exist $\delta\in(0,\infty)$, $n_0\in\N$ such that for all $n\in\N\cap [n_0,\infty)$, $x\in(\eps\rho_n,\infty)$ it holds that
	$\E\big[1-\P(X_0^n\leq x)^{1+\xi^n}\big] \leq \delta^{-1}(\rho_n+1)(1-\P(X_0^n\leq x))$.
	Therefore \eqref{eq:Emax} and Fubini
	yield that for all $\eps\in(0,\infty)$ there exists $\delta\in(0,\infty)$ such that
	\begin{align}
	 \label{eq:l.max.iid.calc.i=>ii}
	\limsup_{n\to\infty}
	\tfrac1{\rho_n} \E\left[\max_{i\in\{0,1,\dots,\xi^n\}}X_i^n\right]
	&=
	\limsup_{n\to\infty}
	\tfrac1{\rho_n}\int_0^\infty \E\big[1-\P(X_0^n\leq x)^{1+\xi^n}\big] dx
	\\ \nonumber
	&\leq 
	\limsup_{n\to\infty}
	\tfrac1{\rho_n}\int_{0}^{\eps\rho_n}1dx +
	\tfrac1{\rho_n}\int_{\eps\rho_n}^{\infty}\tfrac{1+\rho_n}{\delta} (1-P(X_0^n\leq x))dx
	\\ \nonumber
	&\leq
	\eps 
	+ \limsup_{n\to\infty}
	\Big[
	\tfrac1{\delta\rho_n}\int_0^\infty \P(X_0^n>x)dx 
	+ \int_{\eps\rho_n}^{\infty}\tfrac{1}{\delta} P(X_0^n> x)dx
	\Big]
	\\ \nonumber
	&=
	\eps+\limn
	\tfrac1\delta
	\left[\tfrac1{\rho_n}\E[X_0^n]+
	\E[X_0^n\1_{\{X_0^n>\varepsilon\rho_n\}}]
	\right]
	=\eps\,.
	\end{align}
	Since \eqref{eq:l.max.iid.calc.i=>ii} holds for all $\eps\in(0,\infty)$ this shows that ii) holds and finishes this step.
	
	Next we prove that \eqref{eq:l.max.iid.cond.i=>ii} holds for particular distributions of $\xi^n$.
	Whenever for all $n\in\N$ it holds that $\E[\xi^n]=\rho_n$, equation~\eqref{eq:l.max.iid.cond.i=>ii} holds for all $\delta\in(0,1]$ as for all $s\in[0,1]$ it holds that
	\begin{equation}
	\E\big[1-s^{1+\xi^n}\big]
	=(1-s)\E\big[1+s+\dots+s^{\xi^n}\big]\leq(1-s)\E[\xi^n+1]=(1-s)(\rho_n+1)\,.
	\end{equation}
	In particular, the proven part of the lemma guarantees that i) implies ii) if either for all $n\in\N$ it holds that $\Law{\xi^n}=\operatorname{Poi}(\rho_n)$ 
	or if for all $n\in\N$ it holds that
	$\Law{\xi^n}=\operatorname{Bin}(n,\frac{\rho_n}{n})$ and  $\rho_n\in(0,n]$.\\
	The remainder of the proof deals with showing that ii) together with $\sup_{n\in\N}\E[X_0^n]<\infty$ implies i) if 
	either for all $n\in\N$ it holds that $\Law{\xi^n}=\operatorname{Poi}(\rho_n)$
	or if for all $n\in\N$ it holds that 
	$\Law{\xi^n}=\operatorname{Bin}(n,\frac{\rho_n}{n})$ and  $\rho_n\in(0,n]$. It is sufficient to check \eqref{eq:l.max.iid.cond.ii=>i} and
	we may use ii) in doing so. 
	By  ii) it holds for all $\eps\in(0,\infty)$ that
	\begin{equation} \label{eq:l.max.iid.P>1/2}
	\liminf_{n\to\infty} \P(X_0^n\leq\eps\rho_n)
	=\liminf_{n\to\infty} 
	(1-\P(X_0^n>\eps\rho_n))
	\geq
	\limn\Big(
	1-\tfrac{	\E[X_0^n]}{\eps\rho_n}\Big)=1\,,
	\end{equation}
	where we applied Markov's inequality. Hence it holds for all $\eps\in(0,\infty)$ that
	\begin{equation} \label{eq:l.max.iid.n_0}
	n_0(\eps):=\inf\left\{k\in\N\colon\forall n\in\N\cap[k,\infty) \quad\P(X_0^n\leq\eps\rho_n)\geq\tfrac12\right\}
	<\infty\,.
	\end{equation}
	If for all $n\in\N$ it holds that $\xi^n\sim\operatorname{Poi}(\rho_n)$ then for every $\eps\in(0,\infty)$ 
	it holds
	with $\delta=\frac12\exp(-\frac1\eps\sup_{n\in\N}\E[X_0^n])$ 
	for all $n\in\N\cap[n_0(\eps),\infty)$
	and all $s\in[\P(X_0^n\leq\eps\rho_n),1]$ that
	\begin{align} \nonumber
	\E\big[1-s^{1+\xi^n}\big]
	&=(1-s)\E\big[1+s+\dots+s^{\xi^n}\big]
	\geq(1-s)\E\big[(\xi^n+1)s^{\xi^n}\big]
	\\
	=(1-s)\sum_{k=0}^\infty \tfrac{\partial}{\partial s}(s^{k+1})\tfrac{\rho_n^k}{k!}e^{-\rho_n}
	&=(1-s) \tfrac{\partial}{\partial s}(s e^{(s-1)\rho_n})
	= (1-s)(1+s\rho_n)e^{(s-1)\rho_n} \,\,
	\\ \nonumber
	\geq (1-s)(1+\rho_n)se^{-(1-\P(X_0^n\leq \eps\rho_n))\rho_n}
	&\geq (1-s)(1+\rho_n)\tfrac12e^{-\rho_n\P(X_0^n>\eps\rho_n)}
	\geq \delta(1-s)(1+\rho_n)\,,
	\end{align}
	where we applied Markov's inequality. This shows that \eqref{eq:l.max.iid.cond.ii=>i} holds in the Poisson case.
	
	For the rest of the proof assume for all $n\in\N$ that $\rho_n\in(0,n]$ and $\xi^n\sim \operatorname{Bin}(n,\frac{\rho_n}{n})$.
	Note that for all $\eps\in(0,\infty)$, $n\in\N\cap(\frac1\eps \sup_{k\in\N}\E[X_0^k],\infty)$ and  $s \in [\P(X_0^n\leq\eps\rho_n),1]$ Markov's inequality implies
	\begin{equation}
	\begin{split}
	\left((s-1)\tfrac{\rho_n}{n}+1	\right)^{n}
	\geq
	\left(-\P(X_0^n>\eps\rho_n)\tfrac{\rho_n}{n}+1	\right)^{n}
	\geq
	\big(-\tfrac{\E[X_0^n]}{\eps n}+1	\big)^{n}
	\geq
	\big(-\tfrac{\sup_{k\in\N}\E[X_0^k]}{\eps n}+1	\big)^{n}\,,
	\end{split}
	\end{equation}
	which converges to $\exp(-\frac1\eps\sup_{k\in\N}\E[X_0^k])$ as $n\to\infty$ and hence yields that	
	\begin{equation} \label{eq:l.max.iid.n_1}
	\begin{split}
	n_1(\eps)
	:=
	\inf\Big\{k\in\N\colon\forall n\in\N\cap[k,\infty)\, \forall t&\in\left[\P(X_0^n\leq\eps\rho_n),1\right]
	\\ 
	\left((t-1)\tfrac{\rho_n}{n}+1	\right)^{n}
	&\geq
	\tfrac12\exp\left(-\tfrac{\sup_{l\in\N}\E[X_0^l]}{\eps}\right)\Big\}
	<\infty\,.
	\end{split}
	\end{equation}
	For all $s\in[0,1]$ it holds that $((s-1)\frac{\rho_n}{n}+1)\in[0,1]$ and hence that $((s-1)\frac{\rho_n}{n}+1)^{n-1}\geq((s-1)\frac{\rho_n}{n}+1)^{n}$. 
	Therefore for all $\eps\in(0,\infty)$, for all $n\in\N\cap[n_1(\eps)\vee n_0(\eps) ,\infty)$ with $n_0(\cdot)$ from \eqref{eq:l.max.iid.n_0} and for all $s\in[\P(X_0^n\leq\eps\rho_n),1]$ it holds with $\delta:=\frac14\exp(-\frac1\eps\sup_{l\in\N}\E[X_0^l])$ that
	\begin{equation}
	\begin{split} \label{eq:l.max.iid.cond.ii=>i.Bin}
	\E\big[1-s^{1+\xi^n}\big]
	\geq(1-s)\E\big[(\xi^n+1)s^{\xi^n}\big]
	=(1-s)\sum_{k=0}^\infty \tbinom{n}{k}
	\tfrac{\partial}{\partial s}(s^{k+1})(\tfrac{\rho_n}{n})^k(1-\tfrac{\rho_n}{n})^{n-k} 
	\\
	=(1-s)\tfrac{\partial}{\partial s} \left(s(s\tfrac{\rho_n}{n}+1-\tfrac{\rho_n}{n})^n\right)
	=(1-s)\left[((s-1)\tfrac{\rho_n}{n}+1)^n
	+sn\tfrac{\rho_n}{n}((s-1)\tfrac{\rho_n}{n}+1)^{n-1}\right]
	\\
	\geq 
	(1-s) \tfrac12\exp\left(-\tfrac{\sup_{l\in\N}\E[X_0^l]}{\eps}\right)(1+s\rho_n)
	\geq \delta(\rho_n+1)(1-s)
	\,.
	\end{split}
	\end{equation}
	Thus we have shown \eqref{eq:l.max.iid.cond.ii=>i} for the binomial case,
	which finishes the proof of Lemma \ref{l:E.max.iid}.
\end{proof}
	
The decomposition of $L_n$ constructed in Lemma~\ref{l:(X,s,t)solves.MP(Ln)} already suggests that terms like those we treat in the following lemma will be of interest later on.

\begin{lemma}[Moments of displacement in one step] \label{l:Moments.Displacement}
	Given Setting~\ref{setting.KLM} it holds for all $k\in\N$, $(y,\sigma,\tau)\in S_k\times E_k$, and $f\in\C^1_b([0,1],\R)$ that
	\begin{align}
	\label{eq:E[X-x]<|sigma-tau|/eps}
	\qquad
	\left|\E_{y,\sigma,\tau}^k \left[\bar X_1^k-y\right]\right|
	\leq \tfrac { |\sigma-\tau|}{1+\sigma\wedge\tau }\, ,
	\displaybreak[0]
	\\
	\label{eq:nE[(X-x)^2]}
	\left| k\E_{y,\sigma,\tau}^k\left[(\bar X_1^k-y)^2\right] \right|
	\leq 1+ k\big(\tfrac {\sigma-\tau}{1+\sigma\wedge\tau} \big)^2\,,
	\displaybreak[0]
	\\
	\label{eq:nE[(X-x)^3]}
	\left| k\E^k_{y,\sigma,\tau}\left[(\bar X_1^k-y)^3\right]\right|
	\leq
	\tfrac1k 
	+\tfrac{|\sigma-\tau| }{1+\sigma\wedge\tau}
	+k\big(\tfrac{|\sigma-\tau| }{1+\sigma\wedge\tau}\big)^3\,,
	\displaybreak[0]
	\\
	\label{eq:E[X-x]}
	\lim_{n\to\infty}\sup_{x\in S_n} 
	\left| \E\Big[n\E_{x,\sn,\tn}^n \left[\bar X_1^n-x\right]\Big] 
	- x(1-x) \left[\al -\tfrac{\gamma}{2} +p\beta(\tfrac12-x)
	\right] 
	\right|=0 \, ,
	\displaybreak[0]
	\\
	\label{eq:E[E[X-x]^2-(s-t)^2]}
	\limn\E\bigg[\sup_{x\in S_n}  \tfrac{n}{p_n} \Big|\left(\E^n_{x,\sn,\tn}\left[\Xnb[1]-x\right]\right)^2-x^2(1-x)^2(\sn-\tn)^2\Big|\bigg]=0\,,
	\displaybreak[0]
	\\
	\label{eq:E[(X-x)^2].neu}
	\lim_{n\to\infty}\E\bigg[ \sup_{x\in S_n} 
	\left| n\E_{x,\sn,\tn}^n \left[(\bar X_1^n-x)^2\right]
	-  x(1-x)- n(\sn-\tn)^2 x^2(1-x)^2 
	\right|
	\bigg]=0 \,,
	\displaybreak[0]
	\\
	\label{eq:E[int(X-x)^2f]}
	\lim_{n\to\infty}\E\bigg[ \sup_{x\in S_n} 
	\bigg| \tfrac n{p_n}\E_{x,\sn,\tn}^n \bigg[\int_x^{\Xnb[1]}(\Xnb[1]-v)^2f(v)\,dv\bigg]
	\bigg|
	\bigg]=0 \,.
	\end{align}
\end{lemma}
\begin{proof}
	We define a function $\tilde r:=\left([0,1]\times(-1,\infty)^2\ni (x,\sigma,\tau)\mapsto\tfrac{(\sigma x + \tau(1-x))^2}{1+\sigma x + \tau (1-x)}\in[0,\infty)\right)$. 
	We will exploit several properties of $\tilde r$, namely for all $(x,\sigma,\tau)\in[0,1]\times(-1,\infty)^2$ 
	it holds that
	\begin{align} \label{eq:r.tilde}
	\begin{split}
	\tfrac{1}{1+\sigma x+\tau(1-x)}
	&=
	\tfrac{(1+\sigma x+\tau(1-x))(1-(\sigma x+\tau(1-x)))+(\sigma x+\tau(1-x))^2}{1+\sigma x+\tau(1-x)}
	\\
	&=
	1-\sigma x-\tau (1-x)+  \tilde r(x,\sigma,\tau)\,,
	\end{split}
	\\
	\text{and that \qquad}
	\label{eq:r.tilde.bound}
	\tilde r (x,\sigma,\tau)
	&\leq
	\tfrac{|\sigma|+|\tau|}{1+\sigma\wedge\tau\wedge0}\wedge
	\tfrac{(|\sigma|+|\tau|)^2}{1+\sigma\wedge\tau}\,,
	\end{align}
	where we applied $1\leq\frac1{1+\sigma\wedge\tau\wedge0}$. Furthermore, for all $(x,\sigma,\tau)\in[0,1]\times(-1,\infty)^2$ it holds that
	\begin{align}
	\nonumber
	|(1-\sigma x - \tau(1-x)+\tilde r(x,\sigma,\tau))^2-1|
	&\leq 2 (|\sigma|+|\tau|+\tilde r(x,\sigma,\tau))+(|\sigma|+|\tau|+\tilde r(x,\sigma,\tau))^2
	\\ \label{eq:1-(1-sx-tx-r)^2}
	&\leq 4\tfrac{|\sigma|+|\tau|}{1+\sigma\wedge\tau\wedge0}
	+4\left(\tfrac{|\sigma|+|\tau|}{1+\sigma\wedge\tau\wedge0}\right)^2
	\,.
	\end{align}
	Recall Definition~\ref{def:KLM}, in particular recall $q=\Big([0,1]\times(-1,\infty)^2\ni (x,\sigma,\tau)\mapsto\tfrac{(1+\sigma)x}{(1+\sigma)x+(1+\tau)(1-x)}=\tfrac{(1+\sigma)x}{1+\sigma x +\tau(1-x)}\in[0,1]\Big)$ and that for all $n\in\N$ 
	conditioned on $(\Xnb[0],\snb,\tnb)$ the random variable $n\cdot\bar X_1^n$ is binomially distributed with parameters 
 	$n$ and $q(\Xnb[0],\snb,\tnb)$. Therefore and by \eqref{eq:r.tilde} it holds for all $n\in\N$, $(x,\sigma, \tau)\in S_n\times E_n$ that
	\begin{equation}\begin{split}\label{eq:calc.mom1}
	\E_{x,\sigma,\tau}^n\left[\bar X_1^n-x\right]
	&=\tfrac1n n q(x,\sigma,\tau)-x
	= \tfrac{x+\sigma x}{1+\sigma x +\tau(1-x) }-x
	= \tfrac{(\sigma-\tau)x(1-x)}{ 1+\sigma x + \tau(1-x) }\\
	&=x(1-x)(\sigma-\tau)\big(1-\sigma x - \tau(1-x)+\tilde r(x,\sigma,\tau)\big) \\
	&=x(1-x)(\sigma-\tau)\left(1-\sigma \left(\tfrac12-\left(\tfrac12-x\right)\right) - \tau\left(\tfrac12+\left(\tfrac12-x\right)\right)+\tilde r(x,\sigma,\tau)\right) \\
	&=x(1-x)\left[(\sigma-\tau)-\tfrac{\sigma^2-\tau^2}2 + (\sigma-\tau)^2(\tfrac12-x) +(\sigma-\tau)\tilde{r}(x,\sigma,\tau)\right].
	\end{split}
	\end{equation}
	As the denominator in the first line is bounded below by $1+\sigma\wedge\tau$ this proves
	\eqref{eq:E[X-x]<|sigma-tau|/eps}. Since $x\in[0,1]$ implies
	$|x(1-x)|\leq1$ and $|\frac12-x|\leq1$, it holds that \eqref{eq:calc.mom1}, \eqref{eq:sigma-tau->alpha}, \eqref{eq:r.tilde.bound} and \eqref{eq:moments.3.4.sigma.tau} yield
	\begin{equation} \label{eq.mom1.convergence}
	\begin{split}
	0
	&\leq\limsup_{n\to\infty}\sup_{x\in S_n} 
	\left| \E\left[n\E_{x,\sn,\tn}^n \left[\bar X_1^n-x\right]\right] 
	- x(1-x) \left[\al -\tfrac{\gamma}{2} +p\beta(\tfrac12-x)
	\right] 
	\right|
	\\
	&
	\leq \limsup_{n\to\infty}
	\Big[
	\big|	\E\left[n(\sn-\tn)\right]  -\al \big|
	+ \big| \E\left[n\tfrac{(\sn)^2-(\tn)^2}{2}\right]	-\tfrac{\gamma}{2} \big|
	\\
	&
	\quad+\big|\E\left[p_n\tfrac{n}{p_n} (\sn-\tn)^2\right] 
	-p\beta
	\big|		
	+ \sup_{x\in S_n}
	\big|\E[n(\sn-\tn)\tilde r(x,\sn,\tn)]\big|\Big]
	\\
	&\leq0+
	\limsup_{n\to\infty} 
	\E\left[n|\sn-\tn|
	\tfrac{(|\sn|+|\tn|)^2}{1+\sn\wedge\tn}\right]=0
	\,,
	\end{split}
	\end{equation}
	which implies \eqref{eq:E[X-x]}.
	As for all $x\in[0,1]$ it holds that $x^2(1-x)^2\leq 1$
	it holds that \eqref{eq:calc.mom1}, \eqref{eq:1-(1-sx-tx-r)^2} and \eqref{eq:moments.3.4.sigma.tau} imply 
	\begin{equation} \label{eq:calc.nE[X-x]^2.to.p.beta}
	\begin{split}
	& \limsup_{n\to\infty}\E\bigg[ \sup_{x\in S_n} \tfrac{n}{p_n}
	\Big| \Big(\E_{x,\sn,\tn}^n \left[\bar X_1^n-x\right]\Big)^2
	-  x^2(1-x)^2 (\sn-\tn)^2
	\Big|	\bigg]
	\\
	=
	& \limsup_{n\to\infty}\E\bigg[
	\hspace{-1pt}
	\sup_{x\in S_n} 
	\hspace{-2pt}
	 \tfrac{n}{p_n}
	 \hspace{-2pt}
	\left| x^2(1-x)^2 (\sn-\tn)^2\Big((1-\sn x - \tn(1-x)+\tilde r(x,\sn,\tn))^2-1\Big)
	\right|
	\hspace{-1pt}
	\bigg]
	\\
	\leq
	&\limn\E\Big[\tfrac{n}{p_n}
	\Big| (\sn-\tn)^2\Big(
	4\tfrac{|\sn|+|\tn|}{1+\sn\wedge\tn\wedge0}
	+4\left(\tfrac{|\sn|+|\tn|}{1+\sn\wedge\tn\wedge0}\right)^2
	\Big)
	\Big|
	\Big]
	=0 \,,
	\end{split}\end{equation}
	which proves \eqref{eq:E[E[X-x]^2-(s-t)^2]}.
	Using $\Law{n\cdot\bar X_1^n\big|(\Xnb[0],\snb,\tnb)}= \operatorname{Bin}(n,q(\Xnb[0],\snb,\tnb))$ again, we deduce for all $n\in\N$, $(x,\sigma,\tau)\in S_n \times E_n$ that
	\begin{equation}
	\begin{split}
	&\tfrac1n q(x,\sigma,\tau)(1-q(x,\sigma,\tau)) 
	=\E_{x,\sigma,\tau}^n\big[(\bar X_1^n-\E^n_{x,\sigma,\tau}[\bar X_1^n])^2\big]
	\\
	= &\E_{x,\sigma,\tau}^n\big[(\bar X_1^n-x)^2\big] -2\E^n_{x,\sigma,\tau}[\bar X_1^n-x]\E^n_{x,\sigma,\tau}[\bar X_1^n-x]
	+(\E^n_{x,\sigma,\tau}[\bar X_1^n-x])^2	.
	\end{split}
	\end{equation}
	Hence boundedness of $q$ and \eqref{eq:E[X-x]<|sigma-tau|/eps} imply for all $n\in\N$, $(x,\sigma,\tau)\in S_n \times E_n$ that
	\begin{equation} \label{eq:calc.nE[(X-x)^2]}
	n\E_{x,\sigma,\tau}^n[(\bar X_1^n-x)^2]
	=q(x,\sigma,\tau)(1-q(x,\sigma,\tau))
	+n(\E^n_{x,\sigma,\tau}[\bar X_1^n-x])^2 
	\leq 1+ n\tfrac {(\sigma-\tau)^2}{(1+\sigma\wedge\tau)^2} \,,
	\end{equation}
	which reveals \eqref{eq:nE[(X-x)^2]}
	 and implies with \eqref{eq:calc.nE[X-x]^2.to.p.beta} that for proving \eqref{eq:E[(X-x)^2].neu} it suffices to show 
	\begin{equation} \label{eq:e.tilde}
	\lim_{n\to\infty}\E\Big[ \sup_{x\in[0,1]} 
	\big| q(x,\sn,\tn)(1-q(x,\sn,\tn))
	-  x(1-x)
	\big|
	\Big]=0\,.
	\end{equation}
	Reducing by $(1+\tau)>0$ yields for all $(x,\sigma,\tau) \in [0,1]\times(-1,\infty)^2$ that
	\begin{equation} \begin{split} \label{eq:dom.cvg.for.q(1-q)}
	\left|q(x,\sigma,\tau)(1-q(x,\sigma,\tau)) -x(1-x)\right|
	&=\left| \tfrac{x(1-x)(1+\sigma)(1+\tau)}{(1+ \tau + (\sigma -\tau)x )^2} -x(1-x)\right|
	\\
	=x(1-x)
	\Big|
	\tfrac{1+\tau+\sigma-\tau}{1+\tau}
	\big(1+\tfrac{\sigma-\tau }{1 + \tau}x\big)^{-2} -1
	\Big|
	&\leq
	\Big|
	\big(1+\tfrac{\sigma-\tau}{1+\tau}\big)
	\big(1+\tfrac{\sigma-\tau }{1 + \tau}x\big)^{-2} -1
	\Big|\,,
	\end{split}
	\end{equation}
	which convergences uniformly in $x\in[0,1]$ to zero if $\frac{|\sigma-\tau|}{1+\tau}$ converges to zero. 
	Assumption \eqref{eq:moments.3.4.sigma.tau} implies $\frac{|\sn-\tn|}{1+\tn}\xrightarrow{\text{stoch}}0$ and therefore, by dominated convergence with dominating function $\frac12$, that \eqref{eq:e.tilde} and hence \eqref{eq:E[(X-x)^2].neu} hold.
	We use the third central moment of binomial distribution (e.g. \cite{griffiths2013}) to show for all $n\in\N$, $(x,\sigma,\tau)\in S_n\times E_n$ that
	\begin{equation} \label{eq:calc.mom3}
	\begin{split}
	&n^{-2} q(x,\sigma,\tau)(1-q(x,\sigma,\tau)) (1-2q(x,\sigma,\tau))
	=\E^n_{x,\sigma,\tau}\big[(\bar X_1^n-\E^n_{x,\sigma,\tau}[\bar X_1^n])^3\big]
	\\= &\E^n_{x,\sigma,\tau}\big[(\bar X_1^n-x)^3\big] 
	-3\E^n_{x,\sigma,\tau}\big[(\bar X_1^n-x)^2\big]
	\E^n_{x,\sigma,\tau}[\bar X_1^n-x]
	\\&+3\E^n_{x,\sigma,\tau}[\bar X_1^n-x]
	(\E^n_{x,\sigma,\tau}[\bar X_1^n-x])^2
	-(\E^n_{x,\sigma,\tau}[\bar X_1^n-x])^3\,.
	\end{split}
	\end{equation}
	Since $q(\cdot,\cdot,\cdot)\in[0,1]$ it holds for all $n\in\N$, $(x,\sigma,\tau)\in S_n\times E_n$ that $3q(x,\sigma,\tau)(1-q(x,\sigma,\tau))\leq1$ and hence by \eqref{eq:calc.mom3} and \eqref{eq:calc.nE[(X-x)^2]} that 
	\begin{align} \begin{split}
	\left| n\E^n_{x,\sigma,\tau}\left[(\bar X_1^n-x)^3\right]\right|
	&= \big|n^{-1} q(x,\sigma,\tau)(1-q(x,\sigma,\tau)) (1-2q(x,\sigma,\tau)) 
	\\
	&\,\quad
	+3n\E^n_{x,\sigma,\tau}\big[(\bar X_1^n-x)^2\big]\E^n_{x,\sigma,\tau}[\bar X_1^n-x]
	-2n (\E^n_{x,\sigma,\tau}[\bar X_1^n-x])^3
	\big|
	\end{split}
	\\ \nonumber
	&\hspace{-95pt}\leq
	\tfrac1n
	+\big|3\left(q(x,\sigma,\tau)(1-q(x,\sigma,\tau))
	+n(\E^n_{x,\sigma,\tau}[\bar X_1^n-x])^2\right)\E^n_{x,\sigma,\tau}[\bar X_1^n-x]
	-2n(\E^n_{x,\sigma,\tau}[\bar X_1^n-x])^3\big|
	\\ \nonumber
	&\hspace{-95pt}\leq 
	\tfrac1n
	+ \E^n_{x,\sigma,\tau}[\bar X_1^n-x]
	+ n(\E^n_{x,\sigma,\tau}[\bar X_1^n-x])^3
	\,,
	\end{align}
	which implies with \eqref{eq:E[X-x]<|sigma-tau|/eps} that \eqref{eq:nE[(X-x)^3]} holds.
	Due to Fubini it holds for all $\eps\in(0,\infty)$ that 
	\begin{equation}
	\E\left[	\tfrac{1}{p_n}\tfrac{|\snb-\tnb|}{1+\snb\wedge\tnb}\right]
	=
	\int_0^\infty \P\left(\tfrac{|\snb-\tnb|}{p_n(1+\snb\wedge\tnb)}>x\right)dx
	\leq \eps + \int_\eps^\infty \P\left(\tfrac{|\snb-\tnb|}{p_n(1+\snb\wedge\tnb)}>x\right)dx
	\end{equation}
	and hence, due to  \eqref{eq:convergence speed|sigma-tau|} and $\lim_{n\to\infty} np_n=\infty$, that
	\begin{equation}
	\limsup_{n\to \infty} \E\left[	\tfrac{1}{p_n}\tfrac{|\snb-\tnb|}{1+\snb\wedge\tnb}\right]
	\leq \eps+
	\limsup_{n\to \infty} \tfrac1{np_n}\E\left[n\tfrac{|\snb-\tnb|}{1+\snb\wedge\tnb} \1_{\{|\snb-\tnb|>\eps p_n(1+\snb\wedge\tnb)\}}\right] = \eps\,.
	\end{equation}
	Since $\eps\in(0,\infty)$ is arbitrary, as a consequence of \eqref{eq:convergence speed|sigma-tau|} it holds that 
	\begin{equation} \label{eq:|s-t|/(p_n(1+s.wedge.t))}
	\lim_{n\to\infty}\E\left[	\tfrac{1}{p_n}\tfrac{|\snb-\tnb|}{1+\snb\wedge\tnb}\right]=0\,.
	\end{equation}
	Together with \eqref{eq:nE[(X-x)^3]} and \eqref{eq:moments.3.4.sigma.tau} this implies that
	\begin{equation} \label{eq:E[nE(X-x)^3]}
	0\leq
	\limsup_{n\to\infty} \E\left[\sup_{x\in S_n} \big|\tfrac{n}{p_n}\E^n_{x,\sn,\tn}\big[(\Xnb[1]-x)^3\big]\big|\right]=0\,.
	\end{equation} 
	By Jensen's inequality it holds for all $n\in\N$ that 
	\begin{equation} \label{eq:Jensen.for.|s-t|^2/(p_n)}
	\E\big[\tfrac{1}{p_n}|\tfrac{\sn-\tn}{1+\sn\wedge\tn}|^2\big] \leq \big(\tfrac{1}{np_n}\big)^{\frac12}\big(\E\big[\tfrac{n}{p_n}|\tfrac{\sn-\tn}{1+\sn\wedge\tn}|^4\big]\big)^{\frac12}\,.
	\end{equation}
	Using the fourth central moment of binomial distribution  (e.g.~\cite{griffiths2013}) yields for all $n\in\N$, $(x,\sigma,\tau)\in S_n\times E_n$,  that
	\begin{equation}
	\begin{split}
	& n^{-2}q(x,\sigma,\tau)(1-q(x,\sigma,\tau))
	[1+3q(x,\sigma,\tau)(1-q(x,\sigma,\tau))(n-2)]
	\\
	=\,
	&n\E^n_{x,\sigma,\tau}\big[(\Xnb[1]-\E^n_{x,\sigma,\tau}[\Xnb[1]])^4\big]
	\\
	=\,&n\E^n_{x,\sigma,\tau}\big[(\Xnb[1]-x)^4\big]
	-4n\E^n_{x,\sigma,\tau}\big[(\Xnb[1]-x)^3\big]
	\E^n_{x,\sigma,\tau}\big[\Xnb[1]-x\big]
	\\
	&+6n\E^n_{x,\sigma,\tau}\big[(\Xnb[1]-x)^2\big]
	\big(\E^n_{x,\sigma,\tau}\big[\Xnb[1]-x\big]\big)^2
	-3n\big(\E^n_{x,\sigma,\tau}\big[\Xnb[1]-x\big]\big)^4\,.
	\end{split}
	\end{equation}
	Boundedness of $q$ implies boundedness of the first line by $\frac4n$, hence  \eqref{eq:E[X-x]<|sigma-tau|/eps} - \eqref{eq:nE[(X-x)^3]}, \eqref{eq:|s-t|/(p_n(1+s.wedge.t))},  \eqref{eq:Jensen.for.|s-t|^2/(p_n)} and  \eqref{eq:moments.3.4.sigma.tau}   imply
	\begin{equation}\label{eq:E[(X-x)^4]}
	\begin{split}
	0\leq
	& \limsup_{n\to\infty} \E\Big[\sup_{x\in S_n}\tfrac{n}{p_n}\E^n_{x,\sn,\tn}\big[(\Xnb[1]-x)^4\big]\Big]
	\\
	\leq
	& \limsup_{n\to\infty}\E \Big[\sup_{x\in S_n} \Big( \tfrac{4}{np_n}
	+4\tfrac{n}{p_n}\big|\E^n_{x,\sn,\tn}\big[(\Xnb[1]-x)^3\big]
	\E^n_{x,\sn,\tn}\big[\Xnb[1]-x\big]\big|
	\\
	&\quad+6\tfrac{n}{p_n}\E^n_{x,\sn,\tn}\big[(\Xnb[1]-x)^2\big]
	\big(\E^n_{x,\sn,\tn}\big[\Xnb[1]-x\big]\big)^2
	+3\tfrac{n}{p_n}\big(\E^n_{x,\sn,\tn}\big[\Xnb[1]-x\big]\big)^4
	\Big)\Big]
	\\
	\leq & \limsup_{n\to\infty} \E\Big[
	\tfrac{4}{np_n}+\tfrac{4}{p_n}\Big(n^{-1} +\tfrac{|\sn-\tn|}{1+\sn\wedge\tn} +n\left(\tfrac{|\sn-\tn|}{1+\sn\wedge\tn}\right)^3\Big)\tfrac{|\sn-\tn|}{1+\sn\wedge\tn}
	\\&\quad 
	+\tfrac{6}{p_n}\Big(1+n\left(\tfrac{|\sn-\tn|}{1+\sn\wedge\tn}\right)^2\Big)
	\left(\tfrac{|\sn-\tn|}{1+\sn\wedge\tn}\right)^2
	+\tfrac{3n}{p_n}\left(\tfrac{|\sn-\tn|}{1+\sn\wedge\tn}\right)^4
	\Big]
	\\
	=
	&\limsup_{n\to\infty} \E\Big[\tfrac{4}{np_n}
	\Big(1+\tfrac{|\sn-\tn|}{1+\sn\wedge\tn}\Big)
	+\tfrac{10}{p_n}\left(\tfrac{|\sn-\tn|}{1+\sn\wedge\tn}\right)^2
	+\tfrac{13n}{p_n}\left(\tfrac{|\sn-\tn|}{1+\sn\wedge\tn}\right)^4
	\Big]=0\,,
	\end{split}
	\end{equation}
	Now \eqref{eq:E[int(X-x)^2f]} is a consequence of integration by parts, \eqref{eq:E[nE(X-x)^3]} and \eqref{eq:E[(X-x)^4]} as for all $f\in \C^1_b([0,1],\R)$ it holds that
	\begin{align} \begin{split}
	0\leq& \limsup_{n\to\infty} \E\bigg[\sup_{x\in S_n}\Big|\tfrac{n}{p_n}\E^n_{x,\sn,\tn}\bigg[\int_{x}^{\Xnb[1]}(\Xnb[1]-v)^2f(v)\,dv\bigg]\Big|\bigg]
	\\ 
	=& \limsup_{n\to\infty} \E\bigg[\sup_{x\in S_n}\Big|\tfrac{n}{p_n}\E^n_{x,\sn,\tn}\bigg[\tfrac13(\Xnb[1]-x)^3f(x)+\int_{x}^{\Xnb[1]}\tfrac13(\Xnb[1]-v)^3f'(v)\,dv\bigg]\Big|\bigg]
	\end{split}
	\\ \nonumber
	\leq
	&\limn \E\Big[\sup_{x\in S_n}\Big(\|f\|_\infty \big|\tfrac{n}{3p_n}\E^n_{x,\sn,\tn}\big[(\Xnb[1]-x)^3\big]\big|+
	\|f'\|_\infty \tfrac{n}{12p_n}\E^n_{x,\sn,\tn}\big[|\Xnb[1]-x|^4\big]\Big)\Big]=0\,.
	\end{align}
	This finishes the proof of Lemma~\ref{l:Moments.Displacement}.
\end{proof}

Our next lemma deals with derivatives of the first moment of the displacement in one step.
Terms of this kind arise from the "iterated operators" $L_{1,n}h_n$ and $L_{0,n}h_n$, which are "iterated" since $h_n|_{S_n\times E_n}=L_{1,n}f_n$. Their contribution to the SDE~\ref{eq:SDE.CKL} are characterized in~\eqref{eq:generator.A2.stronger}.

\begin{lemma} \label{l:CKL.derivatives}
	Let $F:=\left([0,1]\times(-1,\infty)^2\ni(x,\sigma,\tau)\mapsto \tfrac{(\sigma-\tau)x(1-x)}{ 1+\sigma x + \tau(1-x) }\in[-1,1]\right)$, let Setting~\ref{setting.KLM} be given
	and let $D_n$ be the set defined in Lemma~\ref{l:(X,s,t)solves.MP(Ln)}.\\
	Then for all $n\in\N$, $(y,\sigma,\tau)\in S_n\times E_n$, $f\in\C^4_b([0,1],\R)$ and $m\in\{0,1,2,3\}$  it holds that
	$F(y,\sigma,\tau) =\E^n_{y,\sigma,\tau}[\Xnb[1]-y]$, 
	that $([0,1]\times(-1,\infty)^2\ni(x,\zeta,\eta)\mapsto F(x,\zeta,\eta)f'(x)\in\R)$ $\in\bigcap_{n\in\N} D_n$ 
	and that
	\begin{align}
	\label{eq:d/dx E[X-x]}
	\left|\tfrac{\partial^m}{\partial y^m}
	F(y,\sigma,\tau)
	\right|
	\leq
	m!
	\sum_{l=(m-1)\vee1}^{m\vee1} \left(\tfrac{|\sigma - \tau|}{1+\sigma\wedge\tau}\right)^l\,,
	\displaybreak[0]
	\\
	\label{eq:n^(1/2)d/dx E[X-x]}
	\lim_{n\to\infty}
	\sum_{l=0}^1
	\sup_{x\in S_n}	
	\left|\E\left[ 	n^{\frac12} \tfrac{\partial^l}{\partial x^l}
	F(x,\sn,\tn)
	\right]\right|
	=0 \,,
	\displaybreak[0]
	\\
	\label{eq:E[X-x]d/dxE[X-x]}
	\lim_{n\to\infty}\E\bigg[ \sup_{x\in S_n} 
	\tfrac{n(1-p_n)}{p_n}\Big| \prod_{l=0}^1\tfrac{\partial^l}{\partial x^l}
	F(x,\sn,\tn)
	- (\sn-\tn)^2  x(1-x)(1-2x) 
	\Big|
	\bigg]=0 \,.
	\end{align}
\end{lemma}

\begin{proof}
	In Lemma~\ref{l:Moments.Displacement} we only assumed that Setting~\ref{setting.KLM} is given, therefore we may use all intermediate results from the proof, in particular we will utilize $\tilde r$ defined above~\eqref{eq:r.tilde}. \\
	From the first line of~\eqref{eq:calc.mom1} and $q([0,1]\times(-1,\infty)^2)\subseteq[0,1]$ we infer that F is well-defined and that for all $n\in\N$, $(y,\sigma,\tau)\in S_n\times E_n$ it holds that $F(y,\sigma,\tau) =\E^n_{y,\sigma,\tau}[\Xnb[1]-y]$.\\
	We calculate the appearing derivatives using the quotient rule and apply \eqref{eq:r.tilde}. For all $(x,\sigma,\tau)\in [0,1]\times(-1,\infty)^2$ it holds that
	\begin{align} 
	\label{eq:calc.d/dx E[X-x]}
	\tfrac{\partial}{\partial x}&F(x,\sigma,\tau)
	=
	\tfrac
	{(\sigma-\tau)(1-2x)(1+\sigma x +\tau(1-x))-(\sigma-\tau)^2x(1-x)}
	{(1+\sigma x+ \tau(1-x))^2}	
	\displaybreak[0]
	\\ \nonumber 
	&=\tfrac
	{(\sigma-\tau)\left(1-2x+\sigma x -2\sigma x^2 +\tau-\tau x -2\tau x +2\tau x^2 -\sigma x +\sigma x^2+\tau x-\tau x^2 \right)}
	{(1+\sigma x+ \tau(1-x))^2}
	\displaybreak[0]
	\\ 	\nonumber
	&=(\sigma-\tau)\tfrac
	{1-2x-\sigma x^2 +\tau -2\tau x +\tau x^2}
	{(1+\sigma x+ \tau(1-x))^2}
	\displaybreak[0]
	\\ \nonumber
	&=
	(\sigma-\tau)
	\big[1-2x-\sigma x^2 +\tau(1-x)^2\big]
	\big[1-\sigma x -\tau (1-x) +\tilde r(x,\sigma,\tau)\big]^2
	\displaybreak[0]
	\\ \nonumber
	&=	
	(\sigma-\tau)(1-2x)+(\sigma-\tau)(-\sigma x^2 +\tau(1-x)^2)
	\displaybreak[0]
	\\ \nonumber
	&\quad+ (\sigma-\tau)\big[1-2x-\sigma x^2 +\tau(1-x)^2\big]
	\big[(1-\sigma x -\tau (1-x) +\tilde r(x,\sigma,\tau))^2-1\big]
	\, ,
	\displaybreak[1]
	\\
	\tfrac{\partial^2}{\partial x^2}&F(x,\sigma,\tau)
	=(\sigma-\tau)\left[\tfrac
	{-2-2\tau-2(\sigma-\tau)x}
	{(1+\sigma x + \tau (1-x))^2}
	-
	\tfrac
	{2(1+\tau-(2+2\tau)x -(\sigma-\tau )x^2)(\sigma-\tau)}
	{(1+\sigma x + \tau (1-x))^3}
	\right]
	\displaybreak[0]
	\\ \nonumber
	&= 
	(\sigma-\tau)
	\tfrac	 {-2}	{(1+\sigma x + \tau (1-x))}	
	-(\sigma-\tau)^2\tfrac	 {2(1+\tau -(2+2\tau)x-(\sigma-\tau)x^2)}	 {(1+\sigma x + \tau (1-x))^3}	
	\, ,
	\displaybreak[1]
	\\
	\label{eq:calc.d^3/dx^3 E[X-x]}
	\tfrac{\partial^3}{\partial x^3}&F(x,\sigma,\tau)
	=\tfrac
	{2(\sigma-\tau)^2}
	{(1+\sigma x+\tau(1-x))^2}
	-\tfrac	 {(\sigma-\tau)^2(-4-4\tau-4\sigma x + 4\tau x)}	 {(1+\sigma x+\tau(1-x))^3}	
	+\tfrac	 {(\sigma-\tau)^3 6(1+\tau-(2+2\tau)x-(\sigma-\tau)x^2)}	 {(1+\sigma x+\tau(1-x))^4}	
	\displaybreak[0]
	\\ \nonumber
	&=(\sigma-\tau)^2\tfrac	{2}	{(1+\sigma x+\tau(1-x))^2}
	+(\sigma-\tau)^2\tfrac	 {4} {(1+\sigma x+\tau(1-x))^2}	
	+(\sigma-\tau)^3\tfrac	 {6(1+\tau-(2+2\tau)x-(\sigma-\tau)x^2)} {(1+\sigma x+\tau(1-x))^4}
	\, .
	\end{align}
	Recall $q=([0,1]\times(-1,\infty)^2\ni (x,\sigma,\tau)\mapsto \tfrac{(1+\sigma)x}{(1+\sigma)x+(1+\tau)(1-x)}\to[0,1])$  and note that for all $(x,\sigma,\tau)\in [0,1]\times(-1,\infty)^2$ 
	and all $m\in\{1,2,3\}$ it holds that
	\begin{equation} \label{eq:boundedness.d/dxE[X-x]}
	\begin{split}
	\big|\tfrac	{1-2x-\sigma x^2 +\tau -2\tau x +\tau x^2}	{(1+\sigma x+ \tau(1-x))^{m+1}}\big|
	&\leq
	\tfrac{1}{(1+\sigma\wedge\tau)^{m-1}}
	\big|\tfrac	{(1+\tau)(1-x)^2-(1+\sigma)x^2 }	{((1+\sigma) x+ (1+\tau)(1-x))^2}\big|
	\\
	=
	\tfrac{1}{(1+\sigma\wedge\tau)^{m-1}}
	\big|\tfrac	{(1-q(x,\sigma,\tau))^2 }	{1+\tau}-\tfrac	{(q(x,\sigma,\tau))^2 }	{1+\sigma}\big|
	&\leq
	\tfrac{1}{(1+\sigma\wedge\tau)^{m}}\,.
	\end{split}
	\end{equation} 
	Together with \eqref{eq:calc.d/dx E[X-x]} - \eqref{eq:calc.d^3/dx^3 E[X-x]} this implies \eqref{eq:d/dx E[X-x]}.
	For all $(x,\sigma,\tau)\in [0,1]\times(-1,\infty)^2$ it holds that
	 $1\leq\frac1{1+\sigma\wedge\tau\wedge0}$
	and hence that \eqref{eq:calc.d/dx E[X-x]} and \eqref{eq:1-(1-sx-tx-r)^2} imply 
	\begin{align} \label{eq:d/dxE[X-x]-(s-t)(1-2x)} \nonumber
	&|\tfrac{\partial}{\partial x}F(x,\sigma,\tau)-(\sigma-\tau)(1-2x)|
	\\ 
	\leq\,
	&|\sigma-\tau|(|\sigma|+|\tau|)	
	+|\sigma-\tau|(1+|\sigma|+|\tau|)4\bigg[\tfrac{|\sigma|+|\tau|}{1+\sigma\wedge\tau\wedge0}+\Big(\tfrac{|\sigma|+|\tau|}{1+\sigma\wedge\tau\wedge0}\Big)^2\bigg]
	\\ \nonumber
	\leq\,
	&|\sigma-\tau|\Big[
	\tfrac{5(|\sigma|+|\tau|)}{1+\sigma\wedge\tau\wedge0}
	+\tfrac{8(|\sigma|+|\tau|)^2+4(|\sigma|+|\tau|)^3}{(1+\sigma\wedge\tau\wedge0)^2}
	\Big]
	\leq 8|\sigma-\tau|\sum_{m=1}^{3} \tfrac{(|\sigma|+|\tau|)^m}{(1+\sigma\wedge\tau\wedge0)^{m\wedge2}}
	\,.
	\end{align}
	Together with \eqref{eq:sigma-tau->alpha}, Jensen's inequality and \eqref{eq:moments.3.4.sigma.tau} this implies
	\begin{align} 
	\begin{split}
	0&\leq\limsup_{n\to\infty} \sup_{x\in S_n}\Big|n^{\frac12}\E\Big[
	\tfrac{\partial}{\partial x}F(x,\sn,\tn)
	\Big]\Big|
	\\
	&\leq\limsup_{n\to\infty} \Big(
	\big|n^{\frac12}\E[(\sn[0]-\tn[0])]\big|
	+
	\sup_{x\in S_n}\Big|n^{\frac12}\E\Big[
	\tfrac{\partial}{\partial x}F(x,\sn,\tn) -(\sn-\tn)(1-2x)
	\Big]\Big|
	\Big)
	\end{split}
	\\ \nonumber
	&\leq\,0+8\limsup_{n\to\infty} \Big(
	\Big(n\E\Big[ 	|\sn-\tn|^2 \big(\tfrac{|\sn|+|\tn|}{1+\sn\wedge\tn\wedge0}\big)^2 \Big]	\Big)^{\frac12}
	+
	n^\frac12\E\Big[
	|\sn-\tn|\sum_{m=2}^{3} \tfrac{(|\sn|+|\tn|)^m}{(1+\sn\wedge\tn\wedge0)^{2}}
	\Big]\Big)
	=0\,.
	\end{align}
	Since $\limn\sup_{x\in S_n}\E [n^{\frac12}|F(x,\sn,\tn)|]=0$ is an immediate consequence of \eqref{eq:E[X-x]}, this shows that \eqref{eq:n^(1/2)d/dx E[X-x]} holds.
	Equations \eqref{eq:r.tilde} and \eqref{eq:r.tilde.bound} imply for all $(x,\sigma,\tau)\in[0,1]\times(-1,\infty)^2$ that
	\begin{equation} \label{eq:E[X-x]-(s-t)x(1-x)}
	\begin{split}
	|F(x,\sigma,\tau)-(\sigma-\tau)x(1-x)|
	&=
	\big|\tfrac{(\sigma-\tau)x(1-x)}{ 1+\sigma x + \tau(1-x) }
	-(\sigma-\tau)x(1-x)\big|
	\\
	\leq
	|\sigma-\tau||1-\sigma x-\tau(1-x)+\tilde r (x,\sigma,\tau)-1|
	&\leq
	2|\sigma-\tau|\tfrac{|\sigma|+|\tau|}{1+\sigma\wedge\tau\wedge0}\,.
	\end{split}
	\end{equation}
	Now triangle inequality, boundedness of $F$ by $1\wedge \frac{|\sigma-\tau|}{1+\sigma\wedge\tau}$, \eqref{eq:d/dxE[X-x]-(s-t)(1-2x)} and \eqref{eq:E[X-x]-(s-t)x(1-x)} imply for all $(x,\sigma,\tau)\in [0,1]\times(-1,\infty)^2$ that
	\begin{align} \nonumber
	&\big|F(x,\sigma,\tau)\tfrac{\partial}{\partial x}F(x,\sigma,\tau) -(\sigma-\tau)^2x(1-x)(1-2x)\big|
	\\ \nonumber
	\leq\,
	&\big|F(x,\sigma,\tau)
	\big(\tfrac{\partial}{\partial x}F(x,\sigma,\tau) -(\sigma-\tau)(1-2x)\big)\big|
	+
	\big|\big(F(x,\sigma,\tau)-(\sigma-\tau)x(1-x)\big)(\sigma-\tau) (1-2x)\big|
	\\
	\begin{split} \label{eq:pf.of.E[X-x]d/dxE[X-x]}
	\leq\,&\left(1\wedge\tfrac{|\sigma-\tau|}{1+\sigma\wedge\tau}\right)
	8|\sigma-\tau|\sum_{m=1}^{3} \tfrac{(|\sigma|+|\tau|)^m}{(1+\sigma\wedge\tau\wedge0)^{m\wedge2}}
	+
	2|\sigma-\tau|\tfrac{|\sigma|+|\tau|}{1+\sigma\wedge\tau\wedge0}|\sigma-\tau|
	\\
	\leq\,&
	16|\sigma-\tau|\big(\tfrac{|\sigma|+|\tau|}{1+\sigma\wedge\tau\wedge0}\big)^2
	+ 8 |\sigma-\tau|\tfrac{(|\sigma|+|\tau|)^3}{(1+\sigma\wedge\tau\wedge0)^2}
	+ 2 |\sigma-\tau|\big(\tfrac{|\sigma|+|\tau|}{1+\sigma\wedge\tau\wedge0}\big)^2
	\,,
	\end{split}
	\end{align}
	which shows that \eqref{eq:moments.3.4.sigma.tau} implies \eqref{eq:E[X-x]d/dxE[X-x]}.
	It remains to show that 
	for all 
	$f\in\C^4_b([0,1],\R)$ it holds that 
	$\phi:=\left([0,1]\times(-1,\infty)^2\ni(x,\sigma,\tau)\mapsto F(x,\sigma,\tau)f'(x)\in\R\right)\in\bigcap_{n\in\N} D_n$
	with
	$D_n=\Big\{g\in\C^{3,0}([0,1]\times (-1,\infty)^2,\R)\colon
	\exists K\in(0,\infty)
	\forall (x,\sigma,\tau)\in S_n\times E_n
	\forall m\in\{0,1,2,3\}
	\,\,
	\big|\frac{\partial^m}{\partial x^m}g(x,\sigma,\tau)\big|
	\leq K\Big(1\vee \sum_{l=1}^{m\vee1} \big(\frac{|\sigma-\tau|}{1+\sigma\wedge\tau}\big)^l \Big)
	\Big\}$, $n\in\N$.
	It is obvious from \eqref{eq:calc.d/dx E[X-x]} - \eqref{eq:calc.d^3/dx^3 E[X-x]} that $\phi\in\C^{3,0}([0,1]\times (-1,\infty)^2,\R)$, the boundedness-condition is fulfilled due to \eqref{eq:d/dx E[X-x]} and boundedness of all appearing derivatives of $f\in\C^4_b([0,1],\R)$.
	This finishes the proof of Lemma \ref{l:CKL.derivatives}.
\end{proof}
Our last lemma gathers two statements: The first one, \eqref{eq:LLN.for.L2inverse.KLM}, will be applied for verifying condition~\eqref{eq:L0n-E[L0n]} (see \eqref{eq:f''terms}-\eqref{eq:check.L2inverse.sigma.tau})
whereas the second one, which is a corollary of the proof of the first one, yields that the last assumption of Lemma~\ref{l:occupation.measure} is satisfied for $(Y_t^n)_{t\in[0,\infty)}=\big(\tfrac{n(1-p_n)}{p_n}(\snb[\floor{tn}]-\tnb[\floor{tn}])^2\big)_{t\in[0,\infty)}$.
\begin{lemma} \label{l:LLN.for.limitting.occupation.measure.CKL}
	Let Setting~\ref{setting.KLM} be given,
  let $\phi\in\C_b([0,1],\R)$ be globally Lipschitz continuous
	and for all $n\in\N$
	let 
	$(\sn[t],\tn[t])_{t\in[0,\infty)}:=(\snb[\floor{tn}],\tnb[\floor{tn}])_{t\in[0,\infty)}$.
	Then for all $t\in[0,\infty)$
	it holds that
	\begin{align} \label{eq:LLN.for.L2inverse.KLM}
	\limsupn\,\E\bigg[ \sup_{s\in[0,t]} \bigg| \int_0^{s}
	\phi(X_r^n)n\Big((\sigma_r^n-\tau_r^n)^2-\E\big[(\sigma_r^n-\tau_r^n)^2\big]\Big)
	\,dr \bigg|\bigg]
	=0\,.
	\\ \label{eq:LLN.for.limitting.occupation.measure.KLM}
	\limsupn\,\E\bigg[ \bigg| \int_0^{t}
	\tfrac{n}{p_n}\Big((\sigma_r^n-\tau_r^n)^2-\E\big[(\sigma_r^n-\tau_r^n)^2\big]\Big)
	\,dr \bigg|\bigg]
	=0\,.
	\end{align}
\end{lemma}
\begin{proof}
	For the whole proof fix $u\in(1,\infty)$ satisfying~\eqref{eq:Lp.cond.sigma.tau} and $t\in(0,\infty)$ 
  (the case $t=0$ is trivial).
	As a consequence of \eqref{eq:Lp.cond.sigma.tau} it holds that $\lim_{n\to\infty}\frac1{\floor{tn}}\E\big[\frac{n}{p_n}(\snb-\tnb)^2\big]=0$ and for all $\eps\in(0,\infty)$ that
	\begin{equation*}
	\begin{split}
	0\leq
	\limsup_{n\to\infty}\E\left[\tfrac{n}{p_n}(\snb-\tnb)^2\1_{\left\{	\frac{n}{p_n}(\snb-\tnb)^2>\eps\floor{tn}	\right\}}\right]
	&\leq
	\limsup_{n\to\infty}\tfrac{\E\left[\left|\frac{n}{p_n}(\snb-\tnb)^2\right|^u\right]}{(\eps\floor{tn})^{u-1}}
	\leq
	\tfrac{\sup_{n\in\N}\E\left[\left|\frac{n}{p_n}(\snb-\tnb)^2\right|^u\right]}	{\lim_{k\to\infty}(\eps\floor{tk})^{u-1}}
	=0
	\,,
	\end{split}
	\end{equation*}
	which is condition~$i)$ of Lemma~\ref{l:E.max.iid}.
	Therefore, Lemma~\ref{l:E.max.iid} and boundedness of $\phi$ imply that
	\begin{equation} \label{eq:pf.LLN.KLM.what.remains.when.summing.instead.of.integrating}
	\begin{split}
	\limsup_{n\to\infty}\tfrac1{p_n}\E\bigg[ \sup_{s\in[0,t]} \bigg|
	\int_s^{\ceil{s}_{\T_n}}
	\phi (X_r^n)n\Big[(\sigma_r^n-\tau_r^n)^2-\E\big[(\sigma_r^n-\tau_r^n)^2\big]\Big]
	\,dr \bigg|\bigg]
	\\
	\leq
	\limn \tfrac2n (\|\phi\|_\infty) \E\bigg[\max_{k\in\{0,\dots,\floor{tn}\}}\tfrac{n}{p_n}(\snb[k]-\tnb[k])^2\bigg]
	=0\,.
	\end{split}
	\end{equation}
	Hence
	for proving \eqref {eq:LLN.for.L2inverse.KLM} it is sufficient to show  that
	\begin{equation}  \label{eq:3.25a}
	\begin{split}
	\limsupn\,
	\E\bigg[\sup_{s\in[0,t]}
	\Big|
	\tfrac1n\sum_{k=0}^{\floor{sn}}
	\phi(\Xnb[k])n\Big((\snb[k]-\tnb[k])^2-\E\big[(\snb[k]-\tnb[k])^2\big]\Big)
	\Big|	\bigg]=0\,.
	\end{split}
	\end{equation}
	For all $n\in\N$ let $a_0^n=1$ and assume without loss of generality that, on a possibly larger probability space, there exists a sequence of independent Bernoulli random variables $(a_k^n)_{k\in\N}$ with success probability $p_n$ and an independent sequence of independent and identically distributed random variables $(\snh[k],\tnh[k])_{k\in\N}$ with $\Law{\snh,\tnh}=\pi_n$	
	such that for all $k\in\N$ it holds that
	$(\snb[k+1],\tnb[k+1])=a_{k+1}^n(\snh[k+1],\tnh[k+1])+(1-a_{k+1}^n)(\snb[k],\tnb[k])$. For each $n\in\N$, $k\in\N_0$ let $b_k^n:=\inf\{m\in\N_0\cap[k,\infty]\colon a_{m+1}^n=1\}\in\N_0\cup\infty$ and
	note that the random variables $b_k^n-k$ are geometrically distributed with success probability $p_n$. We calculate the first two moments of the geometric distribution. For all $n\in\N$, $k\in\N_0$ it holds that
	\begin{align} \nonumber
	&\E[b_k^n-k]
	=\sum_{m=0}^{\infty}mp_n(1-p_n)^m
	=-p_n(1-p_n)\tfrac{\del}{\del p_n}\sum_{m=0}^{\infty}(1-p_n)^m
	=	-p_n(1-p_n)\tfrac{\del}{\del p_n} \tfrac{1}{1-(1-p_n)}
	=\tfrac{1-p_n}{p_n}\,,
	\\
	\begin{split}\label{eq:moments.b_k^n}
	&\E[(b_k^n-k)^2]
	=\sum_{m=0}^{\infty}m^2p_n(1-p_n)^m
	=-p_n(1-p_n)\tfrac{\del}{\del p_n} \sum_{m=0}^\infty m(1-p_n)^m
	\\
	&\quad=p_n(1-p_n) \tfrac{\del}{\del p_n}((1-p_n)\tfrac{\del}{\del p_n} \sum_{m=0}^\infty (1-p_n)^m)
	=-p_n(1-p_n) \tfrac{\del}{\del p_n } ((1-p_n)\tfrac1{(p_n)^2})
	\\
	&\quad
	=-p_n(1-p_n)\tfrac{-(p_n)^2-2p_n(1-p_n)}{(p_n)^4}
	=-\tfrac{(1-p_n)(p_n-2)}{(p_n)^2}
	=\tfrac{2-3p_n+(p_n)^2}{(p_n)^2}\,.
	\end{split}
	\end{align}
	For all $n\in\N$ it holds that
	\begin{align} \label{eq:CKL.L2-inverse}
	&\E\bigg[\sup_{s\in[0,t]}
	\Big|
	\tfrac1n\sum_{k=0}^{\floor{sn}}
	\phi(\Xnb[k])\Big(n(\snb[k]-\tnb[k])^2-\E\big[n(\snb[k]-\tnb[k])^2\big]\Big)
	\Big|	\bigg]
	\displaybreak[0]
	\\ 
	&=
	\E\bigg[\sup_{s\in[0,t]}
	\Big|
	\tfrac1n\sum_{k=0}^{\floor{sn}} a_k^n
	\Big(\sum_{l=k}^{b_k^n\wedge\floor{sn}}\phi(\Xnb[l])\Big) \Big(n(\snb[k]-\tnb[k])^2-\E\big[n(\snb[k]-\tnb[k])^2\big]\Big)
	\Big|	\bigg]
	\displaybreak[0]
	\\ \label{eq:CKL.L2-inverse.Hölder}
	&\leq
	\E\bigg[\sup_{s\in[0,t]}
	\Big|
	\tfrac1n\sum_{k=0}^{\floor{sn}} a_k^n
	\Big(\sum_{l=k}^{b_k^n\wedge\floor{tn}}\phi(\Xnb[l])-\phi(\Xnb[k])\Big) \Big(n(\snb[k]-\tnb[k])^2-\E\big[n(\snb[k]-\tnb[k])^2\big]\Big)
	\Big|	\bigg]
	\displaybreak[0]
	\\ \label{eq:CKL.L2-inverse.martingale}
	&\quad+
	\E\bigg[\sup_{s\in[0,t]}\Big|\tfrac1n \sum_{k=0}^{\floor{sn}}\phi(\Xnb[k])a_k^n(b_k^n\wedge\floor{tn}-k+1)\Big(n(\snb[k]-\tnb[k])^2-\E\big[n(\snb[k]-\tnb[k])^2\big]\Big)\Big|\bigg]
	\displaybreak[0]
	\\ \label{eq:CKL.L2-inverse.zero-run}
	&\quad+
	\E\bigg[\sup_{s\in[0,t]}
	\Big|
	\tfrac1n\sum_{k=0}^{\floor{sn}} a_k^n
	\Big(-\sum_{l=b_k^n\wedge\floor{sn}+1}^{b_k^n\wedge\floor{tn}} \phi(\Xnb[l])\Big) \Big(n(\snb[k]-\tnb[k])^2-\E\big[n(\snb[k]-\tnb[k])^2\big]\Big)
	\Big|	\bigg] \,,
	\end{align}
	where the inner sum in the last line is by convention zero if $b_k^n\wedge\floor{sn}=b_k^n\wedge\floor{tn}$.
	We will first take care of \eqref{eq:CKL.L2-inverse.martingale} and therefore show that for all $n\in\N$ the process
	\begin{equation*}
	\begin{split}
	(M_k^n)_{k\in\{0,\dots,\floor{tn}\}}:=
	\left(\tfrac1{p_n}\sum_{m=0}^{k}\phi(\Xnb[m])a_m^n(b_m^n\wedge\floor{tn}-m+1) \left(n(\snb[m]-\tnb[m])^2-\E\big[n(\snb[m]-\tnb[m])^2\big]\right) \right)_{k\in\{0,\dots,\floor{tn}\}}
	\end{split}
	\end{equation*}		
	is a martingale with respect to 
	$\F_k^n:=\sigma(\Xnb[m],\snb[m],\tnb[m],b_m^n, m\in\{0,\dots, k\})$,
  $k\in\{0,\ldots,\lfloor tn\rfloor\}$.
  Adaptedness is obvious, integrability is due to \eqref{eq:Lp.cond.sigma.tau} implied by integrability of $\frac{n}{p_n}(\snb[m]-\tnb[m])^2$ for all $m\in\N_0$ and by boundedness of the remaining terms, martingale property is a consequence of the following observations about dependencies between all involved random variables in $M_{k}^n-M_{k-1}^n$:
	For $n\in\N$, $k\in\{1,\dots,\floor{tn}\}$ it holds that
	$\Xnb[k]$ depends only through $(\Xnb[k-1],\snb[k-1],\tnb[k-1])$ on the other terms, 
	that $b_{k}^n$ depends only on $(a_{m}^n)_{m\in\N\cap[k+1,\infty)}$,
	and that
	$a_{k}^n(n(\snb[k]-\tnb[k])^2-\E[n(\snb[k]-\tnb[k])^2])$ is either zero (if $a_{k}^n=0$) or independent  of $b_{k}^n$ and of $\Xnb[k]$, hence it holds that
	\begin{align*}
	\begin{split}
	\E[M_{k}^n-M_{k-1}^n|\F_{k-1}^n]
	=\,
	&\E\big[	\phi(\Xnb[k])		(b_{k}^n\wedge\floor{tn}-k+1)\big|\F_{k-1}^n\big]
	\E\Big[ a_{k}^n\Big(\tfrac{n}{p_n}(\snb[k]-\tnb[k])^2-\E\Big[\tfrac{n}{p_n}(\snb[k]-\tnb[k])^2	\Big]\Big)
	\Big|\F_{k-1}^n\Big]
	=0\,. 
	\end{split}
	\end{align*}
	Let
	$A_m^n:=\phi(\Xnb[m])a_m^n(b_m^n\wedge\floor{tn}-m+1) \big(\frac{n}{p_n}(\snb[m]-\tnb[m])^2-\E\big[\frac{n}{p_n}(\snb[m]-\tnb[m])^2\big]\big) $,
	$n\in\N$, $m\in\N_0$ and let 
	$v\in(1,2\wedge u)$.
	As a consequence of \eqref{eq:moments.b_k^n} it holds for all $n\in\N$, $m\in\N_0$ that $\E[(b_m^n-m+1)^2]
	=\E[(b_m^n-m)^2]+2\E[b_m^n-m]+1
	=\frac{2-3p_n+(p_n)^2+2p_n(1-p_n)+(p_n)^2}{(p_n)^2}
	=\frac{2-p_n}{(p_n)^2}$.
	Hence it holds by H\"older's, Markov's and Jensen's inequalities, the above observations about dependencies and \eqref{eq:Lp.cond.sigma.tau} that
	\begin{align}
	\nonumber
	&\limsup_{a\to\infty}
	\sup_{n\in\N}\sup_{m\in\{0,\dots,\floor{tn}\}}
	\E\big[
	\big|A_m^n	\big|^{\sqrt{v}} \1_{\{|A_m^n|^{\sqrt{v}}>a\}}\big]
	\leq
	\limsup_{a\to\infty}
	\sup_{n\in\N}\sup_{m\in\{0,\dots,\floor{tn}\}}
	\big(\E\big[
	|A_m^n	|^v\big]\big)^{\frac{1}{\sqrt{v}}}	\P\big(|A_m^n|^{\sqrt{v}}>a\big)^{\frac{\sqrt{v}-1}{\sqrt{v}}}
	\\
	\begin{split}
	&\leq
	\limsup_{a\to\infty}
	\sup_{n\in\N}\sup_{m\in\{0,\dots,\floor{tn}\}}
	\E\big[
	|A_m^n	|^v\big] 
	a^{-\sqrt{v}\frac{\sqrt{v}-1}{\sqrt{v}}}
	\\
	&\leq
	\|\phi\|_\infty^v
	\sup_{n\in\N}\sup_{m\in\{0,\dots,\floor{tn}\}}
	\E[a_m^n]\big(\E\big[|b_m^n-m+1|^2\big]\big)^{\frac v2} 2^v\E\Big[\big|n(\snb[m]-\tnb[m])^2\big|^v\Big]
	\lim_{a\to\infty}a^{1-\sqrt{v}}
	\\
	&=
	\|\phi\|_\infty^v
	\sup_{n\in\N}
	p_n\big(\tfrac{2-p_n}{p_n^2}\big)^{\frac v2} 2^v(p_n)^v\E\Big[\big|\tfrac n{p_n} (\snb[0]-\tnb[0])^2\big|^v\Big]
	\lim_{a\to\infty}a^{1-\sqrt{v}}=0\,.
	\end{split}
	\end{align}
	This implies that the family
	$
	\Big\{ \big|A_m^n
	\big|^{\sqrt{v}},m\in \{0,\dots,\floor{tn}\},n\in\N\Big\}
	$
	is uniformly integrable and that we may apply \cite{Gut1992}'s weak LLN,
  which implies together with Jensen's inequality and Doob's $L^p$-inequality  that
	\begin{align} \label{eq:3.25b.term2}
	\begin{split}
	&\limsupn
	\bigg(\tfrac{1}{p_n}\E\bigg[\sup_{s\in[0,t]}\Big|\tfrac1n \sum_{k=0}^{\floor{sn}}\phi(\Xnb[k])a_k^n(b_k^n\wedge\floor{tn}-k+1)
	\Big(n(\snb[k]-\tnb[k])^2-\E\big[n(\snb[k]-\tnb[k])^2\big]\Big)
	\Big|\bigg]\bigg)^{\sqrt{v}}
	\\
	&\leq
	\limsupn\,
	\E\bigg[\sup_{k\in\{0,\dots,\floor{tn}\}} \big|\tfrac1n M_k^n\big|^{\sqrt{v}}\bigg]
	\leq
	\left(\tfrac{\sqrt{v}}{\sqrt{v}-1}\right)^{\sqrt{v}}
	\limsupn\, 
	\E\Big[ \big|\tfrac1nM_{\floor{tn}}^n\big|^{\sqrt{v}}\Big]
	\\
	&=
	\left(\tfrac{\sqrt{v}}{\sqrt{v}-1}\right)^{\sqrt{v}}
	\limn
	\tfrac{\floor{tn}+1}{n} \left(\tfrac{1}{n}\right)^{\sqrt{v}-1} \E\Bigg[\tfrac{1}{\floor{tn}+1}\Big|\sum_{m=0}^{\floor{tn}}A_m^n\Big|^{\sqrt{v}}\Bigg]
	=0
	\,,
	\end{split}
	\end{align}
	which is convergence of \eqref{eq:CKL.L2-inverse.martingale}.
	Let $K\in(0,\infty)$ be a Lipschitz constant for $\phi$ and 
	for all $k\in\N_0$ and $n\in\N$ let $c_k^n:=\sup\{m\in\N_0\cap[0,k]\colon a_m^n=1\}$. 
	For all $n \in\N$, $k\in\N_0$ and $l\in\{0,\dots,k-1\}$
	it holds due to stationarity of $(\snb[~],\tnb[~])$ that
	\begin{equation}
	\E\big[|\Xnb[k]-\Xnb[l]|\,\big|\,\{c_k^n=l\}\big]
	\leq \sum_{m=l}^{k-1}\E\big[|\Xnb[m+1]-\Xnb[m]|\,\big|\,\{c_k^n=l\}\big]
	\leq
	(k-l)\E\left[\sup_{x\in S_n}\E_{x,\snb,\tnb}^n[|\Xnb[1]-x|]\right]
	\end{equation}
	and hence for all $u\in(1,\infty)$, $n\in \N$ and $k\in\N_0$ due to $p_n\in(0,1]$, $|\Xnb[m+1]-\Xnb[m]|\leq1$, Jensen's inequality, \eqref{eq:nE[(X-x)^2]} and \eqref{eq:moments.b_k^n} that
	\begin{align} \nonumber
	&\E\Big[	\big|p_n\big(\phi(\Xnb[k])-\phi(\Xnb[c_k^n])\big)\big|^{\frac{u}{u-1}}	\Big]
	\leq
	K^{\frac{u}{u-1}}p_n
	\E\big[	\big|	\Xnb[k]-\Xnb[c_k^n]	\big|	\big]
	\leq
	K^{\frac{u}{u-1}}p_n
	\sum_{l=0}^{k}\P(c_k^n=l) \E\Big[	\big|	\Xnb[k]-\Xnb[l]	\big| \Big| \{c_k^n=l\}	\Big]
	\\
	\begin{split}
	&\leq
	K^{\frac{u}{u-1}}p_n
	\sum_{l=0}^{k}\P(c_k^n=l) (k-l)
	\E\Big[\sup_{x\in S_n}	\E_{x,\snb,\tnb}^n\big[	|\Xnb[1]-x|\big] \Big]
	\\
	&\leq
	K^{\frac{u}{u-1}}p_n
	\Big(\E\Big[\sup_{x\in S_n} \E_{x,\snb,\tnb}^n\big[	(\Xnb[1]-x)^2\big] \Big]\Big)^{\frac12}
	\sum_{l=-\infty}^{k} (k-l) p_n(1-p_n)^{k-l}
	\\	&
	\leq
	K^{\frac{u}{u-1}}p_n
	\Big(\E\left[\tfrac{1}{n}+\big(\tfrac{\snb-\tnb}{1+\snb\wedge\tnb}\big)^2\right]\Big)^{\frac12}
	\tfrac{1-p_n}{p_n}
	\,.
	\end{split}
	\end{align}
	Therefore H\"older's inequality with exponent $u$, \eqref{eq:moments.3.4.sigma.tau} and \eqref{eq:Lp.cond.sigma.tau} imply that
	\begin{align} \label{eq:3.25b.term1} 
	\begin{split}
	&\limsupn\,\E\bigg[\sup_{s\in[0,t]}
	\Big|
	\tfrac1n\sum_{k=0}^{\floor{sn}} a_k^n
	\Big(\sum_{l=k}^{b_k^n\wedge\floor{tn}}\phi(\Xnb[l])-\phi(\Xnb[k])\Big) n\Big((\snb[k]-\tnb[k])^2-\E\big[(\snb[k]-\tnb[k])^2\big]\Big)
	\Big|	\bigg]
	\\
	&=\limsupn\,
	\E\bigg[\sup_{s\in[0,t]}
	\Big|
	\tfrac1n\sum_{k=0}^{b^n_\floor{sn}\wedge \floor{tn}}
	\Big(	\phi(\Xnb[k])-\phi(\Xnb[c_k^n]) \Big) n\Big((\snb[k]-\tnb[k])^2-\E\big[(\snb[k]-\tnb[k])^2\big]\Big)
	\Big|	\bigg]
	\\
	&\leq\limsupn\,
	\tfrac1n\sum_{k=1}^{\floor{tn}}
	\Big(	\E\Big[\Big|p_n
	\Big(	\phi(\Xnb[k])-\phi(\Xnb[c_k^n]) \Big)\Big|^{\frac{u}{u-1}}
		\Big]\Big)^{\frac{u-1}u}
	\Big(\E\Big[
	\Big| \tfrac{n}{p_n} \Big((\snb[k]-\tnb[k])^2-\E\big[(\snb[k]-\tnb[k])^2\big]\Big)
	\Big|^u	\Big]\Big)^{\frac{1}u}
	=0\,,
	\end{split}
	\end{align}
	hence we have shown that \eqref{eq:CKL.L2-inverse.Hölder} converges.\\
	For all $n\in\N$ it holds due to independence, \eqref{eq:moments.b_k^n} and \eqref{eq:sigma-tau->alpha} that
	\begin{equation} \label{eq:2.application.l.max.iid.KLM.a}
	\limsupn \tfrac1{\floor{tn}}
	\E\left[b_{0}^n
	\big|n(\snb[0]-\tnb[0])^2-\E\big[n(\snb[0]-\tnb[0])^2\big]
	\big|\right]\hspace{-2pt}
	\leq \hspace{-1pt}
	\limn \tfrac2{\floor{tn}}
	\tfrac{1-p_n}{p_n}
	\E\left[	\big|n(\snb[0]-\tnb[0])^2\big|\right]
	=0.
	\end{equation}
	Furthermore, due to H\"older's inequality, Markov's inequality, independence, Jensen's inequality, \eqref{eq:moments.b_k^n} and \eqref{eq:Lp.cond.sigma.tau} it holds for all  $\eps\in(0,\infty)$ and $v\in(1,2\wedge u)$ that
	\begin{align} \nonumber
	&\limsupn
	\,\E\left[b_{0}^n
	\big|n(\snb[0]-\tnb[0])^2-\E\big[n(\snb[0]-\tnb[0])^2\big]
	\big|
	\1_{\{b_{0}^n
		|n(\snb[0]-\tnb[0])^2-\E[n(\snb[0]-\tnb[0])^2]
		|>\eps\floor{tn}\}}
	\right]
	\\ \nonumber
	\leq\,
	&\limsupn
	\left(\,\E\left[|b_{0}^n
	\big(n(\snb[0]-\tnb[0])^2-\E\big[n(\snb[0]-\tnb[0])^2\big]\big)
	|^v\right]\right)^{\frac{1}{v}}
	(\P(b_{0}^n
	\big|n(\snb[0]-\tnb[0])^2-\E\big[n(\snb[0]-\tnb[0])^2\big]
	\big|>\eps\floor{tn}))^{\frac{v-1}{v}}
	\\
		\begin{split} \label{eq:2.application.l.max.iid.KLM.b}
	\leq\,
	&\limsupn\,
	\E\left[|b_{0}^n	|^v\right]
	\E\left[|	\big(n(\snb[0]-\tnb[0])^2-\E\big[n(\snb[0]-\tnb[0])^2\big]\big)	|^v\right]
	(\eps\floor{tn})^{1-v}
	\\
	\leq\,
	&\limsupn\,
	2^v\big(\E\left[|b_{0}^n	|^2\right]\big)^{\frac v2}
	(p_n)^v
	\big(\E\big[|	\tfrac{n}{p_n}(\snb[0]-\tnb[0])^2	|^u\big]\big)^{\frac vu}
	(\eps\floor{tn})^{1-v}
	\\
	=\,
	&\limn
	2^v\big( \big(2-3p_n+(p_n)^2\big) (p_n)^{-2}\big)^{\frac v2}
	(p_n)^v
	\big(\E\big[|	\tfrac{n}{p_n}(\snb[0]-\tnb[0])^2	|^u\big]\big)^{\frac vu}
	(\eps\floor{tn})^{1-v}
	=0\,.
	\end{split}
	\end{align}
	For all $n\in\N$ and $i\in\N$ let $k_0^n:=0$ and $k_i^n:=\inf\{j\in\N\cap(k_{i-1}^n,\infty)\colon a_j^n=1\}$.
	It turns out that for all $s\in[0,t]$ and $k\in\{0,\dots,\floor{sn}\}\setminus \max_{i\in\N_0}\{k_i^n \colon k_i^n\leq\floor{sn}\}$ it holds that $a_k^n(b_k^n\wedge\floor{tn}-b_k^n\wedge\floor{sn})=0$ and for all $i\in\N_0$, as $b_{k_i^n}^n\geq k_i^n$, that $b_{k_i^n}^n\wedge\floor{tn}-b_{k_i^n}^n\wedge\floor{sn}\leq b_{k_i^n}^n-k_i^n$.
	Hence Lemma~\ref{l:E.max.iid} with $\xi^n=\floor{tn}$, which is applicable due to~\eqref{eq:2.application.l.max.iid.KLM.a} and \eqref{eq:2.application.l.max.iid.KLM.b} to the family
	 $\big((b_{k_i^n}^n-k_i^n)
	\big|n(\snb[k_i^n]-\tnb[k_i^n])^2-\E\big[n(\snb[k_i^n]-\tnb[k_i^n])^2\big]
	\big|\big)_{i\in\N}$ of i.i.d. random variables,
	implies that
	\begin{align}
	\label{eq:3.25b.term3}
	\limsupn\,&\E\bigg[\sup_{s\in[0,t]}
	\Big|
	\tfrac1n\sum_{k=0}^{\floor{sn}} a_k^n
	\Big(\sum_{l=b_k^n\wedge\floor{sn}+1}^{b_k^n\wedge\floor{tn}} \phi(\Xnb[l])\Big) \Big(n(\snb[k]-\tnb[k])^2-\E\big[n(\snb[k]-\tnb[k])^2\big]\Big)
	\Big|	\bigg]
	\\ \nonumber
	\leq\limsupn\,	&
	\|\phi\|_\infty\,
	\E\bigg[\sup_{s\in[0,t]}
	\tfrac1n\sum_{k=0}^{\floor{sn}} a_k^n
	(b_k^n\wedge\floor{tn}-b_k^n\wedge\floor{sn})
	\big|n(\snb[k]-\tnb[k])^2-\E\big[n(\snb[k]-\tnb[k])^2\big]
	\big|	\bigg]
	\\ \nonumber
	\leq\limsupn\,&
	\tfrac{\|\phi\|_\infty\floor{tn}}{n}\,
	\E\bigg[\max_{i\in\{0,\dots,\floor{tn}\}}
	\tfrac1{\floor{tn}}
	(b_{k_i^n}^n-k_i^n)
	\Big|n\big(\snb[k_i^n]-\tnb[k_i^n]\big)^2-\E\Big[n\big(\snb[k_i^n]-\tnb[k_i^n]\big)^2\Big]
	\Big|	\bigg]
	=0\,.
	\end{align}
	which is convergence of \eqref{eq:CKL.L2-inverse.zero-run}.
	Thus we have shown (cf.~\eqref{eq:3.25b.term1} and \eqref{eq:3.25b.term2}) that  \eqref{eq:CKL.L2-inverse}
  converges to zero and consequently that \eqref{eq:3.25a}  holds.
	 This proves \eqref{eq:LLN.for.L2inverse.KLM}.
	 Finally, utilizing \eqref{eq:pf.LLN.KLM.what.remains.when.summing.instead.of.integrating}  with $\phi\equiv1$ and the definition of $(a_k^n)_{n\in\N,k\in\N_0}$ and $(b_k^n)_{n\in\N,k\in\N_0}$, we infer from \eqref{eq:3.25b.term2} with $\phi\equiv1$ that
	 \begin{align} 
	 \begin{split}
	  &\limsupn\,\E\bigg[ \bigg| \int_0^{t}
	 \tfrac{n}{p_n}\Big((\sigma_r^n-\tau_r^n)^2-\E\big[(\sigma_r^n-\tau_r^n)^2\big]\Big)
	 \,dr \bigg|\bigg]
	 \\
	 =\,
	 &\limsupn\,
	 \E\bigg[\Big|\tfrac1n \sum_{k=0}^{\floor{tn}}a_k^n(b_k^n\wedge\floor{tn}-k+1)
	 \Big(\tfrac{n}{p_n}(\snb[k]-\tnb[k])^2-\E\big[\tfrac{n}{p_n}(\snb[k]-\tnb[k])^2\big]\Big)
	 \Big|\bigg]
	 =0\,,
	  \end{split}
	 \end{align}
	 which shows that \eqref{eq:LLN.for.limitting.occupation.measure.KLM} holds as well and finishes the proof of Lemma~\ref{l:LLN.for.limitting.occupation.measure.CKL}.
\end{proof}

We finish this section with the proof of Theorem~\ref{T:KLM}.
\begin{proof}[\textbf{Proof of Theorem \ref{T:KLM}}]
	In order to apply Corollary~\ref{C:stochastic.averaging} we check the assumptions. 
	Fix $u\in(1,\infty)$ for which \eqref{eq:Lp.cond.sigma.tau} holds.
	Without loss of generality assume that Setting~\ref{setting.KLM} holds.
	Set $S:=[0,1]$ and $E:=(-1,\infty)^2$ and equip both of them with the Euclidean distance. 
  For all $n\in\N$ define
	$\T_n:=\frac{\N_0} n$,
	$\theta_n:=\sqrt{np_n}$,
	$g_n:=\big((-1,\infty)^2\ni(\sigma,\tau)\mapsto \frac{n(1-p_n)}{p_n}(\sigma-\tau)^2\in[0,\infty)\big)$
	and let $(L_n,\Dom(L_n))$ and  $((L_{i,n},D_n))_{i\in\{1,2,3\}}$ be the operators defined in
  Lemma~\ref{l:(X,s,t)solves.MP(Ln)}.
  Moreover, define
	$ D_0:=\C^4_b([0,1],\R)$
	and note that $ D_0$ is dense in $\C_b([0,1],\R)$ in the topology of uniform convergence by the Weierstrass approximation theorem.
	Let $A_1, A_2\colon D_0\to \C_b([0,1],\R)$
  be the functions such that for all $f\in D_0$ and $x\in[0,1]$ it holds that
	\begin{equation} \begin{split} \label{eq:def.A_1.A_2.KLM}
	(A_1f)(x)
	&:=x(1-x)\left[\al-\tfrac \gamma 2 +p\beta(\tfrac12-x)\right]f'(x)
	\\&\qquad
	+\tfrac12 \left[x(1-x)+p\beta x^2(1-x)^2\right]f''(x)\,,
	\\
	(A_2f)(x)
	&:=x^2(1-x)^2f''(x)+x(1-x)(1-2x)f'(x)\,.
	\end{split}
	\end{equation}	
	Fix $f\in  D_0$ for the rest of the proof. 
	For all $n\in\N$ define $\pi_n:=\Law{(\snb,\tnb)}$,
	$f_n:=(S\times E\ni(x,\sigma,\tau)\mapsto f(x)\in\R)$ and, with $F$ defined in Lemma~\ref{l:CKL.derivatives},
	$h_n:=\left(S\times E\ni(x,\sigma,\tau)\mapsto \sqrt{\frac n {p_n}} F(x,\sigma,\tau) f'(x)\in\R\right)$.
	We observe that	$f_n\in D_n$,
	Lemma \ref{l:CKL.derivatives} ensures that $h_n \in D_n$ 
	and  that
	$h_n|_{S_n\times E_n}
	=L_{1,n}f_n
	$.

  Assumption \ref{ass:main}.1
  and
  Assumption \ref{ass:main}.4
  are clearly fulfilled. 
  Moreover, 
  Lemma~\ref{l:(X,s,t)solves.MP(Ln)}
  implies that  Assumptions~\ref{ass:main}.2 and \ref{ass:main}.3 hold (with $\tilde S_n:=S$ and $\tilde E_n:=E$).
  By construction (see \eqref{eq:Law.of.sigma,tau}) $(g_n(\sn[t],\tn[t]))_{t\in[0,\infty)}:=(g_n(\snb[\floor{tn}],\tnb[\floor{tn}]))_{t\in[0,\infty)}$ is  stationary,
  due to Assumption~\ref{ass:KLM} it holds that $\sup_{n\in\N}\E[|g_n(\sn[0],\tn[0])|^u]<\infty$ and that 
  $\lim_{n\to\infty}|\E[g_n(\sn[0],\tn[0])]-(1-p)\beta|=0$
  and as a consequence of Lemma~\ref{l:LLN.for.limitting.occupation.measure.CKL} it holds that $\lim_{n\to\infty}\E[|\int_0^t g_n(\sn[s],\tn[s])-\E[g_n(\sn[s],\tn[s])]\,ds|]=0$.
  Therefore $\tilde 5.$~is satisfied as well.

	\textbf{Next we check condition \eqref{eq:hn.converges.to.0}.}
	Using \eqref{eq:E[X-x]<|sigma-tau|/eps} yields for all $n\in\N$, $t\in[0,\infty)$ that
	\begin{equation}\label{eq:checking.fn.to.f}
	\begin{split}
	\E\bigg[\sup_{s\in[0,t]}
	\left|\tfrac{1}{\sqrt{np_n}}h_n(X_s^n,\sigma_s^n,\tau_s^n)\right|\bigg]
	&\leq
	\E\bigg[\sup_{s\in[0,t]}\sup_{x\in S_n}
	\left|\tfrac{1}{p_n}\E^n_{x,\sn[s],\tn[s]}[\Xnb[1]-x]f'(x) \right|\bigg]
		\\
	&\leq
	\|f'\|_\infty
	\tfrac1{np_n}\E\bigg[\max_{k\in\{0,\dots,\floor{tn}\}}
	\left|	n\tfrac{\bar{\sigma}_k^n-\bar{\tau}_k^n}{1+\snb[k]\wedge\tnb[k]}\right|\bigg]\,.
	\end{split}
	\end{equation}
	\noindent
	Note that 
	the number of independent and identically distributed selection regimes throughout generations ${1,\dots,\floor{tn}}$ is binomially distributed with parameters $\floor{tn}$ and~$p_n$. Moreover, from \eqref{eq:|s-t|/(p_n(1+s.wedge.t))} and \eqref{eq:convergence speed|sigma-tau|} condition $i)$ of Lemma~\ref{l:E.max.iid} follows immediately. Hence Lemma~\ref{l:E.max.iid} implies convergence of \eqref{eq:checking.fn.to.f} to zero for $n\to\infty$ and therefore that  \eqref{eq:hn.converges.to.0} holds.\\
	\textbf{Next we check condition \eqref{eq:L0n-E[L0n]}.}
  Applying \eqref{eq:E[(X-x)^2].neu} twice yields for all $t\in[0,\infty)$ that
  \begin{equation}  \begin{split}\label{eq:f''terms}
	&\limsupn\,\E\bigg[   \int_0^{t}\bigg|
	\bigg[
	\tfrac n2 \E^{n,\floor{rn}}_{X_r^n,\sigma_r^n,\tau_r^n} \big[(X^n_{r+\frac1n}-X_r^n)^2\big]
	- \tfrac{X_r^n(1-X_r^n)}{2}-\tfrac{n}{2}(\sigma_r^n-\tau_r^n)^2(X_r^n(1-X_r^n))^2\bigg]f''(X_r^n)
	\\
	&- \int_{E_n} 
	\tfrac n2 \E^{n,\floor{rn}}_{X_r^n,\zeta,\eta} \big[(X^n_{r+\frac1n}-X_r^n)^2\big]
	- \tfrac{X_r^n(1-X_r^n)}{2}-\tfrac{n}{2}(\zeta-\eta)^2(X_r^n(1-X_r^n))^2
	\pi_n (d(\zeta,\eta))f''(X_r^n)
  \bigg|
	\,dr \bigg]
  =0\,.
  \end{split}     \end{equation}
  Then (we use square brackets to gather the arguments per (in-)equality) it holds that 
  [the definition of $L_{0,n}$, Fubini and stationarity of $(\sn[~],\tn[~])$ for every $n\in\N$],
  [  \eqref{eq:f''terms} and \eqref{eq:E[int(X-x)^2f]} applied with $f'''\in\C_b^1([0,1],\R)$]
	 and [Lemma~\ref{l:LLN.for.limitting.occupation.measure.CKL} applied with 
	$\phi:=\left([0,1]\ni x\mapsto \frac12x^2(1-x)^2f''(x)\in\R\right)$]
	yield for all $t\in[0,\infty)$ that
\begin{align} \label{eq:check.L2inverse.sigma.tau}
	\begin{split}
%
%
%
&\color{white}  = \color{black}
	\limsup_{n\to\infty}
	\E\bigg[ \sup_{s\in[0,t]} \bigg| \int_0^{s}
	\bigg[
	\left( L_{0,n}f_n\right)(X_r^n,\sigma_r^n,\tau_r^n)
	-\int_{E_n} \left(L_{0,n}f_n\right)(X_r^n,\zeta,\eta)\pi_n (d(\zeta,\eta))
	\bigg]
	\,dr \bigg|\bigg]
	\end{split}
	\\ \nonumber
	&\leq
	\limsup_{n\to\infty} 
	\,\E\bigg[ \sup_{s\in[0,t]} \bigg| \int_0^{s}
	\tfrac n2\bigg[
	 \E^{n,\floor{rn}}_{X_r^n,\sigma_r^n,\tau_r^n} \big[(X_{r+\frac1n}^n-X_r^n)^2\big]
	- \int_{E_n} 
	 \E^{n,\floor{rn}}_{X_r^n,\zeta,\eta} \big[(X^n_{r+\frac1n}-X_r^n)^2\big]
	\pi_n (d(\zeta,\eta))
  \bigg]
  f''(X_r^n)
	\,dr \bigg|\bigg]\,
	\\ \nonumber
	&\qquad\qquad
  	+2t\limsup_{n\to\infty} 
  	\E\bigg[   \sup_{x\in S_n}\bigg|
  	\bigg[
     n \E^{n}_{x,\sn[0],\tn[0]}
  	\bigg[\int^{\Xnb[1]}_{x} \tfrac12 (\Xnb[1]-v)^2 f'''(v)dv\bigg]	\bigg|
  	 \bigg]
    \\& \nonumber
    =\limsup_{n\to\infty}
	\,\E\bigg[ \sup_{s\in[0,t]} \bigg| \int_0^{s}
	\tfrac{1}{2}(X_r^n(1-X_r^n))^2f''(X_r^n)\Big[n(\sigma_r^n-\tau_r^n)^2-n\int_{E_n}(\zeta-\eta)^2\pi_n(d(\zeta,\eta))\Big]
	\,dr \bigg|\bigg]
	+0
  =0.
	\end{align}
	This proves condition \eqref{eq:L0n-E[L0n]}.
	%
	%
	\textbf{Next we check condition \eqref{eq:generator.A1.stronger}.}
	%
	%
  We apply \eqref{eq:E[X-x]}, \eqref{eq:E[(X-x)^2].neu} and \eqref{eq:E[int(X-x)^2f]} from Lemma \ref{l:Moments.Displacement} together with \eqref{eq:sigma-tau->alpha} to show that
	\begin{align} 
	\begin{split}
	&\color{white}  = \color{black}
	\limsupn\,  \sup_{x\in S_n} \Big| 
	\left({A}_1 f\right)(x) -\int_{E_n} \left(\sqrt{np_n}
	L_{1,n}f_n+L_{0,n}f_n\right)(x,\sigma,\tau)\pi_n (d(\sigma,\tau))
	 \Big|
	\\ 
	&
	\leq
	\limsupn  \sup_{x\in S_n} \Big| 
	x(1-x)\left[\al-\tfrac \gamma 2 +p\beta (\tfrac12-x)\right]f'(x)
	+\tfrac12\left[x(1-x)+p\beta x^2(1-x)^2\right]f''(x)
	\\ 
	&\color{white}  = \color{black}
	-\int_{E_n}  \bigg(
	n\E^n_{x,\sigma,\tau}[	\bar X_1^n-x]f'(x)
	+\tfrac n2\E^n_{x,\sigma,\tau}\big[(\bar X_1^n-x)^2\big]f''(x)
	\\
	&\color{white}  = \color{black}
	+n\E^n_{x,\sigma,\tau}\bigg[\int_x^{\bar X_1^n} \tfrac12 (\bar X_1^n-t)^2f'''(t)dt\bigg]
	\bigg)\pi_n(d(\sigma,\tau))
	\Big|
	=0\,,
	\end{split}
	\end{align}
	which proves \eqref{eq:generator.A1.stronger}.

	\textbf{Next we check
	condition \eqref{eq:generator.A2.stronger}.}
  First we will show that $\frac{1}{\sqrt{np_n}}L_{0,n}h_n$ has no contribution in the limit as $n\to\infty$.
  For all $n\in\N$ it holds that
	\begin{align} \label{eq:L0nL1n.line0}
	&\E\bigg[\sup_{x\in S_n}
	\big|\left(\tfrac{1}{\sqrt{np_n}}L_{0,n}h_n\right)(x,\sn[0],\tn[0])\big|\bigg]
	\\ 
	\label{eq:L0nL1n.line1}
	\leq\,&\E\bigg[\sup_{x\in S_n} \big|
	\tfrac{n(1-p_n)}{2p_n}\E^n_{x,\sn[0],\tn[0]}\big[ (\Xnb[1]-x)^2 \big]
	\tfrac{\partial^2}{\partial x^2}
	\left(F(x,\sn[0],\tn[0])f'(x)\right)\big|\bigg]
	\\
	\label{eq:L0nL1n.line2}
	&+\E\bigg[\sup_{x\in S_n} \big|
	\tfrac{np_n}{2p_n}\E^n_{x,\sn[0],\tn[0]}\big[ (\Xnb[1]-x)^2 \big]
	\int_{E_n} \tfrac{\partial^2}{\partial x^2}
	\left(F(x,\zeta,\eta)f'(x)\right)
	\pi_n(d(\zeta,\eta))\big|\bigg]
	\\
	\label{eq:L0nL1n.line3}
	&+\E\bigg[\sup_{x\in S_n} \Big|
	\tfrac{n(1-p_n)}{2p_n}\E^n_{x,\sn[0],\tn[0]}\Big[ \int_{x}^{\Xnb[1]}(\Xnb[1]-t)^2 
	\tfrac{\partial^3}{\partial t^3}\left(F(t,\sn[0],\tn[0])f'(t)\right) 
	dt\Big]\Big|\bigg]
	\\ 
	\label{eq:L0nL1n.line4}
	&+\E\bigg[\sup_{x\in S_n} \Big|
	\tfrac{np_n}{2p_n}\E^n_{x,\sn[0],\tn[0]}\Big[ \int_{x}^{\Xnb[1]}(\Xnb[1]-t)^2 
	\int_{E_n} \tfrac{\partial^3}{\partial t^3}\left(F(t,\zeta,\eta)f'(t)\right)
	\pi_n(d(\zeta,\eta))dt
	\Big]
	\Big|\bigg]\,.
	\end{align}
	Let $C=\sup_{m\in\{1,2,3\}}\binom{3}{m}\big\|\frac{\partial^{4-m}}{\partial x^{4-m}}f\big\|_{\infty}\in[0,\infty)$. Utilizing $|\Xnb[1]-x|^3\leq|\Xnb[1]-x|^2$ and \eqref{eq:d/dx E[X-x]} yields that 
	\eqref{eq:L0nL1n.line1} as well as \eqref{eq:L0nL1n.line3} are for all $n\in\N$ bounded by
	\begin{align} \label{eq:L0nL1n.lines1&3}
	\begin{split} 
	&
	\tfrac{Cn}{2p_n}\E\bigg[\left(\sup_{x\in S_n} 
	\E^n_{x,\sn[0],\tn[0]}\big[ |\Xnb[1]-x|^2 \big]\right)
	\sum_{m=0}^{3} \sup_{x\in S_n}\big|\tfrac{\partial^m}{\partial x^m}
	F(x,\sn[0],\tn[0]) \big|
	\bigg]
	\\ 	&
	\leq
	\tfrac{8Cn}{p_n}\E\bigg[
	\sum_{m=1}^{3} \sup_{x\in S_n} 
	\E^n_{x,\sn[0],\tn[0]}\big[ (\Xnb[1]-x)^2 \big]
	\big|\tfrac{\sn[0]-\tn[0]}{1+\sn[0]\wedge\tn[0]}\big|^m
	\bigg]\,.
	\end{split} 
	\end{align}
	We estimate \eqref{eq:L0nL1n.lines1&3} using \eqref{eq:nE[(X-x)^2]} for the first two summands and $|\Xnb[1]-x|^2\leq1$ for the third one and apply \eqref{eq:Jensen.for.|s-t|^2/(p_n)}, \eqref{eq:|s-t|/(p_n(1+s.wedge.t))} and \eqref{eq:moments.3.4.sigma.tau} to deduce that
	\begin{equation} \label{L0nL1n.step.2}
	\limsupn\,
	\tfrac{8Cn}{p_n}\E\Big[ 
	\Big(\tfrac1n +\big|\tfrac{\sn[0]-\tn[0]}{1+\sn[0]\wedge\tn[0]}\big|^2\Big)	\Big(\big|\tfrac{\sn[0]-\tn[0]}{1+\sn[0]\wedge\tn[0]}\big|
	+\big|\tfrac{\sn[0]-\tn[0]}{1+\sn[0]\wedge\tn[0]}\big|^2\Big)
	+\big|\tfrac{\sn[0]-\tn[0]}{1+\sn[0]\wedge\tn[0]}\big|^3
	\Big]=0\,,
	\end{equation}
	which shows convergence of \eqref{eq:L0nL1n.line1} and \eqref{eq:L0nL1n.line3}.
	With the same constant $C$ it holds that for all $n\in\N$ we can estimate \eqref{eq:L0nL1n.line2} as well as \eqref{eq:L0nL1n.line4}, using $|\Xnb[1]-x|^3\leq|\Xnb[1]-x|^2\leq1$ together with \eqref{eq:nE[(X-x)^2]} and  \eqref{eq:d/dx E[X-x]}, by
	\begin{align} \begin{split} \label{L0nL1n.lines2+4}
	&
	Cn\E\bigg[\sup_{x\in S_n} 
	\E^n_{x,\sn[0],\tn[0]}\big[ |\Xnb[1]-x|^2 \big]	\bigg]
	\sum_{m=0}^{3} \E\bigg[\sup_{x\in S_n}\big|\tfrac{\partial^m}{\partial x^m}
	F(x,\sn[0],\tn[0]) \big|
	\bigg]
	\\ 
	&\leq
	8Cn	
	\bigg(\E\Big[1\wedge\Big(\tfrac1n +\big|\tfrac{\sn[0]-\tn[0]}{1+\sn[0]\wedge\tn[0]}\big|^2\Big)\Big]\bigg)
	\sum_{m=1}^{3} \E\Big[
	\big|\tfrac{\sn[0]-\tn[0]}{1+\sn[0]\wedge\tn[0]}\big|^m
	\Big]\,.
	\end{split} 
	\end{align}
	Now we utilize the fact that for all real-valued random variables X and all nondecreasing functions $\phi$, $\psi$ it holds that $\E[\phi(X)]\E[\psi(X)]\leq\E[\phi(X)\psi(X)]$ 	
	to infer from \eqref{L0nL1n.lines2+4} and \eqref{L0nL1n.step.2} that \eqref{eq:L0nL1n.line2} and \eqref{eq:L0nL1n.line4} vanish as well in the limit $n\to\infty$.
	Using Fubini and stationarity we conclude for all $t\in[0,\infty)$ that
	\begin{equation}
	\begin{split} \label{eq:lim.L0n.hn.KLM}
	&\limsupn\,\E\bigg[\sup_{s\in[0,t]}
	\big|\int_0^s \left(\tfrac{1}{\sqrt{np_n}}L_{0,n}h_n\right)(X_r^n,\sn[r],\tn[r])dr\big|\bigg]
	\\
	\leq
	t&\lim_{n\to\infty}\E\bigg[\sup_{x\in S_n}
	\big|\left(\tfrac{1}{\sqrt{np_n}}L_{0,n}h_n\right)(x,\sn[0],\tn[0])\big|\bigg]
	=0\,.
	\end{split}
	\end{equation}
	Since it holds for all $n\in\N$, $(x,\sigma,\tau)\in S_n\times E_n$ that   $F(x,\sigma,\tau)=E^n_{x,\sigma,\tau}[\Xnb[1]-x]$
	we infer from \eqref{eq:E[E[X-x]^2-(s-t)^2]}, \eqref{eq:E[X-x]d/dxE[X-x]},	
	Jensen's inequality,  $f\in\C^4_b([0,1],\R)$, \eqref{eq:E[(X-x)^2].neu}, \eqref{eq:sigma-tau->alpha} and \eqref{eq:n^(1/2)d/dx E[X-x]}
	that
	\begin{align}
	\begin{split} \label{eq:A2.sigma.tau}
	&\limsup_{n\to\infty}
	\E\Big[\sup_{x\in S_n} \Big|
	g_n(\sn[0],\tn[0])(A_2f)(x)-\left(L_{1,n}h_n\right)(x,\sn[0],\tn[0])
	\Big|\Big]
	\\
	=\,&
	\limsup_{n\to\infty}
	\E\Big[\sup_{x\in S_n} \Big|
	\tfrac{n(1-p_n)}{p_n}(\sn[0]-\tn[0])^2 \left(x^2(1-x)^2f''(x)+x(1-x)(1-2x)f'(x)\right)
	\\
	&\qquad\qquad-\tfrac n {p_n} \E^n_{x,\sn[0],\tn[0]}[\Xnb[1]-x]
	(1-p_n)\tfrac{\partial}{\partial x}\left(F(x,\sn[0],\tn[0])f'(x)\right)
	\\
	&\qquad\qquad-\tfrac n {p_n} \E^n_{x,\sn[0],\tn[0]}[\Xnb[1]-x]\int_{E_n}p_n \tfrac{\partial}{\partial x}\left(F(x,\zeta,\eta)f'(x)\right)
	\pi_n(d(\zeta,\eta))
	\Big|\Big]
	\end{split}
	\\ \nonumber
	\leq\,&
	\limsup_{n\to\infty}\bigg[
	\E\bigg[\sup_{x\in S_n} \Big|
	\tfrac{n(1-p_n)}{p_n}f''(x) 
	\Big(x^2(1-x)^2(\sn[0]-\tn[0])^2 -\E^n_{x,\sn[0],\tn[0]}[\Xnb[1]-x]	\E^n_{x,\sn[0],\tn[0]}[\Xnb[1]-x]\Big)
	\\ \nonumber
	&\qquad\qquad
	+\tfrac{n(1-p_n)}{p_n}f'(x)\Big(x(1-x)(1-2x)(\sn[0]-\tn[0])^2 -F(x,\sn[0],\tn[0])\tfrac{\partial}{\partial x}F(x,\sn[0],\tn[0])\Big)
	\Big|\bigg] 
	\\\nonumber
	&\qquad\qquad
	+
	\Big(\E\Big[\sup_{x\in S_n}
	\E^n_{x,\sn[0],\tn[0]}\left[ n(\bar X_1^n-x)^2\right]
	\Big]\Big)^{\frac12}
	\\\nonumber
	&\qquad\qquad
	\cdot n^{\frac12}\Big[
	\sup_{x\in S_n}
	\Big|\int_{E_n}
	\left(\tfrac{\partial}{\partial x}F(x,\zeta,\eta)\right)f'(x)
	+F(x,\zeta,\eta)f''(x) 
	\,\pi_n(d(\zeta,\eta))
	\Big|
	\Big]
	\bigg]
	=0\,.
	\end{align}
	Now \eqref{eq:A2.sigma.tau} and \eqref{eq:lim.L0n.hn.KLM} show together with Fubini and stationarity of $(\sn[~],\tn[~])$, $n\in\N$ that \eqref{eq:generator.A2.stronger} holds as well.
	As we have verified all assumptions, Corollary~\ref{C:stochastic.averaging} now implies relative compactness of $(X^n)_{n\in\N}$ and that every limit point $(X_t)_{t\in[0,\infty)}$ 
	%
	 is a solution of the $\D([0,\infty),
	[0,1])$-martingale problem 
	for the pre-generator that is defined for each  
	$\tilde{f}\in D_0$, $x\in[0,1]$ by
	\begin{equation}
	\begin{split}
	\left(A_1\tilde{f}+(1-p)\beta A_2\tilde{f}\right)(x)
	=x(1-x)\left[\al-\tfrac \gamma 2 +p\beta(\tfrac12-x)\right]f'
	+\tfrac12 \left[x(1-x)+p\beta x^2(1-x)^2\right]f''
	\\
	+(1-p)\beta x(1-x)(1-2x)f'
	+(1-p)\beta x^2(1-x)^2f''
	\\
	=
	x(1-x)\left[\al-\tfrac \gamma 2 +(2-p)\beta(\tfrac12-x)\right]f'
	+ \left[x(1-x)+(2-p)\beta x^2(1-x)^2\right]\tfrac{f''}2
	\end{split}
	\end{equation}
	Note that for each limit point $X$
	of the sequence $(X^n)_{n\in\N}$ there exists a modification with continuous sample paths (see, e.g. Prop.~5.3.5 of~\cite{EthierKurtz1986}). 
	In
	particular, we infer from Prop.~4.6 of~\cite{KaratzasShreve1991} existence of a weak solution of the
	SDE~\eqref{eq:SDE.CKL}. In addition, standard Yamada-Watanabe-type
	arguments yield pathwise uniqueness of this SDE; cf., e.g., Theorem 1 in
	\cite{YamadaWatanabe1971}. Therefore, uniqueness
	of a weak solution of the SDE~\eqref{eq:SDE.CKL} follows from a
	Yamada-Watanabe type argument; see, e.g., Proposition~1
	in~\cite{YamadaWatanabe1971}. Finally since any limit point of
	$(X^n)_{n\in\N}$ is a weak solution of the SDE~\eqref{eq:SDE.CKL}, this
	shows that $(X^n)_{n\in\N}$ converges weakly to the unique solution of
	the SDE~\eqref{eq:SDE.CKL}.  
	This finishes the proof of Theorem~\ref{T:KLM}.
\end{proof}

\section*{Acknowledgements}
The authors are deeply indepted to Tom Kurtz for very valuable discussions.
This paper has been partially supported by the DFG Priority Program ``Probabilistic Structures in Evolution'' (SPP 1590), grants HU 1889/4-1 and PF 672/8-1.


\def\cprime{$'$}

\end{document}